\documentclass[12pt,twoside]{amsart}
\usepackage{amssymb}
\usepackage{amscd}
\usepackage{tikz} 
\usepackage[abbrev,alphabetic]{amsrefs}
\usepackage{hyperref}
\usepackage{here}
\usepackage{xypic}
\usepackage{xcolor}
\usepackage{color}

\usepackage[normalem]{ulem} 
\newcommand\redout{\bgroup\markoverwith
{\textcolor{red}{\rule[.5ex]{2pt}{0.4pt}}}\ULon}

\title[Invariants of algebraic varieties over imperfect fields]
{Invariants of algebraic varieties over imperfect fields} 
\author{Hiromu Tanaka} 
\address{Graduate School of Mathematical Sciences, 
The University of Tokyo, 
3-8-1 Komaba, Meguro-ku, Tokyo 153-8914, JAPAN} 
\email{tanaka@ms.u-tokyo.ac.jp}
\newcommand{\rank}[0]{{\operatorname{rank}}}

\newcommand{\red}[0]{{\operatorname{red}}}

\newcommand{\Ker}[0]{{\operatorname{Ker}}}

\newcommand{\Proj}[0]{{\operatorname{Proj}}}
\newcommand{\Spec}[0]{{\operatorname{Spec}}}

\newcommand{\Bs}[0]{{\operatorname{Bs}}}
\newcommand{\Supp}[0]{{\operatorname{Supp}}}

\newcommand{\edim}[0]{{\operatorname{edim}}}
\newcommand{\gedim}[0]{{\operatorname{gedim}}}

\newcommand{\length}[0]{{\operatorname{length}}}
\newtheorem{thm}{Theorem}[section]
\newtheorem{lem}[thm]{Lemma}
\newtheorem{cor}[thm]{Corollary}
\newtheorem{prop}[thm]{Proposition}

\newtheorem{step}{Step}

\theoremstyle{definition}
\newtheorem{ex}[thm]{Example}

\newtheorem{dfn}[thm]{Definition}

\newtheorem{rem}[thm]{Remark}

\makeatletter
  
  \@addtoreset{equation}{thm}
  \makeatother

\newcommand{\MO}{\mathcal{O}}

\newcommand{\Q}{\mathbb{Q}}
\newcommand{\Z}{\mathbb{Z}}
\newcommand{\F}{\mathbb{F}}

\newcommand{\m}{\mathfrak{m}}

\begin{document}

\maketitle

\begin{abstract}
We introduce four invariants of algebraic varieties over imperfect fields, 
each of which measures either geometric non-normality or geometric non-reducedness. 
The first objective of this article is to establish fundamental properties of these invariants. 
We then apply our results to curves over imperfect fields. 
In particular, we establish a genus change formula and prove the boundedness of 
non-smooth regular curves of genus one. 
We also compute our invariants for some explicit examples. 
\end{abstract}

\tableofcontents

\section{Introduction}

Study of fibrations $\pi:V \to W$ is one of main tools in algebraic geometry 
to reduce many problems to the ones of lower dimension. 
To apply this technique, 
it is important 
to derive good properties on fibres from the ambient space $V$. 
For example, we may use the generic smoothness if the base field is of characteristic zero, 
i.e. if $\pi:V \to W$ is a morphism of smooth varieties, then also general fibres of $f$ are smooth. 
In positive characteristic, the same statement is not true in general. 
We naturally encounter such pathological phenomena 
in classification theory of algebraic varieties, 
e.g. quasi-elliptic fibrations \cite{BM76} and wild conic bundles \cite{Kol91}, \cite{MS03}. 
Although general fibres behave wildly in positive characteristic, 
a similar statement holds for generic fibres. 
More specifically, if the ambient space $V$ is a smooth variety, then the generic fibre $V \times_W \Spec\, K(W)$ is a regular variety, i.e. any stalk is a regular local ring. 
On the other hand, the base field $K(W)$ is no longer a perfect field 
even if the original base field is algebraically closed.  
Therefore, it is important to study algebraic varieties over imperfect fields.

The purpose of this article is to give a systematic approach 
to study algebraic varieties over imperfect fields by using invariants. 
Especially, 
we focus on geometric non-normality and geometric non-reducedness. 
Given a field $k$ of characteristic $p>0$ and 
a normal variety  $X$ such that $k$ is algebraically closed in $K(X)$, 
we will introduce the following four invariants: 
\begin{enumerate}
\renewcommand{\labelenumi}{(\roman{enumi})}
\item 
The {\em capacity of denormalising base changes} $\gamma(X/k)$ of $X/k$ 
(Definition \ref{d-gamma}). 
\item 
The {\em Frobenius length of geometric non-normality} $\ell_F(X/k)$ of $X/k$ 
(Definition \ref{d-F-length}). 
\item 
The {\em Frobenius length of geometric non-reducedness}  $m_F(X/k)$ of $X/k$ 
(Definition \ref{d-mF}). 
\item 
The {\em thickening exponent} $\epsilon(X/k)$ of $X/k$ 
(Definition \ref{d-epsilon}). 
\end{enumerate}
All the invariants $\gamma(X/k), \ell_F(X/k), m_F(X/k),$ and $\epsilon(X/k)$ are non-negative integers. 
The first two invariants are related to geometric normality, 
whilst the latter two are related to geometric reducedness. 
Furthermore, the following properties hold.

\begin{enumerate}
\item[(v)] 
The following are equivalent (Proposition \ref{p-char-gnor}). 
\begin{enumerate}
\item $X$ is geometrically normal over $k$. 
\item $\gamma (X/k)=0$. 
\item $\ell_F(X/k)=0$. 
\end{enumerate}
\item[(vi)] The following are equivalent (Proposition \ref{p-char-gred}). 
\begin{enumerate}
\item $X$ is geometrically reduced over $k$. 
\item $\epsilon (X/k)=0$. 
\item $m_F(X/k)=0$. 
\end{enumerate}
\item[(vii)] $m_F(X/k) \leq \ell_F(X/k) \leq \gamma(X/k)$ (Subsection \ref{ss-rel-bet-inv}). 
\item[(viii)] If $[k:k^p] <\infty$, then the following holds (Subsection \ref{ss-rel-bet-inv}): 
\[
\epsilon(X/k) \leq  m_F(X/k) (\log_p [k:k^p]-1).
\]
\end{enumerate}

In order to establish fundamental properties of our invariants, 
we utilise a certain class of purely inserapable base changes, 
which we call elemental extensions (Definition \ref{d-el-ext}). 
This notion is useful for not only general theory 
but also explicit situations (cf. Subsection \ref{ss2-ex}). 
This kind of base changes is originally used 
by Schr\"{o}er (cf. \cite[Lemma 1.3, Proposition 1.5]{Sch10}).

\subsection{Applications}\label{intro-ss-cano}

%The initial motivation of the author 
%was to introduce invariants related to canonical divisors $K_X$. 
\cite{Tan18} and \cite{PW} have found 
relations between canonical divisors of $X$ and $Y$, 
where $Y$ denotes the normalisation of the reduced structure $(X \times_k \overline k)_{\red}$. 
More specifically, \cite{Tan18} has proven that if $X$ is a proper normal variety 
over a field $k$ of characteristic $p>0$, then 
\begin{itemize}
\item 
the linear equivalence $K_Y+C \sim f^*K_X$ holds for some effective Weil divisor $C$ on $Y$, 
and 
\item 
$C$ can be chosen to be nonzero if $X$ is not geometrically normal. 
\end{itemize}
Then \cite{PW} has shown that $C \sim (p-1)D$ for some effective Weil divisor $D$ on $Y$. 
Using these results, 
they succeeded to establish some results related to minimal model program. 
A motivation for introducing 
the invariant $\gamma(X/k)$ is to give a refinement of these results. 

\begin{thm}[cf. Theorem \ref{t-cano-gamma}]\label{intro-cano-gamma}
Let $k$ be a field of characteristic $p>0$ 
and let $X$ be a proper normal variety over $k$ such that $k$ is algebraically closed in $K(X)$. 
For $Y:=(X \times_k k^{1/p^{\infty}})_{\red}^N$, 
let $f:Y \to X$ be the induced morphism. 
If $X$ is not geometrically normal over $k$,  
then there exist nonzero effective Weil divisors $C_1, ..., C_{\gamma(X/k)}$ on $Y$ such that 
the linear equivalence 
\[
K_Y + (p-1)\sum_{i=1}^{\gamma(X/k)} C_i \sim f^*K_X
\]
holds. 
The same statement holds after replacing $\gamma(X/k)$ by $\ell_F(X/k)$ 
(see Proposition \ref{p-Qgor-index}(2)). 
\end{thm}

This result can be considered as a higher dimensional generalisation of Tate's genus change formula (\cite{Sch09}, \cite{Tat52}). 
Tate's genus change formula only treats geometrically integral case. 
Using the  the thickening exponent $\epsilon(X/k)$, 
we formulate a genus change formula for the general case.

\begin{thm}[Theorem \ref{t-genus-change}]\label{intro-genus-change}
Let $k$ be a field of characteristic $p>0$ 
and let $X$ be a regular projective curve over $k$ with $H^0(X, \MO_X)=k$. 
Let $(k', Y)$ be one of 
$(k^{1/p^{\infty}}, (X \times_k k^{1/p^{\infty}})_{\red}^N)$ and 
$(\overline k, (X \times_k \overline k)_{\red}^N)$. 
Set $g_X:=\dim_k H^1(X, \MO_X)$ and $g_Y:=\dim_{k'} H^1(Y, \MO_Y)$. 
If $X$ is not geometrically normal, 
then there exist nonzero effective Cartier divisors $D_1, ..., D_{\gamma(X/k)}$ on $Y$ 
such that  
\[
2g_Y-2+(p-1)\sum_{i=1}^{\gamma(X/k)}  \deg_{k'} D_i = \frac{2g_X-2}{p^{\epsilon(X/k)}}.
\]
\end{thm}

As direct consequences, we obtain criteria for geometric integrity of curves.

\begin{cor}[cf. Corollary \ref{c-genus-change}]\label{intro-genus-change}
Let $k$ be a field of characteristic $p>0$ 
and let $X$ be a regular projective curve over $k$ with $H^0(X, \MO_X)=k$. 
Set $g_X:=\dim_k H^1(X, \MO_X)$. 
Then the following hold. 
\begin{enumerate}
\item 
If $2g_X-2$ is not divisible by $p$, then $X$ is geometrically integral over $k$. 
%\item 
%If $X$ is not smooth, then the inequality
%\[
%2g_X-2  \geq p.
%\]
%holds. 
\item 
If $2g_X-2 < p(p-3)$, then $X$ is geometrically integral over $k$. 
\end{enumerate}
\end{cor}

\begin{rem}
For any prime number $p$, 
there exists a field of characteristic $p>0$ and 
a regular projective curve $X$ over $k$ 
such that $H^0(X, \MO_X)=k$, $X$ is not geometrically reduced over $k$, and 
$2g_X-2=p(p-3)$ (Proposition \ref{p-Fermat}). 
In this sense, the inequality in Corollary \ref{intro-genus-change}(2) 
is optimal. 
\end{rem}

Since it is useful to compare $X$ and $Y=(X \times_k k^{1/p^{\infty}})_{\red}^N$ in \cite{Tan18} and \cite{PW}, 
it is a fundamental problem to derive some good properties on $Y$ from $X$. 
As a result in this direction, 
the Frobenius length of geometric non-normality $\ell_F(X/k)$ 
gives a bound of $\Q$-Gorenstein index of $Y$ as follows.

\begin{thm}[cf. Theorem \ref{t-Qgor-index}]\label{intro-Qgor-index} 
Let $k$ be a field of characteristic $p>0$ and 
let $X$ be a regular variety over $k$ such that $k$ is algebraically closed in $K(X)$. 
Set $\ell:=\ell_F(X/k)$, 
$Y:=(X \times_k k^{1/p^{\infty}})_{\red}^N$, and $Z:=(X \times_k \overline k)_{\red}^N$. 
Then $p^{\ell} K_Y$ and $p^{\ell} K_Z$ are Cartier. 
\end{thm}

%As above, our invariants are used to establish several results related to canonical divisors. 
%Conversely, canonical divisors will be used to define $\gamma(X/k)$ 
%(cf. the proof of Theorem \ref{t-termination}). 

Apart from the results on canonical divisors discussed above, 
the invariants introduced in this paper will be used as follows. 

\begin{enumerate}
\renewcommand{\labelenumi}{(\Roman{enumi})}
\item 
We will apply our results to boundedness problems of varieties over imperfect fields. 
\item 
In a forthcoming article \cite{BT}, Bernasconi and the author give an explicit upper bound 
for torsion indices of relatively torsion line bundles on three-dimensional del Pezzo fibrations $V \to B$. 
\end{enumerate}

(I) 
For instance, if $k$ is a field of characteristic $p >0$ such that $[k:k^p]<\infty$, 
then we prove the boundedness for non-smooth regular curves over $k$ of genus one (Theorem \ref{t-g1-bdd-F-fin}). 
We also prove that if $X$ is a regular curve of genus one 
which is geometrically integral and not smooth, 
then $X$ is a cubic curve when $p=3$ and 
$X$ is a complete intersection of two quadric surfaces when $p=2$ (Theorem \ref{t-gen-g-1}). 
Moreover, in a forthcoming article \cite{Tan2}, the author shows the boundedness of 
geometrically integral regular del Pezzo surfaces. 
This was the original motivation for the author to introduce the invariants. 

(II) 
In \cite{BT}, we apply our results to the generic fibre $X:=V \times_B \Spec\,K(B)$. 
The strategy is to use the following Frobenius factorisation theorem.

\begin{thm}[Theorem \ref{t-Frob-factor1}]\label{intro-F-factor}
Let $k$ be a field of characteristic $p>0$. 
Let $X$ be a normal variety over $k$ such that $k$ is algebraically closed in $K(X)$. 
Set $\ell:=\ell_F(X/k)$. 
Then the $\ell$-th iterated absolute Frobenius morphism of $X \times_k k^{1/p^{\infty}}$ 
factors through the induced morphism 
$g:(X \times_k k^{1/p^{\infty}})_{\red}^N \to X \times_k k^{1/p^{\infty}}$: 
\[
F^{\ell}_{X \times_k k^{1/p^{\infty}}}:X \times_k k^{1/p^{\infty}}\to (X \times_k k^{1/p^{\infty}})_{\red}^N \xrightarrow{g} X \times_k k^{1/p^{\infty}}.
\]
\end{thm}

For example, if $(X \times_k k^{1/p^{\infty}})_{\red}^N \simeq \mathbb P^2_{k^{1/p^{\infty}}}$, then any numerically trivial inversible sheaf $L$ on $X$ 
satisfies $L^{p^{\ell}} \simeq \MO_X$. 
Indeed, this follows from Theorem \ref{intro-F-factor} and the flat base change theorem. 
%For the general case, we apply the same idea together with geometric information. 

\medskip

\textbf{Acknowledgements:} 
The author would like to thank Fabio Bernasconi, Kentaro Mitsui, 
and Jakub Witaszek for useful comments and answering questions. 
The author is also grateful to the referee for reading the paper carefully and for many useful comments. 
The author was funded by 
the Grant-in-Aid for Scientific Research (KAKENHI No. 18K13386).

\section{Preliminaries}

In this section, we summarise notation and fundamental results 
on base fields and base changes. 
Although some of them are well known, 
we give proofs of them for the reader's convenience.

\subsection{Notation}\label{ss-notation}

In this subsection, we summarise notation we will use in this paper. 

\begin{enumerate}
\item We will freely use the notation and terminology in \cite{Har77}. 
In particular, $D_1 \sim D_2$ means linear equivalence of Weil divisors. 
\item For an integral scheme $X$, 
we define the {\em function field} $K(X)$ of $X$ 
as $\MO_{X, \xi}$ for the generic point $\xi$ of $X$. 
For an integral domain $A$, $K(A)$ denotes the function field of $\Spec\,A$. 
\item A scheme $X$ is {\em regular} (resp. {\em normal}) if 
the local ring $\MO_{X, x}$ at any point $x \in X$ is regular (resp. an integrally closed domain). 
In particular, normal schemes are reduced. 
\item 
For a reduced noetherian scheme $X$, its normalisation is denoted by $X^N$. 
Note that if $X$ is an integral scheme, then $X^N$ is the usual normalisation. 
For the general case, if $X=\bigcup_{i \in I} X_i$ is the decomposition into the irreducible components and we equip $X_i$ with the reduced scheme structure, then $X^N = \coprod_{i \in I} X_i^N$. 
\item 
For a scheme $X$, its {\em reduced structure} $X_{\red}$ 
is the reduced closed subscheme of $X$ such that the induced closed immersion 
$X_{\red} \to X$ is surjective. 
\item 
We say that two schemes $X$ and $Y$ over a field $k$ are $k$-{\em isomorphic} 
if there exists an isomorphism $\theta\colon X \to Y$ of schemes 
such that both $\theta$ and $\theta^{-1}$ commutes with the structure morphisms: 
$X \to \Spec\,k$ and $Y \to \Spec\,k$. 
\item For a field $k$, 
we say that $X$ is a {\em variety over} $k$ or a $k$-{\em variety} if 
$X$ is an integral scheme that is separated and of finite type over $k$. 
We say that $X$ is a {\em curve} over $k$ or a $k$-{\em curve} 
if $X$ is a $k$-variety of dimension one. 
\item For a field $k$, let $\overline k$ be an algebraic closure of $k$. 
If $k$ is of characteristic $p>0$ and $e$ is a non-negative integer, 
then we set $k^{p^e}:=\{x^{p^e}\,|\, x \in k\}$ and 
$k^{1/p^e}:=\{x \in \overline k\,|\, k^{p^e} \in k\}$. 
We set $k^{1/p^{\infty}}:=\bigcup_{e=0}^{\infty} k^{1/p^e}
=\bigcup_{e=0}^{\infty} \{x \in \overline k\,|\, x^{p^e} \in k\}$. 
\item 
For a normal variety $X$ over a field $k$, 
we define the {\em canonical divisor} $K_X$ as a Weil divisor on $X$ such that 
$\MO_X(K_X) \simeq \omega_{X/k}$, where $\omega_{X/k}$ denotes the dualising sheaf (cf. \cite[Section 2.3]{Tan18}). 
Canonical divisors are unique up to linear equivalence. 
Note that $\omega_{X/k} \simeq \omega_{X/k'}$ for any 
field extension $k \subset k'$ that induces 
a factorisation $X \to \Spec\,k' \to \Spec\,k$ (\cite[Lemma 2.7]{Tan18}). 
\item 
For an $\F_p$-scheme $X$ and a positive integer $e$, 
we denote by $F^e_X:X \to X$ the {\em absolute Frobenius morphism} of $X$, 
i.e. the homeomorphic affine morphism induced by the $p^e$-th power map: $\MO_X \to \MO_X, a \mapsto a^{p^e}$. 
Set $F_X:=F_X^1$. 
\item 
For a ring $R$ and $R$-algebras $A_1$ and $A_2$, 
the induced ring homomorphism $A_1 \to A_1 \otimes_R A_2,\,\, a \mapsto a \otimes 1$ is 
called the {\em first coprojection}. 
The {\em second coprojection} is defined by the same way. 
\item 
Let $X$ be a regular projective curve over a field $k$. 
We set 
\[
g_X:= \frac{\dim_k H^1(X, \MO_X)}{\dim_k H^0(X, \MO_X)}, 
\]
which is called the {\em genus} of $X$. 
Note that $g_X$ is independent of the choice of $k$, 
i.e. $g_X$ does not change even if we replace $k$ 
by $k'$ for the Stein factorisation $X \to \Spec\,k' \to \Spec\,k$ 
of the structure morphism. 
\item 
For a field $k$ and a $k$-scheme $X$, 
we have the induced morphism $\theta: k \to H^0(X, \MO_X)$. 
By abuse of notation, 
the equation $H^0(X, \MO_X)=k$ means that $\theta$ is an isomorphism. 
\end{enumerate}

\subsection{Base fields and base changes}

In this subsection, we summarise fundamental results on 
base fields and base changes, which are probably well known for experts.

\begin{prop}\label{p-H0=k}
Let $k$ be a field and let $X$ be a proper variety over $k$. 
Then the following hold. 
\begin{enumerate}
\item 
If $k$ is algebraically closed in $K(X)$, 
then $H^0(X, \MO_X)=k$. 
\item 
Assume that $X$ is normal. 
Then $k$ is algebraically closed in $K(X)$ if and only if $H^0(X, \MO_X)=k$. 
\end{enumerate}
\end{prop}

\begin{proof}
The assertion (2) is verified 
by taking the Stein factorisation of the structure morphism $X \to \Spec\,k$. 
The assertion (1) holds, because 
the following composite injection 
\[
k \hookrightarrow H^0(X, \MO_X) \hookrightarrow  H^0(X^N, \MO_{X^N})
\]
is bijective by (2). 
\end{proof}

\begin{prop}\label{p-prelim-girre}
Let $k$ be a field and let $X$ be a variety such that $k$ is algebraically closed in $K(X)$. 
Then $X$ is geometrically irreducible over $k$. 
\end{prop}

\begin{proof}
We may assume that $X$ is a proper normal variety over $k$. 
By Proposition \ref{p-H0=k}, $X$ is geometrically connected over $k$. 
Therefore, the assertion follows from \cite[Lemma 2.2(1)]{Tan18}. 
\end{proof}

\begin{prop}\label{p-bc-comp-field}
Let $k$ be a field. 
Let $X$ be a variety over $k$ such that $k$ is algebraically closed in $K(X)$. 
Let $k \subset l$ be a field extension. 
Then the following hold. 
\begin{enumerate}
\item $K((X \times_k l)_{\red})$ is isomorphic  to the composite field of $K(X)$ and $l$. 
\item If $k \subset l$ is an algebraic extension, then  
$(X \times_k l)_{\red}^N$ is isomorphic to the normalisation 
of $X$ in $K((X \times_k l)_{\red})$. 
\end{enumerate}
\end{prop}

\begin{proof}
We may assume that $X=\Spec\,A$. 
It follows from Proposition \ref{p-prelim-girre} 
that $X$ is geometrically irreducible over $k$. 

Let us show (1). 
Let $L$ be the composite field of $K(X)$ and $l$. 
By the commutative diagram of ring homomorphisms
\[
\begin{CD}
K(X) @>>> L\\
@AAA @AAA\\
k @>>> l, 
\end{CD}
\]
we obtain a ring homomorphism: 
\[
\theta:(K(X) \otimes_k l)_{\red} \to L. 
\]
Since $(K(X) \otimes_k l)_{\red}$ and $L$ are fields, 
$\theta$ is injective. 
We consider $K(X), l$, and $(K(X) \otimes_k l)_{\red}$ 
as subfields of $L$. 
Since $L$ is the composite field of $K(X)$ and $l$, 
it follows from $K(X) \cup l \subset (K(X) \otimes_k l)_{\red}$ 
that $(K(X) \otimes_k l)_{\red} = L$ i.e. $\theta$ is surjective. 
Thus, (1) holds.

Let us show (2). 
Let $Y$ be the normalisation of $X$ in $K((X \times_k l)_{\red})$. 
Then we obtain a finite surjective morphism $Y \to (X \times_k l)_{\red}^N$ 
of normal varieties over $l$. 
It is enough to show that $K((X \times_k l)_{\red}^N)=K(Y)$, 
which follows from the definition of $Y$. 
\end{proof}

\begin{cor}\label{c-bc-comp-field}
Let $k$ be a field. 
Let $X$ be a variety over $k$ such that $k$ is algebraically closed in $K(X)$. 
Let $k \subset l$ be an algebraic field extension. 
Set $X':=(X \times_k l)_{\red}^N$ and 
let $k'$ be the algebraic closure of $l$ in $K(X')$. 
Then the induced morphism 
\[
(X \times_k k')_{\red}^N \to (X \times_k l)_{\red}^N = X'
\]
is an isomorphism. 
\end{cor}

\begin{proof}
By Proposition \ref{p-bc-comp-field}(1), 
we have $K((X \times_k k')_{\red}) = K(X \times_k l)_{\red})$. 
Then the assertion follows from Proposition \ref{p-bc-comp-field}(2). 
\end{proof}

\begin{prop}\label{p-sm-sep}
Let $k$ be a field and 
let $f:X \to Y$ be a dominant $k$-morphism of varieties over $k$ 
such that $k$ is algebraically closed in $K(Y)$ and 
$K(Y) \subset K(X)$ is a (not necessarily algebraic) separable extension. 
Let $k_X$ be the algebraic closure of $k$ in $K(X)$. 
Then $k \subset k_X$ is a finite separable extension. 
\end{prop}

\begin{proof}
Since $K(Y) \subset K(X)$ is a finitely generated separable extension, 
this is a separably generated extension \cite[(27.F)]{Mat80}. 
Therefore, it suffices to treat the following two cases: $K(Y) \subset K(X)$ is purely transcendental extension and $K(Y) \subset K(X)$ is a finite separable extension. 

Assume that $K(Y) \subset K(X)$ is a purely transcendental extension. 
In this case, $K(Y)$ is algebraically closed in $K(X)$. 
Since $k$ is algebraically closed in $K(Y)$, 
the algebraic closure $k_X$ of $k$ in $K(X)$ is equal to $k$. 

It is enough to treat the case when $K(Y) \subset K(X)$ is a finite separable extension. 
Take $\alpha \in k_X$ and let 
\[
f(t)=t^n+a_1t^{n-1}+\cdots+a_n \in K(Y)[t]
\]
be the minimal polynomial of $\alpha$ over $K(Y)$. 
For the minimal polynomial $g(t) \in k[t]$ of $\alpha$ over $k$, 
we have that $g(t)=f(t)h(t)$ for some $h(t) \in K(Y)[t]$. 
In particular, each $a_i$ is algebraic over $k$. 
Since $k$ is algebraically closed in $K(Y)$, we have that $a_i \in k$ for any $i$. 
In other words, we obtain $f(t)=g(t)$. 
Since $K(X)/K(Y)$ is separable, $f(t)$ is a separable polynomial. 
This shows that $k_X$ is separable over $k$. 
\end{proof}

\begin{lem}\label{l-sep-rn-commu}
Let $k \subset l$ be a (not necessarily algebraic) separable field extension and 
let $X$ be a $k$-scheme of finite type over $k$. 
Then the following hold. 
\begin{enumerate}
\item If $X$ is regular, then so is $X \times_k l$. 
\item If $X$ is reduced, then so is $X \times_k l$. 
\item If $X$ is normal, then so is $X \times_k l$. 
\end{enumerate}
\end{lem}

\begin{proof}
In this proof, we use formally smooth morphisms (for its definition, 
we refer to \cite[Ch. IV, D\'efinition 17.1.1]{Gro67}). 
It is sufficient to show (1). 
We may assume that $k$ is of characteristic $p>0$. 
Assume that $X$ is regular. 
Then $X$ is geometrically regular over $\F_p$. 
It follows from \cite[Ch. 0, Proposition 22.6.7(ii)]{EGAIV1} 
that $X$ is formally smooth over $\F_p$. 
Since $l$ is formally smooth over $k$ \cite[Theorem 26.9]{Mat89}, 
it holds that $X \times_k l \to X$ is formally smooth 
\cite[Ch. IV, D\'efinition 17.1.3(iii)]{Gro67}. 
Then $X \times_k l$ is formally smooth over $\F_p$ 
\cite[Ch. IV, D\'efinition 17.1.3(ii)]{Gro67}. 
Hence, $X \times_k l$ is regular \cite[Ch. 0, Proposition 22.6.7(ii)]{EGAIV1}. 
\end{proof}

\begin{prop}\label{p-sep-rn-commu}
Let $k \subset l$ be a (not necessarily algebraic) separable field extension and 
let $X$ be a scheme of finite type over $k$. 
Then the following hold. 
\begin{enumerate}
\item The induced morphism $(X \times_k l)_{\red} \to X_{\red} \times_k l$ is an isomorphism. 
\item The induced morphism $(X \times_k l)^N_{\red} \to X^N_{\red} \times_k l$ is an isomorphism. 
\end{enumerate}
\end{prop}

\begin{proof}
Since $k \subset l$ is separable, 
$X_{\red} \times_k l$ is reduced (Lemma \ref{l-sep-rn-commu}(2)). 
Note that $(X \times_k l)_{\red} \to X \times_k l$ and $X_{\red} \times_k l \to X \times_k l$ 
are homeomorphisms. 
Therefore, also $(X \times_k l)_{\red} \to X_{\red} \times_k l$ is a homeomorphism. 
Then $(X \times_k l)_{\red} \to X_{\red} \times_k l$ 
is a surjective closed immersion of reduced schemes, which is automatically an isomorphism. 
Thus (1) holds. 

Let us show (2). 
Since $k \subset l$ is separable, 
it holds that $X^N_{\red} \times_k l$ is normal (Lemma \ref{l-sep-rn-commu}(3)). 
Then $(X \times_k l)^N_{\red} \to X^N_{\red} \times_k l$ 
is a finite surjective morphism between normal schemes. 
Since it is also birational, (2) holds. 
\end{proof}

\begin{cor}\label{c-sep-rn-commu}
Let $k$ be a field and let $X$ be a variety over $k$ 
such that $k$ is algebraically closed in $K(X)$. 
Assume that $k \subset l$ is  a (not necessarily algebraic) separable field  extension. 
Then the following hold. 
\begin{enumerate}
\item $X \times_k l$ is a variety over $l$. 
\item If $X$ is normal, then also $X \times_k l$ is normal. 
\item $l$ is algebraically closed in $K(X \times_k l)$. 
\end{enumerate}
\end{cor}

\begin{proof}
By Proposition \ref{p-prelim-girre} and Proposition \ref{p-sep-rn-commu}, 
(1) and (2) hold. 
To prove (3), 
we may assume that $X$ is a proper normal variety over $k$. 
Since $H^0(X, \MO_X)=k$ (Proposition \ref{p-H0=k}), 
we obtain $H^0(X \times_k l, \MO_{X \times_k l})=l$. 
Again by Proposition \ref{p-H0=k}, (3) holds. 
\end{proof}

\begin{lem}\label{l-Fbc-alg-cl}
Let $k$ be a field of characteristic $p>0$. 
Let $X$ be a %normal 
variety over $k$ such that $k$ is algebraically closed in $K(X)$. 
Fix a non-negative integer $n$. 
Then $k^{1/p^n}$ is algebraically closed in $K((X \times_k k^{1/p^n})_{\red})$. 
\end{lem}

\begin{proof}
For $F:=K(X)$, we have $K((X \times_k k^{1/p^n})_{\red})=(F \otimes_k k^{1/p^n})_{\red}$. 
Let $l$ be the algebraic closure of $k^{1/p^n}$ in $(F \otimes_k k^{1/p^n})_{\red}$. 
It is enough to show that $l \subset k^{1/p^n}$. 
Fix $x \in l$. 
We have 
\[
x^{p^n} \in l^{p^n} \subset ((F \otimes_k k^{1/p^n})_{\red})^{p^n} \subset F, 
\]
where the last inclusion follows from the fact that 
$(F \otimes_k k^{1/p^n})_{\red}$ coincides with the composite field of $F$ and $k^{1/p^n}$ 
(Proposition \ref{p-bc-comp-field}(1)). 
Since $x^{p^n}$ is algebraic over $k$ and $k$ is algebraically closed in $F=K(X)$, 
we obtain $x^{p^n} \in k$. 
This implies $x \in k^{1/p^n}$.  
\end{proof}

\subsection{Geometric reducedness}

In this subsection, we summarise fundamental results on 
geometric reducedness and geometric normality. 
Most of them are probably well known for experts.

\begin{prop}\label{p-geom-ht1}
Let $k$ be a field of characteristic $p>0$ and let $X$ be a scheme of finite type over $k$. 
Then the following hold. 
\begin{enumerate}
\item The following are equivalent. 
\begin{enumerate}
\item $X$ is geometrically regular over $k$, i.e. smooth over $k$. 
\item $X \times_k k^{1/p}$ is regular. 
\end{enumerate}
\item The following are equivalent. 
\begin{enumerate}
\item $X$ is geometrically reduced over $k$. 
\item $X \times_k k^{1/p}$ is reduced. 
\end{enumerate}
\item The following are equivalent. 
\begin{enumerate}
\item $X$ is geometrically normal over $k$. 
\item $X \times_k k^{1/p}$ is normal. 
\end{enumerate}
\end{enumerate}
\end{prop}

\begin{proof}
The assertion (1) holds by \cite[Chap. 0, Th\'eor\`eme 22.5.8]{EGAIV1}. 
Note that, for any field extension $k \subset k'$, 
$X$ is $S_r$ if and only if so is $X \times_k k'$ \cite[Corollary immediately after Theorem 23.3]{Mat89}. 
Thus, the assertion (2) follows from (1) and 
the fact that being reduced is equivalent to being $S_0$ and $R_1$. 
Thanks to Serre's criterion for normality, we may apply the same argument to (3). 
\end{proof}

\begin{lem}\label{l-geom-red-birat}
Let $k$ be a field of characteristic $p>0$.  
Let $A$ be a $k$-algebra which is an integral domain 
and let $F$ be the function field of $A$. 
Let $k \subset l$ be a field extension. 
Then $A \otimes_k l$ is reduced if and only if $F \otimes_k l$ is reduced. 
\end{lem}

\begin{proof}
Since $l$ is flat over $k$, the injective ring homomorphism $A \hookrightarrow F$ induces 
an injective ring homomorphism $A \otimes_k l \hookrightarrow F \otimes_k l$. 
Therefore, if $F \otimes_k l$ is reduced, then so is $A \otimes_k l$. 

Let us show the inverse implication. 
Assuming that $A \otimes_k l$ is reduced, let us prove that $F \otimes_k l$ is reduced. 
We have the following commutative diagram of injective ring homomorphisms:
\[
\begin{CD}
A @>{\rm injective}>> F\\
@VV{\rm injective}V @VV{\rm injective}V\\
A \otimes_k l @>{\rm injective}>> F \otimes_k l.
\end{CD}
\]
Take an element $\zeta \in F \otimes_k l$ such that $\zeta^m=0$ for some $m \geq 1$. 
There exists $a \in A\setminus \{0\}$ such that $\zeta \cdot (a \otimes 1) \in A \otimes_k l$. 
It holds that 
\[
(\zeta \cdot (a \otimes 1))^m = \zeta^m \cdot (a \otimes 1)^m=0.
\]
Since $A \otimes_k l$ is reduced, we have that $\zeta \cdot (a \otimes 1)=0$. 
The equation $(a \otimes 1)(a^{-1} \otimes 1)=1$ holds in $F \otimes_k l$, 
hence we obtain $\zeta =0$. 
\end{proof}

\begin{prop}\label{p-geom-red-birat}
Let $k$ be a field of characteristic $p>0$.  
Let $X$ be a variety over $k$. 
Let $Y$ be another variety over $k$ which is birational to $X$. 
Then $X$ is geometrically reduced over $k$ if and only if $Y$ is geometrically reduced over $k$
\end{prop}

\begin{proof}
The assertion follows from Lemma \ref{l-geom-red-birat}. 
\end{proof}

Let us show that if $X$ has a rational point around which $X$ is regular, 
then $X$ is reduced (Corollary \ref{c-rat-red}). 
We start with the following auxiliary result.

\begin{prop}\label{p-rat-red}
Let $k$ be a field and let $X$ be a scheme of finite type over $k$. 
Let $P$ be a $k$-rational point of $X$. 
If $\MO_{X, P}$ is a regular local ring, 
then there exists an open neighbourhood $U$ of $P \in X$ such that $U$ is smooth over $k$. 
\end{prop}

\begin{proof}
Set $d:=\dim \MO_{X, x}$. 
We may assume that $X=\Spec\,A$ and that 
there exist elements $x_1, ..., x_d \in A$ such that $\m_P=(x_1, ..., x_d)$, 
where $\m_P$ denotes the maximal ideal corresponding to $P$. 
In particular, $A/\m_P \simeq k$. 
For any field extension $k \subset k'$ and $A':=A \otimes_k k'$, we have that $A'/\m_P A' \simeq k'$. 
This implies that $A'$ is regular at the maximal ideals lying over $\m_P$. 
\end{proof}

\begin{cor}\label{c-rat-red}
Let $k$ be a field and let $X$ be a variety over $k$. 
If there exists a $k$-rational point $P$ of $X$ such that $\MO_{X, P}$ is regular, 
then $X$ is geometrically reduced over $k$. 
\end{cor}

\begin{proof}
By Proposition \ref{p-rat-red}, 
there exists an open neighbourhood $U$ of $P \in X$ such that $U$ is smooth over $k$. 
Hence, $U$ is geometrically reduced over $k$. 
It follows from Proposition \ref{p-geom-red-birat} that also $X$ is geometrically reduced over $k$. 
\end{proof}

\subsection{Pullbacks of Weil divisors}

In this subsection, we recall how to define pullbacks of Weil divisor 
for finite morphisms and morphisms induced by base changes. 

\begin{dfn}\label{d-pull-div}
Let $X$ and $Y$ be integral separated normal excellent schemes and 
let $f:Y \to X$ be a dominant morphism. 
Let $X'$ be the non-regular locus of $X$, which is a closed subset of codimension $\geq 2$. 
Assume that also the codimension of $f^{-1}(X')$ is at least two. 
For a Weil divisor $D$, 
we define its {\em pullback} $f^*D$ as follows. 
First we have a morphism $f_1:Y_1 \to X_1$ for $X_1:=X \setminus X'$ and $Y_1:=Y \setminus f^{-1}(X')$. 
Then the pullback $f_1^*(D|_{X_1})$ is defined, since $f$ is dominant and $D$ is a Cartier divisor. 
Set $f^*(D)$ to be the closure of $f_1^*(D|_{X_1})$. 
Note that if $D_1 \sim D_2$, then $f^*D_1 \sim f^*D_2$. 
\end{dfn}

Let $k$ be a field and let $X$ be a normal variety over $k$ 
such that $k$ is algebraically closed in $K(X)$. 
It holds that $X$ is geometrically irreducible over $k$ 
(Proposition \ref{p-prelim-girre}). 
For a field extension $k \subset k'$, 
the base change $X \times_k k'$ is a separated irreducible scheme of finite type over $k'$. 
Therefore, $(X \times_k k')_{\red}$ is a variety over $k'$ and 
its normalisation $Y$ is a normal variety over $k'$. 
We have a sequence of the induced morphisms 
\[
f:Y \xrightarrow{\nu} (X \times_k k')_{\red} \xrightarrow{\rho} X \times_k k' \xrightarrow{\beta} X.
\] 
For any Weil divisor $D$ on $X$, we define $f^*D$ as in Definition \ref{d-pull-div}. 
In this way, we define a Weil divisor $f^*K_X$ up to linear equivalence.

\section{Elemental extensions}\label{s-el-ext}

The purpose of this section is to establish general theory on elemental extensions. 
In Subsection \ref{s1-el-ext}, we give the definition of elemental extensions 
$(k \subset l \subset k', X' \to X)$. 
In Subsection \ref{s1.1-el-ext}, 
we study some properties on internal degrees $[l:k]$ and external degrees $[k':l]$. 
In Subsection \ref{s1.3-el-ext}, 
we study two typical cases: $[l:k]=p$ and $k'=k^{1/p}$. 
In Subsection \ref{s2-el-ext}, the reduced normalisation of any base change is 
decomposed into elemental extensions and cartesian base changes (Theorem \ref{t-dec-elem}). 
Using this result, we prove, in Subsection \ref{s3-el-ext}, that the monotonicity of canonical divisors 
for reduced normalisation of base changes (Theorem \ref{t-T-PW}). 
In Subsection \ref{s4-el-ext}, we prove the termination of denormalising base changes (Theorem \ref{t-termination}). 
In Subsection \ref{s5-el-ext}, we study behaviour of 
elemental extensions under base changes and morphisms. 

\subsection{Definition and basic properties}\label{s1-el-ext}

In this subsection, we introduce elemental extensions (Definition \ref{d-el-ext}). 
We then establish some fundamental properties 
(Lemma \ref{l-el-ext}, Lemma \ref{l-el-ext-proper}).

\begin{dfn}\label{d-el-ext}
Let $k$ be a field of characteristic $p>0$ 
and let $X$ be a variety over $k$ such that $k$ is algebraically closed in $K(X)$. 
We say that a pair $(k \subset l \subset k', \varphi:X' \to X)$ 
is an {\em elemental extension} of $(k, X)$ if $k \subset l \subset k'$ are field extensions, 
$\varphi:X' \to X$ is a morphism of schemes, and 
the following properties hold. 
\begin{enumerate}
\item $k \subset l$ is a purely inseparable field extension such that $l^p \subset k$ 
and $X \times_k l$ is an integral scheme. 
\item $X' \to X \times_k l$ is the normalisation of the integral scheme $X \times_k l$. 
\item $\varphi:X' \to X$ coincides with the induced morphism. 
\item $k'$ is the algebraic closure of $l$ in $K(X')$. 
\end{enumerate}
In this case, 
by abuse of notation, also the pair $(k', X')$ and the 
following  induced  commutative diagram 
\[
\begin{CD}
X @<\varphi << X'\\
@VV\alpha V @VV\alpha'V\\
\Spec\,k @<\psi<< \Spec\,k'.
\end{CD}
\]
are called an {\em elemental extension} of $(k, X)$. 
%\begin{enumerate}
%\renewcommand{\labelenumi}{(\roman{enumi})}
%\item We call $[l:k]$ the {\em internal degree}.
%\item We call $[k':l]$ the {\em external degree}. 
%\item We call $[k':k]$ the {\em (total) degree}. 
%\end{enumerate}
\end{dfn}

\begin{rem}\label{r-el-ext}
Let $k$ be a field of characteristic $p>0$ 
and let $X$ be a variety over $k$ such that $k$ is algebraically closed in $K(X)$. 
For an elemental extension $(k \subset l \subset k', X' \to X)$ of $(k, X)$, 
it holds by definition that $X'= (X \times_k l)^N$. 
By Corollary \ref{c-bc-comp-field}, we also have that 
the induced morphism $(X \times_k k')_{\red}^N \to (X \times_k l)^N = X'$ 
is an isomorphism. 
\end{rem}

\begin{rem}\label{r-el-ext2}
Let $k$ be a field of characteristic $p>0$ 
and let $X$ be a variety over $k$ such that $k$ is algebraically closed in $K(X)$. 
For an elemental extension $(k \subset l \subset k', X' \to X)$ of $(k, X)$, 
also $(k \subset l \subset k', X' \to X^N)$ is an elemental extension of $(k, X^N)$. 
Indeed, we have the induced morphism $X' \to X^N$ and 
$X^N \times_k l$ is an integral scheme (Lemma \ref{l-geom-red-birat}). 
Then we obtain $X' \simeq (X^N \times_k l)^N$ by Corollary \ref{c-bc-comp-field}. 
\end{rem}

\begin{rem}\label{r-triv-ext}
Let $k$ be a field of characteristic $p>0$ 
and let $X$ be a variety over $k$ such that $k$ is algebraically closed in $K(X)$. 
Let $(k \subset l \subset k', \varphi:X' \to X)$ be an elemental extension 
of $(k, X)$. 
\begin{enumerate}
\item 
If $X$ is proper over $k$, 
then we have the equation $H^0(X, \MO_X)=k$ (Proposition \ref{p-H0=k}(1)) 
and 
the flat base change theorem implies that $H^0(X \times_k l, \MO_{X \times_k l})=l$. 
\item 
If $X$ is proper over $k$ and $X \times_k l$ is normal, 
then (1) and Proposition \ref{p-H0=k}(2) 
imply that $X'=X \times_k k'$, $l=k'$, and $H^0(X', \MO_{X'})=k'$. 
\end{enumerate}
\end{rem}

\begin{lem}\label{l-el-ext}
Let $k$ be a field of characteristic $p>0$ 
and let $X$ be a variety over $k$ such that $k$ is algebraically closed in $K(X)$. 
Let $(k \subset l \subset k', \varphi:X' \to X)$ be an elemental extension 
of $(k, X)$. 
Then the following hold. 
\begin{enumerate}
\item $k'^p \subset k$. 
\item $l \subset k'$ is a finite purely inseparable extension. 
\end{enumerate}
\end{lem}

\begin{proof}
Let us show (1). 
Fix $x \in k'$. 
We have 
\[
x^p \in k'^p \subset K(X')^p=K(X \times_k l)^p =(K(X) \otimes_k l)^p \subset K(X), 
\]
hence $x^p \in K(X)$ and $x^p$ is algebraic over $k$. 
Since  $k$ is algebraically closed in $K(X)$, we obtain $x^p \in k$.  
Thus, (1) holds. 

Let us show (2). 
Since the field extension $l \subset k'$ is algebraic and finitely generated, 
we obtain $[k':l] <\infty$. 
Then (1) implies (2). 
\end{proof}

\begin{lem}\label{l-el-ext-proper}
Let $k$ be a field of characteristic $p>0$ 
and let $X$ be a variety over $k$ such that $k$ is algebraically closed in $K(X)$. 
Let $(k \subset l \subset k', \varphi:X' \to X)$ be an elemental extension of $(k, X)$. 
Then the following are equivalent. 
\begin{enumerate}
\item $X \times_k k'$ is normal. 
\item $X \times_k l$ is normal and $k' = l$. 
\end{enumerate}
Furthermore, if $X$ is proper over $k$, then each of (1) and (2) is equivalent to (3). 
\begin{enumerate}
\item[(3)] $X \times_k l$ is normal. 
\end{enumerate}
\end{lem}

\begin{proof}
Let us show that (2) implies (1). 
Assume (2). 
Then it holds that $X' = X \times_k l = X \times_k k'$. 
Therefore, $X \times_k k'$ is normal. 
Hence (1) holds.

Let us show that (1) implies (2). 
Assume (1). 
Then $X \times_k l$ is normal. 
Therefore, the induced morphism 
\[
X \times_k k' = (X \times_k k')_{\red}^N \to (X \times_k l)^N = X \times_k l
\]
is an isomorphism (Remark \ref{r-el-ext}). 
Since $X \to \Spec\,k$ is faithfully flat, we obtain $k' = l$. 
Hence (1) holds. 
To summarise, we have completed the proof of the equivalence between (1) and (2).

Suppose that $X$ is proper over $k$. 
It is enough to show that (3) implies (2), 
which follows from Remark \ref{r-triv-ext}(2). 
%Assume (3). 
%Then we have $X'=X \times_k l$. 
%Since $H^0(X, \MO_X)=k$ (Proposition \ref{p-H0=k}), 
%the flat base change theorem implies that 
%\[
%H^0(X', \MO_{X'}) 
%= H^0(X \times_k l, \MO_{X \times_k l}) 
%= H^0(X, \MO_X) \otimes_k l 
%= l.
%\]
%This deduces the equation: $k'=l$. 
%Hence (2) holds. 
\end{proof}

\subsection{Internal degrees and external degrees}\label{s1.1-el-ext}

In this subsection, we further introduce terminologies: internal degrees and external degrees. 

\begin{dfn}\label{d-el-ext2}
Let $k$ be a field of characteristic $p>0$ 
and let $X$ be a variety over $k$ such that $k$ is algebraically closed in $K(X)$. 
Let $(k \subset l \subset k', \varphi:X' \to X)$ be an elemental extension of $(k, X)$. 
\begin{enumerate}
\renewcommand{\labelenumi}{(\roman{enumi})}
\item We call $[l:k]$ the {\em internal degree}.
\item We call $[k':l]$ the {\em external degree}. 
\item We call $[k':k]$ the {\em (total) degree}. 
\end{enumerate}
\end{dfn}

The external degree $[k':l]$ is always finite (Lemma \ref{l-el-ext}), 
whilst the internal degree is not necessarily finite (Example \ref{e-external-infty}).

\begin{ex}\label{e-external-infty}
If $k$ is a field of characteristic $p>0$ and $X :=\Spec\,k$, 
then $(k^{1/p}, \Spec\,k^{1/p})$ is an elemental extension of $(k, \Spec\,k)$. 
Indeed, $l:=k^{1/p}$ satisfies the conditions listed in Definition \ref{d-el-ext}. 
In particular, if $[k:k^p]=\infty$, then $[l:k]=[k^{1/p}:k]=\infty$. 
Hence, the internal degree is not necessarily finite. 
\end{ex}

\begin{rem}
Let $k$ be a field of characteristic $p>0$ 
and let $X$ be a variety over $k$ such that $k$ is algebraically closed in $K(X)$. 
Let $(k', X')$ be an elemental extension of $(k, X)$. 
Although the choice of $l$ is not necessarily unique  (cf. Proposition \ref{p-Fermat}), 
we shall later prove that the internal degree $[l:k]$ and the external degree $[k':l]$ 
do not depend on the choice of $l$ (Proposition \ref{p-epsilon-compute}). 
We do not use this fact in this paper. 
\end{rem}

\subsection{Examples}\label{s1.3-el-ext}

We shall frequently use two typical elemental extensions: 
$[l:k]=p$ (Example \ref{e-int-p}) and $k'=k^{1/p}$ (Theorem \ref{t-schroer}). 
Both of them are essentially obtained by Schr\"{o}er. 

\begin{ex}\label{e-int-p}
Let $k$ be a field of characteristic $p>0$ 
and let $X$ be a variety over $k$ such that $k$ is algebraically closed in $K(X)$. 
Let $k \subset k'$ be a purely inseparable extension of degree $p$. 
Then $X \times_k l$ is automatically an integral scheme by \cite[Lemma 3.3]{Tan18}. 
Let $X'$ be the normalisation of $X \times_k l$ and 
let $k'$ be the algebraic closure of $l$ in $K(X')$. 
Then $(k \subset l \subset k', X' \to X)$ is an elemental extension of $(k, X)$ 
of internal degree $p$. 
\end{ex}

\begin{thm}\label{t-schroer}
Let $k$ be a field of characteristic $p>0$ and 
let $X$ be a variety over $k$ such that $k$ is algebraically closed in $K(X)$. 
Then the following hold. 
\begin{enumerate}
\item 
There exists an intermediate field $k \subset l \subset k^{1/p}$ such that 
$X \times_k l$ is an integral scheme and the algebraic closure of $l$ in $K(X \times_k l)$ 
is equal to $k^{1/p}$. 
In particular, if $X' := (X \times_k k^{1/p})_{\red}^N$, 
then $(k \subset l \subset k^{1/p}, \varphi:X' \to X)$ is an elemental extension of $(k, X)$ 
for the induced morphism $\varphi$ (cf. Remark \ref{r-el-ext}).
\item 
Fix an intermediate field $k \subset l \subset k^{1/p}$ such that 
$X \times_k l$ is an integral scheme and the algebraic closure of $l$ in $K(X \times_k l)$ 
is equal to $k^{1/p}$. 
Then the following are equivalent. 
\begin{enumerate}
\item 
$X$ is geometrically reduced over $k$. 
\item 
$l=k^{1/p}$. 
\end{enumerate}
\end{enumerate}
\end{thm}

\begin{proof}
Let us show (1). 
We may assume that $X$ is normal. 
If $X$ is proper and $H^0(X, \MO_X)=k$, then 
the assertion follows from \cite[Proposition 1.5]{Sch10}. 
The general case is reduced to this case 
by taking an open immersion $j:X \hookrightarrow \overline X$ to a proper normal variety $\overline X$. %replacing $X$ with a proper normal variety birational to $X$. 
Thus (1) holds. 

Let us show (2). 
If (b) holds, then $X \times_k k^{1/p}$ is reduced, 
hence Proposition \ref{p-geom-ht1}(2) implies that $X$ is geometrically reduced over $k$. 
Thus (b) implies (a). 
Conversely, assume (a). 
Take a non-empty affine open subset $Y$  of $X$ such that $Y \times_k l$ is normal. 
In particular, the structure morphism $Y \times_k l \to \Spec\,l$ factors through $\Spec\,k^{1/p}$. 
Hence, there is an injective ring homomorphism 
\[
k^{1/p} \hookrightarrow H^0(Y \times_k l, \MO_{Y \times_k l}),
\]
which induces another injective ring homomorphism:
\[
k^{1/p} \otimes_l k^{1/p} \hookrightarrow H^0(Y \times_k k^{1/p}, \MO_{Y \times_k k^{1/p}}).
\]
Since $H^0(Y \times_k k^{1/p}, \MO_{Y \times_k k^{1/p}})$ is reduced, 
so is $k^{1/p} \otimes_l k^{1/p}$. 
Then (b) holds. 
\end{proof}

\subsection{Decomposition into elemental extensions}\label{s2-el-ext}

The purpose of this subsection is to prove the decomposition theorem (Theorem \ref{t-dec-elem}).
As a consequence, 
an arbitrary variety reaches to a geometrically normal variety after taking finitely many elemental extensions (Corollary \ref{c-dec-elem2}). 
Finally, we give a criterion of geometric normality and geometric reducedness 
by using decomposition as in Theorem \ref{t-dec-elem} (Theorem \ref{t-dec-elem3}). 
We start with the following auxiliary result.

\begin{lem}\label{l-N-descend}
Let $k$ be a field of characteristic $p>0$ 
and let $X$ be a variety over $k$ such that $k$ is algebraically closed in $K(X)$. 
Let $k \subset l$ be a field extension. 
Then there exists an intermediate field $k \subset k' \subset l$ such that 
\begin{enumerate}
\item $k'$ is a finitely generated field over $k$ 
(in particular, if $k \subset l$ is algebraic, then we have $[k':k]<\infty$), and 
\item for an arbitrary intermediate field $k' \subset k'' \subset l$ and the induced diagram 
\[
\begin{CD}
(X \times_k k'')_{\red}^N @<<< (X \times_k l)_{\red}^N\\
@VVV @VVV\\
(X \times_k k'')_{\red} @<<< (X \times_k l)_{\red}\\
@VVV @VVV\\
X \times_k k'' @<<< X \times_k l\\
@VVV @VVV\\
\Spec\,k'' @<<< \Spec\,l,\\
\end{CD}
\]
each square is cartesian. 
\end{enumerate}
\end{lem}

\begin{proof}
We can find an intermediate field $k \subset k' \subset l$, schemes $Y$ and $Z$ of finite type over $k'$, and a commutative diagram 
\[
\begin{CD}
Z @<<< (X \times_k l)_{\red}^N\\
@VVV @VVV\\
Y @<<< (X \times_k l)_{\red}\\
@VVV @VVV\\
X \times_k k' @<<< X \times_k l\\
@VVV @VVV\\
\Spec\,k' @<<< \Spec\,l\\
\end{CD}
\]
such that $k'$ is a finitely generated field over $k$ and each square is cartesian. 
By faithfully flat descent, we can show that 
the induced morphisms $Y \to (X \times_k l)_{\red}$ and 
$Z \to (X \times_k l)_{\red}^N$ are isomorphisms. 
Hence (1) holds. 
It is clear that (2) holds. 
\end{proof}

\begin{thm}\label{t-dec-elem}
Let 
\[
\begin{CD}
X @<\varphi << Y\\
@VV\alpha V @VV\beta V\\
\Spec\,k @<\psi << \Spec\,l
\end{CD}
\]
be a commutative diagram of schemes such that 
\begin{enumerate}
\item 
$\psi$ is induced by a purely inseparable field extension $k \subset l$, 
\item 
$X$ is a normal variety over $k$ such that $k$ is algebraically closed in $K(X)$, 
\item 
$Y$ is a normal variety over $l$ such that $l$ is algebraically closed in $K(Y)$, 
\item 
$Y=(X \times_k l)_{\red}^N$, and $\varphi$ is the induced morphism. 
\end{enumerate}
Set $X_0:=X$ and $k_0:=k$. 
Then there exists a commutative diagram of schemes 
\[
\begin{CD}
X_0 @<\varphi_1 << X_1 @<\varphi_2 << \cdots @<\varphi_n << X_n @<\widetilde{\varphi} << Y\\
@VV\alpha_0 V @VV\alpha_1 V @. @VV\alpha_n V @VV\beta V\\
\Spec\,k_0 @<\psi_1<< \Spec\,k_1 @<\psi_2<< \cdots @<\psi_n<< \Spec\,k_n @<\widetilde{\psi} <<  \Spec\,l
\end{CD}
\] 
such that the square 
\[
\begin{CD}
X_{i-1} @<\varphi_i <<  X_i\\
@VV\alpha_{i-1}V @VV\alpha_iV\\
\Spec\,k_{i-1} @<\psi_i<< \Spec\,k_i
\end{CD}
\]
is an elemental extension of internal degree $p$ 
(cf. Definition \ref{d-el-ext} and Definition \ref{d-el-ext2}) for any $i \in \{1, ..., n\}$ and 
the rightmost square 
\[
\begin{CD}
X_n @<\widetilde{\varphi} << Y\\
@VV\alpha_n V @VV\beta V\\
\Spec\,k_n @<\widetilde{\psi} <<  \Spec\,l
\end{CD}
\]
is cartesian. 
\end{thm}

\begin{proof}
For an open immersion $j:X \hookrightarrow \overline X$ to a proper normal variety over $k$, 
we may replace $X$ by $\overline X$ in order to prove the assertion. 
Therefore, the problem is reduced to the case when $X$ is proper over $k$.

By Lemma \ref{l-N-descend}, 
there exists an intermediate field $k \subset \widetilde k \subset l$ such that 
\begin{enumerate}
\item[(a)] $[\widetilde k:k] < \infty$, and 
\item[(b)] for the induced diagram 
\[
\begin{CD}
Z:=(X \times_k \widetilde k)_{\red}^N @<<< (X \times_k l)_{\red}^N=Y\\
@VVV @VVV\\
(X \times_k \widetilde k)_{\red} @<<< (X \times_k l)_{\red}\\
@VVV @VVV\\
X \times_k \widetilde k @<<< X \times_k l\\
@VVV @VVV\\
\Spec\,\widetilde k @<<< \Spec\,l,\\
\end{CD}
\]
each square is cartesian. 
\end{enumerate}
By (3) and Proposition \ref{p-H0=k}, we have $H^0(Y, \MO_Y)=l$. 
Then it holds that 
\[
H^0(Z, \MO_Z) \otimes_{\widetilde k} l = 
H^0(Y, \MO_Y)=l.
\]
Therefore, we deduce that $H^0(Z, \MO_Z)=\widetilde{k}$.

%Set $\kappa:=H^0(Z, \MO_Z)$. 
%Then it follows from the flat base change theorem that 
%\[
%\kappa \otimes_{\widetilde k} l =H^0(Z, \MO_Z) \otimes_{\widetilde k} l = H^0(Y, \MO_Y)=:\lambda.
%\]
%In particular, we have that $Z \times_{\kappa} \lambda \simeq Y$. 

If $k_0=\widetilde{k}$, then the assertion is clear. 
Assume that $k_0 \subsetneq \widetilde{k}$. 
We now construct a new triple $(l_1, k_1, X_1)$ as follows. 
Take an intermediate field $k_0 \subset l_1 \subset \widetilde{k}$ such that $[l_1:k_0]=p$. 
Then the scheme $X \times_k l_1 = X_0 \times_{k_0} l_1$ is a variety over $l_1$ \cite[Lemma 3.3]{Tan18}. 
Set $X_1$ to be the normalisation of $X_1 \times_{k_0} l_1$. 
Finally, we define $k_1$ as the algebraic closure of $l_1$ in $K(X_1)$. 
Then $Z \to X_0$ factors through the induced morphism $X_1 \to X_0$: 
\[
X_0 \leftarrow X_1 \leftarrow Z. 
\]
In particular, we have field extensions $K(X_0) \subset K(X_1) \subset K(Z)$. 
Then 
we obtain $k_1 \subset \widetilde{k}$, as $\widetilde{k}$ is algebraically closed in $K(Z)$. 
Repeating this procedure, we obtain a commutative diagram of schemes 
\[
\begin{CD}
X_0 @<\varphi_1 << X_1 @<\varphi_2 << \cdots @<\varphi_n << X_n @<<< Z\\
@VV\alpha_0 V @VV\alpha_1 V @. @VV\alpha_n V @VVV\\
\Spec\,k_0 @<\psi_1<< \Spec\,k_1 @<\psi_2<< \cdots @<\psi_n<< \Spec\,k_n @= \Spec\,\widetilde{k}
\end{CD}
\] 
such that %$k_n=\widetilde{k}$ and 
the diagram 
\[
\begin{CD}
X_{i-1} @<\varphi_i <<  X_i\\
@VV\alpha_{i-1}V @VV\alpha_iV\\
\Spec\,k_{i-1} @<\psi_i<< \Spec\,k_i
\end{CD}
\]
is an elemental extension of internal degree $p$ for any $i \in \{1, ..., n\}$.

In order to prove the remaining assertion, it is enough to show that $Z \simeq X_n$. 
Since both $X_n \to X_0$ and $Z \to X_0$ are finite surjective morphisms, 
it is enough to show that $K(X_n) \supset K(Z)$. 
Since $K(Z)$ is the composite field of $K(X)$ and $\widetilde{k}$ (Proposition \ref{p-bc-comp-field}), 
it suffices to prove that $\widetilde{k} \subset K(X_n)$. 
This follows from $\widetilde{k} = k_n \subset K(X_n)$. 
\end{proof}

\begin{cor}\label{c-dec-elem2}
Let $k_0$ be a field of characteristic $p>0$ and 
let $X_0$ be a normal variety over $k_0$ such that $k_0$ is algebraically closed in $K(X_0)$. 
%Set $Y:=(X \times_k k^{1/p^{\infty}})_{\red}^N$. 
Then there exists a commutative diagram of schemes 
\[
\begin{CD}
X_0 @<\varphi_1 << X_1 @<\varphi_2 << \cdots @<\varphi_n << X_n\\% @<\widetilde{\varphi} << Y\\
@VV\alpha_0 V @VV\alpha_1 V @. @VV\alpha_n V \\ %@VV\beta V\\
\Spec\,k_0 @<\psi_1<< \Spec\,k_1 @<\psi_2<< \cdots @<\psi_n<< \Spec\,k_n \\%@<\widetilde{\psi} <<  \Spec\,k^{1/p^{\infty}}
\end{CD}
\] 
such that the diagram 
\[
\begin{CD}
X_{i-1} @<\varphi_i <<  X_i\\
@VV\alpha_{i-1}V @VV\alpha_iV\\
\Spec\,k_{i-1} @<\psi_i<< \Spec\,k_i
\end{CD}
\]
is an elemental extension of internal degree $p$ 
(cf. Definition \ref{d-el-ext} and Definition \ref{d-el-ext2}) for any $i \in \{1, ..., n\}$ and 
$X_n$ is geometrically normal over $k_n$. 
\end{cor}

\begin{proof}
The assertion (1) holds 
by applying Theorem \ref{t-dec-elem} for $l:=k^{1/p^{\infty}}$ and 
$Y:=(X \times_k k^{1/p^{\infty}})_{\red}^N$. 
\end{proof}

\begin{thm}\label{t-dec-elem3}
Let $k_0$ be a field of characteristic $p>0$ and 
let $X_0$ be a normal variety over $k_0$ such that $k_0$ is algebraically closed in $K(X_0)$. 
Let $(k_{i-1} \subset l_i \subset k_i, \varphi_i:X_i \to X_{i-1})$ 
be an elemental extension (cf. Definition \ref{d-el-ext}) for any $i \in \{1, ..., n\}$. 
Let 
\[
\begin{CD}
X_0 @<\varphi_1 << X_1 @<\varphi_2 << \cdots @<\varphi_n << X_n\\
@VV\alpha_0 V @VV\alpha_1 V @. @VV\alpha_n V \\ 
\Spec\,k_0 @<\psi_1<< \Spec\,k_1 @<\psi_2<< \cdots @<\psi_n<< \Spec\,k_n \\
\end{CD}
\] 
be the induced commutative diagram. 
Then the following hold. 
\begin{enumerate}
\item 
Assume that $X_n$ is geometrically normal over $k_n$. 
Then the following are equivalent. 
\begin{enumerate}
\item 
$X_0$ is geometrically normal over $k_0$. 
\item 
$X_{i-1} \times_{k_{i-1}} l_i$ is normal and $l_i=k_i$ for any $i \in \{1, ..., n\}$. 
\end{enumerate}
Furthermore, if $X_0$ is proper over $k_0$, then each of (a) and (b) is equivalent to (c). 
\begin{enumerate}
\item[(c)]  
$X_{i-1} \times_{k_{i-1}} l_i$ is normal for any $i \in \{1, ..., n\}$. 
\end{enumerate}
\item 
Assume that $X_n$ is geometrically reduced over $k_n$. 
Then the following are equivalent. 
\begin{enumerate}
\item[(a)'] 
$X_0$ is geometrically reduced over $k_0$. 
\item[(b)']  
$l_i=k_i$ for any $i \in \{1, ..., n\}$. 
\end{enumerate}
\end{enumerate}
\end{thm}

\begin{proof}
Let us show (1). 
We first prove that (a) implies (b). 
Assume (a). 
%If $X_0$ is geometrically normal over $k_0$, 
Then it is clear that $X_0 \times_{k_0} k_1$ is normal. 
It follows from Lemma \ref{l-el-ext-proper} that 
$X_0 \times_{k_0} l_1$ is normal and $l_1=k_1$. 
Hence, it holds that $X_1=X_0 \times_{k_0} k_1$, which is again geometrically normal over $k_1$. 
Repeating this argument, we see that (b) holds. 

Let us prove that (b) implies (a). 
Assume (b). 
Then we have that $X_1=X_0 \times_{k_0} l_1=X_0 \times_{k_0} k_1$. 
Repeating the same argument, we have that $X_n=X_0 \times_{k_0} k_n$. 
Since $X_n$ is geometrically normal over $k_n$, $X_0$ is geometrically normal over $k_0$. 
Therefore, (a) holds. 
Finally, if $X_0$ is proper over $k_0$, 
it follows from Lemma \ref{l-el-ext-proper} that (b) and (c) are equivalent. 
This completes the proof of (1). 

The assertion (2) follows from (1). 
Indeed, after replacing $X_0$ by a suitable non-empty open subset of $X_0$, 
(a) and (a)' are equivalent, and so are (b) and (b)' (cf. Proposition \ref{p-geom-red-birat}). 
\end{proof}

\subsection{Canonical divisors and base changes}\label{s3-el-ext}

The purpose of this subsection is to prove the monotonicity of canonical divisors 
under base changes (Theorem \ref{t-T-PW}).  
The result is known when $X$ is proper and $l=k^{1/p^{\infty}}$ ((\cite[Theorem 1.1]{PW}, \cite[Theorem 4.2]{Tan18}). 
Indeed, our strategy is very similar to the one in \cite[Theorem 4.2]{Tan18} and depends on \cite{PW}. 
We include its proof for the sake of completeness.

\begin{thm}\label{t-T-PW}
Let $k$ be a field of characteristic $p>0$ and 
let $X$ be a normal variety over $k$ such that $k$ is algebraically closed in $K(X)$. 
Let $k \subset l$ be a field extension. 
Set $Y:=(X \times_k l)_{\red}^N$. 
Then there exists an effective Weil divisor $C$ on $Y$ such that 
the linear equivalence 
\begin{equation}\label{e1-T-PW}
K_Y + (p-1) C \sim f^*K_X
\end{equation}
holds, where $f:Y \to X$ denotes the induced morphism.  
Furthermore, $C$ can be chosen to be nonzero if one of the following hold. 
\begin{enumerate}
\item $X \times_k l$ is reduced and non-normal. 
\item $X$ is proper over $k$ and $X \times_k l$ is not normal. 
\end{enumerate}
\end{thm}

\begin{proof}
The proof consists of three steps.

\begin{step}\label{s1-T-PW}
In order to prove the assertion of Theorem \ref{t-T-PW}, 
we may assume that $k \subset l$ is a purely inseparable extension. 
\end{step}

\begin{proof} (of Step \ref{s1-T-PW})
There is the decomposition $k \subset k' \subset k'' \subset l$ such that 
$k \subset k'$ is a purely transcendental extension, 
$k' \subset k''$ is a separable algebraic extension, 
and $k'' \subset l$ is a purely inseparable extension. 
Then $X'':=X \times_k k''$ is a normal variety over $k''$
such that $k''$ is algebraically closed in $K(X'')$ 
(cf. Proposition \ref{p-H0=k}). 
Replacing $(k'', X'')$ by $(k, X)$, 
the problem is reduced to the case when $k \subset l$ is a purely inseparable extension. 
This completes the proof of Step \ref{s1-T-PW}. 
\end{proof}

\begin{step}\label{s2-T-PW}
If (1) holds, 
then there exists a nonzero effective Weil divisor $C$ on $Y$ such that 
the linear equivalence (\ref{e1-T-PW}) holds. 
\end{step}

\begin{proof}(of Step \ref{s2-T-PW})
Fix an open immersion $X \hookrightarrow \overline X$ to a proper normal variety $\overline X$. 
It follows from \cite[Theorem 1.1 and Theorem 1.2]{PW} that 
there is an effective Weil divisor $\overline C$ on $\overline Y$ 
such that $K_{\overline Y} +(p-1)\overline C \sim \overline{f}^*K_{\overline X}$ 
and $\Supp\,\overline C$ set-theoretically equal to the closed subscheme 
defined by the conductor ideal. 
In particular, $\overline C \cap Y \neq \emptyset$. 
This completes the proof of Step \ref{s2-T-PW}. 
\end{proof}

\begin{step}\label{s3-T-PW}
Assume that $k \subset l$ is a purely inseparable extension. 
Then there exists an effective Weil divisor $C$ on $Y$ such that 
\begin{enumerate}
\item[(I)] 
the linear equivalence 
\[
K_Y + (p-1) C \sim f^*K_X
\]
holds, and 
\item[(II)] 
if $X$ is proper and $X \times_k l$ is not normal, then we have $C \neq 0$ . 
\end{enumerate}
\end{step}

\begin{proof}(of Step \ref{s3-T-PW}) 
Set $k_0:=k$ and $X_0:=X$. 
Let 
\[
\begin{CD}
X_0 @<\varphi_1 << X_1 @<\varphi_2 << \cdots @<\varphi_n << X_n @<\widetilde{\varphi} << Y\\
@VV\alpha_0 V @VV\alpha_1 V @. @VV\alpha_n V @VV\beta V\\
\Spec\,k_0 @<\psi_1<< \Spec\,k_1 @<\psi_2<< \cdots @<\psi_n<< \Spec\,k_n @<\widetilde{\psi} <<  \Spec\,l
\end{CD}
\] 
be a commutative diagram satisfying the properties in Theorem \ref{t-dec-elem}. 
For any $i$, 
there exists an intermediate field $k_{i-1} \subset l_i \subset k_i$ 
such that $(k_{i-1} \subset l_i \subset k_i, \varphi_i:X_i \to X_{i-1})$ is an elemental extension. 
In particular, $X_{i-1} \times_{k_{i-1}} l_i$ is an integral scheme and 
$\varphi_i$ is decomposed into two morphisms: 
\[
\varphi_i:X_i \xrightarrow{\nu_i} X_{i-1} \times_{k_{i-1}} l_i \xrightarrow{\beta_i} X_{i-1},
\]
where $\beta_i$ is the projection and $\nu_i$ is the normalisation of $X_{i-1} \times_{k_{i-1}} l_i$. 
In particular, it follows from \cite[Theorem 1.1 and Theorem 1.2]{PW} that 
the linear equivalence 
\begin{equation}\label{e2-T-PW}
K_{X_i} + (p-1)C_i \sim \varphi_i^* K_{X_{i-1}}
\end{equation}
holds for some effective Weil divisor $C_i$ on $X_i$ such that $\Supp\,C_i$ 
is set-theoretically equal to the conductor of $\nu_i$. 
In particular, $C_i=0$ if and only if $X_{i-1} \times_{k_{i-1}} l_i$ is normal. 
For the induced morphism $h_i:Y \to X_i$, (\ref{e2-T-PW}) implies that 
\[
K_Y+(p-1) \sum_{i=1}^n h_i^*C_i\sim f^*K_X. 
\]
Set $C:=\sum_{i=1}^n h_i^*C_i$. 
Then (I) holds. 

Let us show (II). 
%If $X \times_k l$ is normal, then it is clear that $C=0$. 
Assuming that $C=0$ and $X$ is proper, 
it suffices to prove that $X \times_k l$ is normal. 
Since $C=\sum_{i=1}^n h_i^*C_i=0$, we have $C_i=0$ for any $i$. 
Then $X_{i-1} \times_{k_{i-1}} l_i$ is normal for any $i$. 
It follows from Lemma \ref{l-el-ext-proper} that $l_i=k_i$. 
To summarise, we have that 
\[
X_i = X_{i-1} \times_{k_{i-1}} l_i = X_{i-1} \times_{k_{i-1}} k_i
\]
for any $i$. 
By induction on $i$, it holds that 
\[
X \times_k l = X_0 \times_{k_0} l = (X_0 \times_{k_0} k_n) \times_{k_n} l = X_n\times_{k_n} l = Y.
\]
Therefore, $X \times_k l$ is normal. 
Hence, (II) holds. 
This completes the proof of Step \ref{s3-T-PW}. 
\end{proof}
Step \ref{s1-T-PW}, Step \ref{s2-T-PW}, and 
Step \ref{s3-T-PW} complete the proof of Theorem \ref{t-T-PW}. 
\end{proof}

\subsection{Termination of denormalising base changes}\label{s4-el-ext}

The purpose of this subsection is to show termination of sequences of denormalising base changes (Theorem \ref{t-termination}).

\begin{thm}\label{t-termination}
Let $k_0$ be a field of characteristic $p>0$ and 
let $X_0$ be a normal variety over $k_0$ such that $k_0$ is algebraically closed in $K(X_0)$. 
Then there exists a positive integer $\gamma_0$ that satisfies the following property: 
if 
\[
\begin{CD}
X_0 @<\varphi_1 << X_1 @<\varphi_2 << \cdots @<\varphi_n << X_n @<\varphi_{n+1}<< \cdots\\
@VV\alpha_0 V @VV\alpha_1 V @. @VV\alpha_n V @.\\
\Spec\,k_0 @<\psi_1<< \Spec\,k_1 @<\psi_2<< \cdots @<\psi_n<< \Spec\,k_n @<\psi_{n+1} <<  \cdots
\end{CD}
\] 
is a commutative diagram such that 
\begin{enumerate}
\item 
the lower sequence is induced by a sequence of field extensions
\[
k_0 \subset k_1 \subset k_2 \subset \cdots, 
\]
\item 
$X_n$ is a normal variety over $k_n$ such that $k_n$ is algebraically closed in $K(X_n)$ for any $n$, 
\item 
the induced morphism $X_n \to (X_{n-1} \times_{k_{n-1}} k_n)^N_{\red}$ 
is an isomorphism, and $\varphi_n$ coincides 
with the induced morphism for any $n$, 
\end{enumerate}
then it holds that 
\[
\#\{n \in \Z_{>0}\,|\, \theta_n:X_n \to X_{n-1} \times_{k_{n-1}} k_n\text{ is not an isomorphism}\} \leq \gamma_0,
\]
where $\theta_n$ denotes the induced morphism. 
\end{thm}

\begin{proof}
We first reduce the problem to the case when $X_0$ is proper over $k_0$. 
Fix an open immersion $j:X_0 \hookrightarrow Y_0$ 
to a proper normal variety $Y_0$ over $k_0$. 
Set $Y_n:=(Y_{n-1} \times_{k_{n-1}} k_n)_{\red}^N$ for any $n$. 
%Since each $X_n$ coincides with $(X_0 \times_{k_0} k_n)_{\red}^N$, 
Then we obtain a commutative diagram with the induced morphisms: 
\[
\begin{CD}
X_0 @<\varphi_1 << X_1 @<\varphi_2 << \cdots @<\varphi_n << X_n @<\varphi_{n+1}<< \cdots\\
@VVj_0 V @VVj_1 V @. @VVj_n V @.\\
Y_0 @<\varphi'_1 << Y_1 @<\varphi'_2 << \cdots @<\varphi'_n << Y_n @<\varphi'_{n+1}<< \cdots\\
@VV\beta_0 V @VV\beta_1 V @. @VV\beta_n V @.\\
\Spec\,k_0 @<\psi_1<< \Spec\,k_1 @<\psi_2<< \cdots @<\psi_n<< \Spec\,k_n @<\psi_{n+1} <<  \cdots. 
\end{CD}
\] 
Note that each $j_n:X_n \to Y_n$ is an open immersion and all the upper squares are cartesian. 
Replacing $X_0$ by $Y_0$, 
the problem is reduced to the case when $X_0$ is proper over $k_0$. 

We now reduce the problem to the case when $X_0$ is projective over $k_0$.  
Note that $X_0$ is assumed to be proper over $k_0$. 
Applying Chow's lemma, there exists a projective birational morphism 
$\mu_0:Z_0 \to X_0$ from a projective normal variety over $k_0$. 
Set $Z_n:=(Z_{n-1} \times_{k_{n-1}} k_n)_{\red}^N$ for any $n$. 
Then we obtain a commutative diagram with the induced morphisms: 
\[
\begin{CD}
Z_0 @<\varphi''_1 << Z_1 @<\varphi''_2 << \cdots @<\varphi''_n << Z_n @<\varphi''_{n+1}<< \cdots\\
@VV\mu_0 V @VV\mu_1 V @. @VV\mu_n V @.\\
X_0 @<\varphi_1 << X_1 @<\varphi_2 << \cdots @<\varphi_n << X_n @<\varphi_{n+1}<< \cdots\\
@VV\beta_0 V @VV\beta_1 V @. @VV\beta_n V @.\\
\Spec\,k_0 @<\psi_1<< \Spec\,k_1 @<\psi_2<< \cdots @<\psi_n<< \Spec\,k_n @<\psi_{n+1} <<  \cdots. 
\end{CD}
\] 
Fix a positive integer $n$ such that the induced morphism 
\[
\widetilde{\theta}_n:Z_n \to Z_{n-1} \times_{k_{n-1}} k_n
\]
is an isomorphism. 
Let us prove that also $X_n \to X_{n-1} \times_{k_{n-1}} k_n$ 
is an isomorphism. 
The equation $(\mu_{n-1})_* \MO_{Z_{n-1}} = \MO_{X_{n-1}}$ is stable under flat base changes, 
hence we have 
\[
\mu'_*\MO_{Z_{n-1} \times_{k_{n-1}} k_n} = \MO_{X_{n-1} \times_{k_{n-1}} k_n}, 
\]
where $\mu':Z_{n-1} \times_{k_{n-1}} k_n \to X_{n-1} \times_{k_{n-1}} k_n$ denotes 
the induced morphism. 
Therefore, $X_{n-1} \times_{k_{n-1}} k_n$ is a normal variety. 
In particular,  $\theta_n:X_n \to X_{n-1} \times_{k_{n-1}} k_n$ is an isomorphism. 
Thus, the problem is reduced to the case when $X_0$ is projective over $k_0$.

From now on, we treat the case when $X_0$ is projective over $k_0$. 
Let $\kappa$ be an algebraic closure of $k_0$ and set $V:=(X_0 \times_{k_0} \kappa)^N_{\red}$. 
Fix an ample Cartier divisor $H$ on $V$ and set 
\begin{equation}\label{e1-termination}
\gamma_0:=\max\{1, (f^*K_{X_0} - K_V)\cdot H^{\dim V-1}\},
\end{equation}
where $f:V \to X_0$ denotes the induced morphism. 
Take a commutative diagram as in the statement. 
For the algebraic closure $\kappa'$ of $\bigcup_{n=0}^{\infty} k_n$, 
set $V':=V \times_{\kappa} \kappa'$. 
Let $f'_n:V' \to X_n$ be the induced morphism. 
By (\ref{e1-termination}), it holds that 
\begin{equation}\label{e2-termination}
\gamma_0=\max\{1, (f'^*K_{X_0} - K_{V'})\cdot H'^{\dim V'-1}\},
\end{equation}
where $H'$ denotes the pullback of $H$ to $V'$. 

We have that 
\begin{equation}\label{e3-termination}
K_{X_n} + D_n \sim \varphi_n^*K_{X_{n-1}}
\end{equation}
for some effective Weil divisor on $D_n$ on $X_n$ (Theorem \ref{t-T-PW}). 
Note that the following are equivalent. 
\begin{enumerate}
\renewcommand{\labelenumi}{(\roman{enumi})}
\item $D_n=0$. 
\item The induced morphism $\theta_n:X_n \to X_{n-1} \times_{k_{n-1}} k_n$ is an isomorphism. 
\end{enumerate}
Indeed, it is clear that (ii) implies (i). 
Conversely, assuming (i), we have that $X_{n-1} \times_{k_{n-1}} k_n$ is normal  (Theorem \ref{t-T-PW}), 
hence (ii) holds.

Fix an arbitrary positive integer $m$. 
By (\ref{e3-termination}), we obtain 
\begin{equation}\label{e4-termination}
K_{X_m} + \sum_{n=1}^m D_{n, m} \sim g_m^*K_{X_0}, 
\end{equation}
where $g_m:X_m \to X_0$ denotes the induced morphism and 
$D_{n, m}$ is the pullback of $D_n$ to $X_m$. 
We obtain 
\begin{equation}\label{e5-termination}
K_{V'}+E_m \sim f_m'^*K_{X_m} 
\end{equation}
for some effective Weil divisor $E_m$ on $V'$ (Theorem \ref{t-T-PW}). 
It follows from (\ref{e4-termination}) and (\ref{e5-termination}) that 
\[
K_{V'}+E_m + \sum_{n=1}^m f_m'^*D_{n, m}\sim
f_m'^*\left(K_{X_m} + \sum_{n=1}^m D_{n, m}\right) \sim 
f_m'^*g_m^*K_{X_0}= f_0'^*K_{X_0}. 
\]
It holds that 
{\small 
\begin{eqnarray*}
\gamma_0
%&:=&E_0 \cdot H^{\dim V-1} \\
&=& 
\max\{1, (f_0'^*K_{X_0}-K_{V'}) \cdot H'^{\dim V'-1}\} \\
&\geq& 
(f_0'^*K_{X_0}-K_{V'}) \cdot H'^{\dim V'-1}\\
&=& \left(E_m + \sum_{n=1}^m f_m'^*D_{n, m}\right) \cdot H'^{\dim V'-1}\\
&\geq& \left(\sum_{n=1}^m f_m^*D_{n, m}\right) \cdot H^{\dim V'-1}\\
&\geq& \#\{n \in \{1, ..., m\}\,|\, D_n \neq 0\}\\
&=& \#\{n \in \{1, ..., m\}\,|\, \theta_n:X_n \to X_{n-1} \times_{k_{n-1}} k_n\text{ is not an isomorphism}\},
\end{eqnarray*}
}
where the first equality holds by (\ref{e2-termination}) 
and the last equality follows from the equivalence between (i) and (ii). 
Since $m$ is chosen to be an arbitrary positive integer, 
the assertion holds.  
\end{proof}

\subsection{Relations to base changes and morphisms}\label{s5-el-ext}

\begin{prop}\label{p-el-ext-sep-bc}
Let $k$ be a field of characteristic $p>0$ and 
let $X$ be a variety over $k$ such that $k$ is algebraically closed in $K(X)$. 
Let $(k \subset l \subset k', \varphi:X' \to X)$ be an elemental extension of $(k, X)$. 
Let $k \subset k_1$ be a (not necessarily algebraic) separable extension. 
Set $l_1:=l \otimes_k k_1$, $k'_1:=k' \otimes_k k_1$, $X_1:=X \times_k k_1$, and 
$X'_1:=X' \times_k k_1$. 
Let $\varphi_1:X'_1 \to X_1$ be the induced morphism. 
Then the following hold. 
\begin{enumerate}
\item $X_1$ is a variety over $k_1$ such that $k_1$ is algebraically closed in $K(X_1)$. 
Furthermore, if $X$ is normal, then $X_1$ is normal. 
\item $l_1$ and $k'_1$ are fields. 
The induced ring extensions $k_1 \subset l_1 \subset k'_1$ are purely inseparable field extensions such that $k'^p_1 \subset k_1$. 
\item 
$[l_1:k_1]=[l:k],\,\, [k':l]=[k'_1:l_1],\,\, [k':k]=[k'_1:k_1]$. 
\item $(k_1 \subset l_1 \subset k_1', \varphi_1:X'_1 \to X_1)$ is 
an elemental extension of $(k_1, X_1)$. 
\end{enumerate}
\end{prop}

\begin{proof}
The assertion (1) follows from Corollary \ref{c-sep-rn-commu}. 
Let us show (2). 
Since $k \subset k'$ is purely inseparable and $k \subset k_1$ is separable, 
we have that $k'_1=k' \otimes_k k_1$ is a field. %as it is reduced and univ homeo to Spec k
The same argument implies that $l_1$ is a field. 
It is clear that $k_1'^p \subset k_1$ 
and the extensions $k_1 \subset l_1 \subset k'_1$ are purely inseparable. 
Hence, (2) holds. 
The assertion (3) is obvious. 
 
Let us show (4). 
Since $X \times_k l$ is reduced, so is its separable base change $X_1 \times_{k_1} l_1$ (Lemma \ref{l-sep-rn-commu}(2)). 
Hence, $X_1 \times_{k_1} l_1 = X \times_k l_1$ is an integral scheme (Proposition \ref{p-prelim-girre}). 
Then it follows from 
Proposition \ref{p-sep-rn-commu}(2) that 
$X'_1 \to X_1 \times_{k_1} l_1$ is the normalisation of $X_1 \times_{k_1} l_1$. 
Thus (4) holds. 
\end{proof}

\begin{prop}\label{p-el-ext-alg-fib}
Let $k$ be a field of characteristic $p>0$ and 
let $f:X \to Y$ be a proper $k$-morphism of normal varieties over $k$ 
such that $f_*\MO_X=\MO_Y$ and 
$k$ is algebraically closed in $K(Y)$ 
(this condition automatically implies that $k$ is algebraically closed in $K(X)$). 
Let $(k \subset l \subset k_{X'}, \varphi:X' \to X)$ be an elemental extension of $(k, X)$. 
Then 
\begin{enumerate}
\item $Y \times_k l$ is an integral scheme. 
\end{enumerate}
Furthermore, if $Y' \to Y \times_k l$ is the normalisation of $Y \times_k l$ and $k_{Y'}$ denotes the 
algebraic closure of $l$ in $K(Y')$, 
then the following hold. 
\begin{enumerate}
\setcounter{enumi}{1}
\item $(k \subset l \subset k_{Y'}, \psi:Y' \to Y)$ is an elemental extension of $(k, Y)$, 
where $\psi$ denotes the induced morphism.  
\item $l \subset k_{X'}$ factors through $l \subset k_{Y'}$: 
\[
l \subset k_{Y'} \subset k_{X'}.
\]
\item If $Y \times_k l$ is not normal, then $X \times_k l$ is not normal. 
\end{enumerate}
\end{prop}

\begin{proof}
Let us show (1). 
Note that the induced morphism $f \times_k l:X \times_k l \to Y \times_k l$ satisfies 
$(f \times_k l)_*\MO_{X \times_k l}=\MO_{Y \times_k l}$. 
Since $X \times_k l$ is an integral scheme, so is $Y \times_k l$, 
hence (1) holds. 
If $X \times_k l$ is normal, then so is $Y \times_k l$. 
Thus (4) holds. 
The remaining assertions (2) and (3) are obvious. 
\end{proof}

\begin{prop}\label{p-el-ext-bir}
Let $k$ be a field of characteristic $p>0$ and 
let $f:X \dashrightarrow Y$ be a birational map over $k$ 
of varieties over $k$ such that $k$ is algebraically closed in $K(X)=K(Y)$. 
Let $(k \subset l \subset k', \varphi:X' \to X)$ be an elemental extension of $(k, X)$. 
Then 
\begin{enumerate}
\item 
$Y \times_k l$ is an integral scheme which is birational to $X \times_k l$.  
\end{enumerate}
Furthermore, if $Y' \to Y \times_k l$ is the normalisation of $Y \times_k l$ 
and $\psi:Y' \to Y$ denotes the induced morphism, then the following hold. 
\begin{enumerate}
\setcounter{enumi}{1}
\item 
$(k \subset l \subset k', \psi:Y' \to Y)$ is an elemental extension of $(k, Y)$. 
%This is called the elemental extension of $(k, Y)$ corresponding to $(k \subset l \subset k', \varphi:X' \to X)$. 
\item 
If the rational map $f:X \dashrightarrow Y$ is a morphism, then 
the following diagram of the induced morphisms 
\[
\begin{CD}
X @<<< X \times_k l @<<< X'\\
@VVfV @VVf \times_k lV @VVf'V\\
Y @<<< Y \times_k l @<<< Y'
\end{CD}
\] 
is commutative. 
Furthermore, $f':X' \to Y'$ is birational. 
\item 
If $f$ is an open immersion, then so are $f \times_k l$ and $f'$. 
In particular, if $Y \times_k l$ is normal, then so is $X \times_k l$. 
\item 
Assume that $Y$ is normal and $f$ is a proper morphism. 
Then it holds that 
$(f \times_k l)_*\MO_{X \times_k l} = \MO_{Y \times_k l}$ and $f'_*\MO_{X'}=\MO_{Y'}$. 
In particular, if $X \times_k l$ is normal, then so is $Y \times_k l$. 
\end{enumerate}
\end{prop}

\begin{proof}
The assertion (1) follows from Lemma \ref{l-geom-red-birat}. 
The remaining assertions (2)--(5) are clear. 
\end{proof}

\begin{prop}\label{p-el-ext-sm}
Let $k_Y$ be a field of characteristic $p>0$ and 
let $f:X \to Y$ be a smoorh $k_Y$-morphism of normal 
varieties over $k_Y$ such that $k_Y$ is algebraically closed in $K(Y)$. 
Let $k_X$ be the algebraic closure of $k_Y$ in $K(X)$. 
Let $(k_Y \subset l_Y \subset k'_Y, \psi:Y' \to Y)$ be an elemental extension of $(k_Y, Y)$. 
Set $l_X:=l_Y \otimes_{k_Y} k_X$ and $X':=Y' \times_Y X$. 
Then the following hold. 
\begin{enumerate}
\item 
$l_X$ and $k'_Y \otimes_{k_Y} k_X$ are fields. 
\item 
$X \times_k l$ is an integral scheme and the induced morphism 
$X' \to X \times_{k_X} l_X$ is the normalisation of $X \times_{k_X} l_X$. 
%$X'$ is a normal variety such that $k'_X$ is algebraically closed in $K(X')$. 
\item 
$(k_X \subset l_X \subset k'_X, \varphi:X' \to X)$ is an elemental extension of $(k_X, X)$, 
where $\varphi$ denotes the induced morphism and $k'_X$ is the algebraic closure of $l_X$ in $K(X')$. 
\item 
$[l_X:k_X]=[l:k]$ and $[k'_X:l_X] \geq [k'_Y:l_Y]$. 
\item 
If $Y \times_{k_Y} l_Y$ is normal, then $X \times_{k_X} l_X$ is normal. 
\item 
If $f$ is surjective and $X \times_{k_X} l_X$ is normal, 
then $Y \times_{k_Y} l_Y$ is normal.
\end{enumerate} 
\end{prop}

\begin{proof}
By Proposition \ref{p-sm-sep}, $k_X$ is a finite separable extension of $k_Y$. 
Hence, (1) holds. 
Let us show (2). 
Since $X \times_{k_X} l_X \simeq X \times_{k_Y} l_Y$, 
we have the commutative diagram 
\begin{equation}\label{e1-el-ext-sm}
\begin{CD}
X @<<< X \times_{k_X} l_X @<<< X'\\
@VVf V @VV f \times_{k_Y} l_YV @VVf'V\\
Y @<<< Y \times_{k_Y} l_Y @<<< Y'\\
\end{CD}
\end{equation}
where both the squares are cartesian. 
In particular, the vertical arrows are smooth morphisms. 
Since $Y$ is normal, 
it follows from \cite[Lemma 2.2(2)]{Tan18} that 
the horizontal arrows are finite universal homeomorphisms.  
Therefore, $X \times_{k_X} l_X$ is an integral scheme and 
$X' \to X \times_{k_X} l_X$ is the normalisation of $X \times_{k_X} l_X$. 
Thus, (2) holds.

The assertion (3) follows from (1) and (2). 
It is clear that $[l_X:k_X]=[l:k]$. 
Since we have injection $k_X \otimes_{k_Y} k'_Y \hookrightarrow k'_X$, we obtain  
\[
[k'_X:l_X] \geq [k_X \otimes_{k_Y} k'_Y:k_X \otimes_{k_Y} l_Y] 
=[k'_Y:l_Y].
\]
Thus (4) holds. 
The assertions (5) and (6) follow from the fact that 
the left square of (\ref{e1-el-ext-sm}) is cartesian. 
\end{proof}

\section{Capacity of denormalising base changes $\gamma(X/k)$}

In this section, we first introduce capacity of denormalising base changes $\gamma(X/k)$ 
in Subsection \ref{ss1-gamma} (Definition \ref{d-gamma}). 
We then discuss behaviour of $\gamma(X/k)$ 
under base changes (Subsection \ref{ss2-gamma}) and morphisms (Subsection \ref{ss3-gamma}). 
A prominent property of $\gamma(X/k)$ is the following formula: 
\[
K_Y+(p-1) \sum_{i=1}^{\gamma(X/k)} C_i \sim f^*K_X
\]
where we set $Y:=(X \times_k k^{1/p^{\infty}})_{\red}^N$ and $f:Y \to X$ 
denotes the induced morphism (Theorem \ref{t-cano-gamma}). 
%On the other hand, it seems to be difficult to compute $\gamma(X/k)$, 
%as it is defined as a maximum in an abstract situation. 
%Because of this reason, the author does not know whether 
%$\gamma(X/k)$ changes even under algebraic separable base changes (cf. Proposition \ref{p-gamma-sep-bc})

\subsection{Definition and basic properties}\label{ss1-gamma}

In this subsection, we first introduce capacity of denormalising base changes $\gamma(X/k)$ 
(Definition \ref{d-gamma}). 
We then describe an alternative definition in terms of elemental extensions (Proposition \ref{p-def-gamma2}). 
Finally, we give a formula for canonical divisors (Theorem \ref{t-cano-gamma}).

\begin{dfn}\label{d-gamma}
Let $k$ be a field of characteristic $p>0$ 
and let $X$ be a normal variety over $k$ such that $k$ is algebraically closed in $K(X)$. 
Set $k_0:=k$ and $X_0:=X$. 
We define $\gamma(X/k)$ as the minimum non-negative integer $\gamma$ such that 
if 
\[
\begin{CD}
X_0 @<\varphi_1 << X_1 @<\varphi_2 << \cdots @<\varphi_n << X_n @<\varphi_{n+1}<< \cdots\\
@VV\alpha_0 V @VV\alpha_1 V @. @VV\alpha_n V @.\\
\Spec\,k_0 @<\psi_1<< \Spec\,k_1 @<\psi_2<< \cdots @<\psi_n<< \Spec\,k_n @<\psi_{n+1} <<  \cdots
\end{CD}
\] 
is a commutative diagram such that 
\begin{enumerate}
\item 
the lower sequence is induced by a sequence of purely inseparable extensions
\[
k_0 \subset k_1 \subset k_2 \subset \cdots, 
\]
\item 
$X_n$ is a normal variety over $k_n$ such that $k_n$ is algebraically closed in $K(X_n)$ for any $n$, 
\item 
the induced morphism $X_n \to (X_{n-1} \times_{k_{n-1}} k_n)^N_{\red}$ 
is an isomorphism, and $\varphi_n$ coincides 
with the induced morphism for any $n$, 
\end{enumerate}
then it holds that 
\[
\#\{n \in \Z_{>0}\,|\, \theta_n:X_n \to X_{n-1} \times_{k_{n-1}} k_n\text{ is not an isomorphism}\} \leq \gamma,
\]
where $\theta_n$ denotes the induced morphism. 
It follows from Theorem \ref{t-termination} that such $\gamma(X/k)$ exists. 
We call $\gamma(X/k)$ the {\em capacity of denormalising extensions} of $X/k$. 
\end{dfn}

\begin{prop}\label{p-def-gamma}
Let $k$ be a field of characteristic $p>0$ 
and let $X$ be a normal variety over $k$ such that $k$ is algebraically closed in $K(X)$. 
Then the following are equivalent. 
\begin{enumerate}
\item $X$ is geometrically normal over $k$. 
\item $\gamma(X/k)=0$. 
\end{enumerate}
\end{prop}

\begin{proof}
It is clear that (1) implies (2). 
Let us prove that (2) implies (1). 
Assume (2). 
Then $(X \times_k k^{1/p^{\infty}})_{\red}^N \to X \times_k k^{1/p^{\infty}}$ 
is an isomorphism by Definition \ref{d-gamma}. 
Therefore, (1) holds. 
\end{proof}

\begin{prop}\label{p-def-gamma2}
Let $k$ be a field of characteristic $p>0$ 
and let $X$ be a normal variety over $k$ such that $k$ is algebraically closed in $K(X)$. 
Set $k_0:=k$, $X_0:=X$, and $\gamma:=\gamma(X/k)=\gamma(X_0/k_0)$. 
Then there exists a commutative diagram 
\[
\begin{CD}
X_0 @<\varphi_1 << X_1 @<\varphi_2 << \cdots @<\varphi_{\gamma} << X_{\gamma}\\
@VV\alpha_0 V @VV\alpha_1 V @. @VV\alpha_{\gamma} V @.\\
\Spec\,k_0 @<\psi_1<< \Spec\,k_1 @<\psi_2<< \cdots @<\psi_{\gamma}<< \Spec\,k_{\gamma}
\end{CD}
\] 
such that 
\begin{enumerate}
\item 
the lower sequence is induced by a sequence of purely inseparable extensions
\[
k_0 \subset k_1 \subset \cdots \subset k_{\gamma}, 
\]
\item 
$X_n$ is a normal variety over $k_n$ such that $k_n$ is algebraically closed in $K(X_n)$ for any $n$, 
\item 
the induced morphism $X_n\to (X_{n-1} \times_{k_{n-1}} k_n)^N_{\red}$ 
is an isomorphism, $\varphi_n$ coincides 
with the induced morphism for any $n$, and 
\item 
for any $1 \leq n \leq \gamma$, there is a commutative diagram 
\[
\begin{CD}
X_{n-1} @<\varphi'_n<< X'_{n-1} @<\varphi''_n<< X_n\\
@VV\alpha_{n-1}V @VV\beta_{n-1}V @VV\alpha_n V\\
\Spec\,k_{n-1} @<\psi'_n<< \Spec\,k'_{n-1} @<\psi''_n<< \Spec\,k_n
\end{CD}
\]
such that $\varphi_n=\varphi'_n \circ \varphi''_n$, 
$\psi_n=\psi'_n \circ \psi''_n$, 
the left square is a cartesian, and 
the right square is an elemental extension of internal degree $p$ which is not cartesian. 
\end{enumerate}
\end{prop}

\begin{proof}
By Definition \ref{d-gamma}, there exists a commutative diagram 
\[
\begin{CD}
X_0 @<\varphi_1 << X_1 @<\varphi_2 << \cdots @<\varphi_{\gamma} << X_{\gamma}\\
@VV\alpha_0 V @VV\alpha_1 V @. @VV\alpha_{\gamma} V @.\\
\Spec\,k_0 @<\psi_1<< \Spec\,k_1 @<\psi_2<< \cdots @<\psi_{\gamma}<< \Spec\,k_{\gamma}
\end{CD}
\] 
such that (1)--(3) and the following property $(4)'$ hold. 
\begin{enumerate}
\item[$(4)'$] 
the induced morphism $X_n \to (X_{n-1} \times_{k_{n-1}} k_n)_{\red}^N$ is not an isomorphism 
for any $n$. 
\end{enumerate}
We apply Theorem \ref{t-dec-elem} to each square 
\[
\begin{CD}
X_{n-1} @<\varphi_n<< X_n\\
@VV\alpha_{n-1}V @VV\alpha_nV\\
\Spec\,k_{n-1} @>\psi_n>> \Spec\,k_n. 
\end{CD}
\]
By Definition \ref{d-gamma}, this square is decomposed into the following commutative diagram: 
\[
\begin{CD}
X_{n-1} @<<< Y_{n-1} @<<< Z_n @<<< X_n\\
@VV\alpha_{n-1}V @VVV @VVV @VV\alpha_n V\\
\Spec\,k_{n-1} @<<< \Spec\,k^Y_{n-1} @<<< \Spec\,k^Z_n @<<< \Spec\,k_n
\end{CD}
\]
where the left and the right squares are cartesian and 
the middle square is an elemental extension of internal degree $p$ which is not cartesian. 
After replacing $(k_n, X_n)$ suitably, the assertion holds. 
\end{proof}

\begin{thm}\label{t-cano-gamma}
Let $k$ be a field of characteristic $p>0$ 
and let $X$ be a normal variety over $k$ such that $k$ is algebraically closed in $K(X)$. 
For $\kappa \in \{k^{1/p^{\infty}}, \overline k\}$, 
set $Y:=(X \times_k \kappa)_{\red}^N$ and let $f:Y \to X$ be the induced morphism. 
If $X$ is not geometrically normal over $k$, 
then there exist effective Weil divisors $C_1, ..., C_{\gamma(X/k)}$ such that 
the linear equivalence 
\begin{equation}\label{e1-cano-gamma}
K_Y + (p-1)\sum_{i=1}^{\gamma(X/k)} C_i \sim f^*K_X
\end{equation}
holds. 
%Furthermore, if $X$ is proper over $k$, then 
Furthermore, each $C_i$ can be chosen to be nonzero if 
$X$ is geometrically reduced or $X$ is proper over $k$. 
\end{thm}

\begin{proof}
If $\kappa=k^{1/p^{\infty}}$, 
then the assertion follows from Theorem \ref{t-T-PW} and Proposition \ref{p-def-gamma2}. 
If $\kappa= \overline k$, 
then the problem is reduced to the case when $\kappa=k^{1/p^{\infty}}$, 
since $Y=(X \times_k \overline k)_{\red}^N \to X$ factors through 
$(X \times_k k^{1/p^{\infty}})_{\red}^N \to X$. 
\end{proof}

\subsection{Base changes}\label{ss2-gamma}

\begin{prop}\label{p-gamma-sep-bc}
Let $k$ be a field of characteristic $p>0$ 
and let $X$ be a normal variety over $k$ such that $k$ is algebraically closed in $K(X)$. 
Let $k \subset l$ be a (not necessarily algebraic) separable extension. 
Then $\gamma(X/k) \leq \gamma(X \times_k l/l)$. 
\end{prop}

\begin{proof}
The assertion follows from Proposition \ref{p-def-gamma2} and 
the fact that a separable base change of an elemental extension is again 
an elemental extension (Proposition \ref{p-el-ext-sep-bc}). 
\end{proof}

\begin{rem}
In the situation of Proposition \ref{p-gamma-sep-bc}, 
it is natural to hope the equation: $\gamma(X/k) \geq \gamma(X \times_k l/l)$. 
The author does not know how to prove this. 
\end{rem}

\begin{prop}\label{p-gamma-F-bc}
Let $k$ be a field of characteristic $p>0$ 
and let $X$ be a normal variety over $k$ such that $k$ is algebraically closed in $K(X)$. 
Assume that $X$ is not geometrically normal over $k$. 
Set $X':=(X \times_k k^{1/p})_{\red}^N$. 
Then $\gamma(X/k)  \geq \gamma(X'/k^{1/p})+1$ 
(note that $k^{1/p}$ is algebraically closed in $K(X')$ by Lemma \ref{l-Fbc-alg-cl}). 
\end{prop}

\begin{proof}
Since $X' \to (X \times_k k^{1/p})_{\red}^N$ is not an isomorphism (Proposition \ref{p-geom-ht1}(3)), 
the assertion follows from Definition \ref{d-gamma}.  
\end{proof}

\subsection{Behaviour under morphisms}\label{ss3-gamma}

\begin{prop}\label{p-gamma-sm}
Let $k$ be a field of characteristic $p>0$ 
and let $f:X \to Y$ be a smooth surjective $k$-morphism of normal $k$-varieties such that 
$k$ is algebraically closed in $K(Y)$. 
Let $k_X$ be the algebraic closure of $k_Y$ in $K(X)$. 
Then $\gamma(X/k_X) \geq \gamma(Y/k)$. 
\end{prop}

\begin{proof}
Set $k_0:=k$, $Y:=Y_0$, and $\gamma:=\gamma(Y/k)=\gamma(Y_0/k_0)$. 
Proposition \ref{p-def-gamma2} implies that there is a commutative diagram 
\[
\begin{CD}
Y_0 @<\varphi_1 << Y_1 @<\varphi_2 << \cdots @<\varphi_{\gamma} << Y_{\gamma}\\
@VV\alpha_0 V @VV\alpha_1 V @. @VV\alpha_{\gamma} V @.\\
\Spec\,k_0 @<\psi_1<< \Spec\,k_1 @<\psi_2<< \cdots @<\psi_{\gamma} << \Spec\,k_{\gamma}
\end{CD}
\] 
that satisfies the properties as in Proposition \ref{p-def-gamma2}. 
In particular, for any $1 \leq n \leq \gamma$, there is a commutative diagram 
\[
\begin{CD}
Y_{n-1} @<\varphi'_n<< Y'_{n-1} @<\varphi''_n<< Y_n\\
@VV\alpha_{n-1}V @VV\beta_{n-1}V @VV\alpha_n V\\
\Spec\,k_{n-1} @<\psi'_n<< \Spec\,k'_{n-1} @<\psi''_n<< \Spec\,k_n
\end{CD}
\]
such that $\varphi_n=\varphi'_n \circ \varphi''_n$, 
$\psi_n=\psi'_n \circ \psi''_n$, 
the left square is a cartesian, and 
the right square is an elemental extension of internal degree $p$ which is not cartesian. 

We inductively define $X_n$ and $X'_n$ by base changes, i.e. 
$X_n$ and $X'_n$ complete the following diagram consisting of cartesian diagrams: 
\[
\begin{CD}
X_{n-1} @<\widetilde{\varphi}'_n<< X'_{n-1} @<\widetilde{\varphi}''_n<< X_n\\
@VVf_{n-1}V @VVf'_{n-1}V @VVf_nV\\
Y_{n-1} @<\varphi'_n<< Y'_{n-1} @<\varphi''_n<< Y_n.\\
\end{CD}
\]
Hence, the horizontal arrows are finite universal homeomorphisms 
and the vertical arrows are smooth. 
Therefore, $X_n$ and $X'_n$ are normal varieties for any $n$. 
Let $k_{X_n}$ and $k_{X'_n}$ be the algebraic closure of $k$ in $K(X_n)$ and $K(X'_n)$, 
respectively. 
It follows from Proposition \ref{p-el-ext-sm} 
that $(X_n, k_{X_n})$ is an elemental extension 
of $(k_{X'_{n-1}}, X'_{n-1})$ such that $X'_{n-1} \times_{k_{X'_{n-1}}} k_{X_n}$ is not normal. 
Therefore, Definition \ref{d-gamma} implies that 
$\gamma(X/k) \geq \gamma=\gamma(Y/k)$. 
\end{proof}

\begin{prop}\label{p-gamma-opem-imm}
Let $k$ be a field of characteristic $p>0$ 
and let $f:X \to Y$ be an open immersion over $k$ of normal $k$-varieties such that 
$k$ is algebraically closed in $K(Y)$. 
Then $\gamma(X/k) \leq \gamma(Y/k)$. 
\end{prop}

\begin{proof}
We can apply a similar argument to the one of Proposition \ref{p-gamma-sm}. 
\end{proof}

\section{Frobenius length of geometric non-normality $\ell_F(X/k)$}

In this section, 
we first introduce Frobenius length of geometric non-normality $\ell_F(X/k)$ 
in Subsection \ref{ss1-lF} (Definition \ref{d-F-length}). 
We then discuss behaviour of $\ell_F(X/k)$ 
under base changes (Subsection \ref{ss2-lF}) and morphisms (Subsection \ref{ss3-lF}). 
In Subsection \ref{ss4-lF}, 
we establish Frobenius factorisation property (Theorem \ref{t-Frob-factor1}). 
In Subsection \ref{ss5-lF}, 
we give a bound on $\Q$-Gorenstein index of $(X \times_k k^{1/p^{\infty}})_{\red}^N$ 
under mild assumptions (Theorem \ref{t-Qgor-index}).

\subsection{Definition and basic properties}\label{ss1-lF}

\begin{dfn}\label{d-F-length}
Let $k$ be a field of characteristic $p>0$. 
Let $X$ be a normal variety over $k$ such that $k$ is algebraically closed in $K(X)$. 
We define $\ell_F(X/k)$ by 
{\small 
\[
\ell_F(X/k):=\min\{\ell \in \Z_{\geq 0}\,|\, 
(X \times_k k^{1/p^{\ell}})_{\red}^N \text{ is geometrically normal over }k^{1/p^{\ell}}\}.
\]
}
Note that the minimum in the right hand side exists (Remark \ref{r-F-length}). 
We call $\ell_F(X/k)$ the {\em Frobenius length of geometric non-normality} of $X/k$. 
\end{dfn}

\begin{rem}\label{r-F-length}
We use notation as in Definition \ref{d-F-length}. 
Lemma \ref{l-N-descend} implies that $(X \times_k k^{1/p^{\ell}})_{\red}^N$ 
is geometrically normal over $k^{1/p^{\ell}}$ for any $\ell \gg 0$. 
Furthermore, we have that 
\[
(X \times_k k^{1/p^{\ell_1}})_{\red}^N \times_{k^{1/p^{\ell_1}}} k^{1/p^{\ell_2}} 
\simeq (X \times_k k^{1/p^{\ell_2}})_{\red}^N
\]
for any integers $\ell_1, \ell_2$ such that $\ell_F(X/k) \leq \ell_1 < \ell_2 \leq \infty$.  
\end{rem}

\begin{rem}\label{r-F-length2}
We use notation as in Definition \ref{d-F-length}. 
Proposition \ref{p-geom-ht1}(3) implies 
that the following are equivalent. 
\begin{enumerate}
\item $X$ is geometrically normal over $k$. 
\item $\ell_F(X/k)=0$. 
\end{enumerate}
\end{rem}

\subsection{Base changes}\label{ss2-lF}

\begin{prop}\label{p-length-sep-bc}
Let $k$ be a field of characteristic $p>0$. 
Let $X$ be a normal variety over $k$ such that $k$ is algebraically closed in $K(X)$. 
Let $k \subset l$ be a (not necessarily algebraic) separable field extension. 
Then $\ell_F(X/k) =\ell_F(X\times_k l/l)$. 
\end{prop}

\begin{proof}
Fix $e \in \Z_{\geq 0}$. 
Set $Y:=X \times_k l$. 
Since $k^{1/p^e} \subset l^{1/p^e}$ is a separable extension, 
it follows from Proposition \ref{p-sep-rn-commu}(2) that 
\[
(X \times_k k^{1/p^e})_{\red}^N \times_{k^{1/p^e}} l^{1/p^e} 
\simeq 
(X \times_k k^{1/p^e} \times_{k^{1/p^e}} l^{1/p^e})_{\red}^N 
\simeq (Y \times_l l^{1/p^e})_{\red}^N. 
\]
Therefore, $(X \times_k k^{1/p^e})_{\red}^N$ is geometrically normal over $k^{1/p^e}$ 
if and only if 
$(Y \times_l l^{1/p^e})_{\red}^N$ is geometrically normal over $l^{1/p^e}$. 
\end{proof}

\begin{prop}\label{p-length-alg-bc}
Let $k$ be a field of characteristic $p>0$. 
Let $X$ be a normal variety over $k$ such that $k$ is algebraically closed in $K(X)$. 
Let $k \subset k'$ be a field extension 
such that $k'$ is algebraically closed in $X':=(X \times_k k')_{\red}^N$. 
Then the inequality $\ell_F(X/k) \geq \ell_F(X'/k')$ holds. 
\end{prop}

\begin{proof}
By Proposition \ref{p-length-sep-bc}, we may assume that $k \subset k'$ is purely inseparable. 
Set $\ell:=\ell_F(X/k)$. 
It follows from Definition \ref{d-F-length} that $(X \times_k k^{1/p^{\ell}})_{\red}^N$ is geometrically normal over $k^{1/p^{\ell}}$. 
It holds that 
\[
(X \times_k k^{1/p^{\ell}})_{\red}^N \times_{k^{1/p^{\ell}}} k'^{1/p^{\ell}}
\simeq 
(X' \times_{k'} k'^{1/p^{\ell}})_{\red}^N, 
\]
because each of the both hand sides coincides with the normalisation of $X$ 
in the composite field of $K(X)$ and $k'^{1/p^{\ell}}$ (Proposition \ref{p-bc-comp-field}). 
Therefore, $(X' \times_{k'} k'^{1/p^{\ell}})_{\red}^N$ is geometrically normal over $k'^{1/p^{\ell}}$. 
Then it holds that $\ell_F(X/k) = \ell \geq \ell_F(X'/k')$, as desired. 
\end{proof}

\begin{prop}\label{p-length-Frob-bc}
Let $k$ be a field of characteristic $p>0$. 
Let $X$ be a normal variety over $k$ such that $k$ is algebraically closed in $K(X)$. 
Let $r$ be a non-negative integer such that $\ell_F(X/k) \geq r$. 
Set $X':=(X \times_k k^{1/p^r})_{\red}^N$. 
Then it holds that $\ell_F(X'/k^{1/p^r}) =\ell_F(X/k)-r$. 
\end{prop}

\begin{proof} 
Fix $r$ as in the statement.  
Set $k':=k^{1/p^r}$. 
For any integer $\ell$ such that $\ell \geq r$, 
it holds that 
\[
(X' \times_{k'} k'^{1/p^{\ell-r}})_{\red}^N 
=
((X \times_k k^{1/p^r})_{\red}^N \times_{k^{1/p^r}} k^{1/p^{\ell}})_{\red}^N 
\simeq  (X \times_k k^{1/p^{\ell}})_{\red}^N,
\]
where the last isomorphism holds, 
because each of both hand sides is isomorphic to the normalisation 
of $X$ in the composite field of $K(X)$ and $k^{1/p^{\ell}}$ 
(Proposition \ref{p-bc-comp-field}). 
Then the assertion follows from Definition \ref{d-F-length}. 
\end{proof}

\subsection{Behaviour under morphisms}\label{ss3-lF}

\begin{prop}\label{p-ell-sm-des}
Let $k$ be a field of characteristic $p>0$ 
and let $f:X \to Y$ be a smooth $k$-morphism of normal $k$-varieties such that 
$k$ is algebraically closed in $K(Y)$. 
Let $k_X$ be the algebraic closure of $k$ in $K(X)$. 
Then the following hold. 
\begin{enumerate}
\item 
$\ell_F(X/k_X) \leq \ell_F(Y/k)$.
\item 
If $f$ is surjective, then $\ell_F(X/k_X) = \ell_F(Y/k)$.
\end{enumerate}
\end{prop}

\begin{proof}
Set $\ell(X):=\ell_F(X/k_X)$ and $\ell(Y):=\ell_F(Y/k)$. 
For any positive integer $e$, the following three morphisms are smooth: 
\begin{itemize}
\item $X \times_k k^{1/p^e} \to Y \times_k k^{1/p^e}$. 
\item $(X \times_k k^{1/p^e})_{\red} \to (Y \times_k k^{1/p^e})_{\red}$. 
\item $(X \times_k k^{1/p^e})_{\red}^N \to (Y \times_k k^{1/p^e})_{\red}^N$. 
\end{itemize}
Since $k \subset k_X$ is a finite separable extension (Proposition \ref{p-sm-sep}), 
it holds that 
$
k_X^{1/p^e} \simeq k_X \otimes_k k^{1/p^e}. 
$
Hence, we obtain 
\[
(X \times_k k^{1/p^e})_{\red}^N \simeq (X \times_{k_X} k_X^{1/p^e})_{\red}^N. 
\]
%Therefore, if $(Y \times_k k^{1/p^e})_{\red}^N$ is normal, then $(X \times_{k_X} k_X^{1/p^e})_{\red}^N$ is normal. 

Let us show (1). 
By Definition \ref{d-F-length}, $(Y \times_k k^{1/p^{\ell(Y)}})^N_{\red}$ is geometrically normal over 
$k^{1/p^{\ell(Y)}}$. 
Then $(X \times_k k^{1/p^{\ell(Y)}})_{\red}^N$ is geometrically normal over $k^{1/p^{\ell(Y)}}$. 
This implies that $(X \times_{k_X} k_X^{1/p^{\ell(Y)}})_{\red}^N$ 
is geometrically normal over $k_X^{1/p^{\ell(Y)}}$. 
Therefore, (1) holds.

Let us show (2). 
By Definition \ref{d-F-length}, $(X \times_{k_X} k_X^{1/p^{\ell(X)}})_{\red}^N$ is geometrically normal over 
$k_X^{1/p^{\ell(X)}}$. 
Then $(X \times_k k^{1/p^{\ell(X)}})_{\red}^N$ is geometrically normal over $k^{1/p^{\ell(X)}}$. 
Therefore, $(Y \times_k k^{1/p^{\ell(X)}})_{\red}^N$ is geometrically normal over $k^{1/p^{\ell(X)}}$. 
Hence, it holds that $\ell(X) \leq \ell(Y)$. 
Since we have already shown $\ell(X) \geq \ell(Y)$ in (1), 
the assertion (2) holds. 
\end{proof}

\begin{prop}
Let $k$ be a field of characteristic $p>0$ 
and let $f:X \to Y$ be a proper birational $k$-morphism of normal $k$-varieties such that 
$k$ is algebraically closed in $K(X)$. 
Then it holds that $\ell_F(X/k) \geq \ell_F(Y/k)$.
\end{prop}

\begin{proof}
Set $\ell(X):=\ell_F(X/k)$. 
The induced morphism 
\[
f':X':=(X \times_k k^{1/p^{\ell(X)}})_{\red}^N \to (Y \times_k k^{1/p^{\ell(X)}})_{\red}^N=:Y'
\] 
is a proper birational morphism of normal varieties over $k^{1/p^{\ell(X)}}$. 
In particular, we have that $f'\MO_{X'}=\MO_{Y'}$. 
Since $X'$ is geometrically normal over $k^{1/p^{\ell(X)}}$, 
$Y'$ is geometrically normal over $k^{1/p^{\ell(X)}}$. 
This implies $\ell_F(X/k) = \ell(X) \geq \ell_F(Y/k)$. 
\end{proof}

\subsection{Frobenius factorisation}\label{ss4-lF}

The purpose of this subsection is to prove the following theorem. 

\begin{thm}\label{t-Frob-factor1}
Let $k$ be a field of characteristic $p>0$. 
Let $X$ be a normal variety over $k$ such that $k$ is algebraically closed in $K(X)$. 
Set $\ell:=\ell_F(X/k)$. 
Then the $\ell$-th iterated absolute Frobenius morphism of $X \times_k k^{1/p^{\infty}}$ 
factors through the induced morphism 
$g:(X \times_k k^{1/p^{\infty}})_{\red}^N \to X \times_k k^{1/p^{\infty}}$: 
\[
F^{\ell}_{X \times_k k^{1/p^{\infty}}}:X \times_k k^{1/p^{\infty}}\to (X \times_k k^{1/p^{\infty}})_{\red}^N \xrightarrow{g} X \times_k k^{1/p^{\infty}}.
\]
\end{thm}

\begin{proof}
We first prove that the $\ell$-th iterated 
absolute Frobenius morphism $F_{X \times_k k^{1/p^{\ell}}}^{\ell}$ of $X \times_k k^{1/p^{\ell}}$ 
factors through the induced morphism $\widetilde{g}:(X \times_k k^{1/p^{\ell}})_{\red}^N \to 
X \times_k k^{1/p^{\ell}}$: 
\begin{equation}\label{e1-Frob-factor1}
F_{X \times_k k^{1/p^{\ell}}}^{\ell}:
X \times_k k^{1/p^{\ell}} \xrightarrow{\widetilde{h}}
(X \times_k k^{1/p^{\ell}})_{\red}^N
\xrightarrow{\widetilde{g}} 
X \times_k k^{1/p^{\ell}}.
\end{equation}
By Proposition \ref{p-bc-comp-field}, we have a factorisation: 
$F^{\ell}_X:X \xrightarrow{h_1} (X \times_k k^{1/p^{\ell}})_{\red}^N \xrightarrow{g_1} X$, 
where $g_1$ is the induced morphism. 
For the first projection $\pi:X \times_k k^{1/p^{\ell}} \to X$, we obtain the following diagram 
\begin{center}
\begin{tikzpicture}[auto]
\node (X1) at (0, 0) {$X$};
\node (Z) at (4, 0) {$(X \times_k k^{1/p^{\ell}})_{\red}^N$};
\node (X2) at (8, 0) {$X$};
\node (Y1) at (0, 3) {$X \times_k k^{1/p^{\ell}}$}; 
\node (Y2) at (8, 3) {$X \times_k k^{1/p^{\ell}}$};

\draw (2,2) circle (0.15);
\draw (6.7,1) circle (0.15);
\draw (3.5,-0.6) circle (0.15);
\node (A) at (4, 0.85) {$\vartriangle$};

\draw[->] (Y1) to node {$\small{F_{X \times_k k^{1/p^{\ell}}}^{\ell}}$} (Y2);
\draw[->] (Y2) to node {$\scriptstyle \pi$} (X2);
\draw[->] (Y1) to node {$\scriptstyle \pi$} (X1);
\draw[->, bend right] (X1) to node {$\scriptstyle{F_X^{\ell}}$} (X2);
\draw[->] (X1) to node {$\scriptstyle{\zeta}$} (Y2);
\draw[->] (X1) to node {$\scriptstyle{h_1}$} (Z);
\draw[->] (Z) to node {$\scriptstyle{g_1}$} (X2);
\draw[->] (Z) to node {$\scriptstyle{\widetilde{g}}$} (Y2);
\end{tikzpicture}
%\end{equation}
\end{center}
where we can check that 
there is the morphism $\zeta$ such that 
$F^{\ell}_{X \times_k k^{1/p^{\ell}}} = \zeta \circ \pi$ and 
$F^{\ell}_X = \pi \circ \zeta$. 
Clearly, the above three diagrams with circles are commutative.

%SInce \MO_X \to \MO_{X \times_k k^{1/p}}$ is inje, it suffices to prove the commutativity 
%after compositing with pi, then the assertion is clear. 
We now prove that also the diagram with the triangle is commutative. 
To this end, we may assume that $X$ is affine. 
Since all the morphisms are affine, all the schemes in the diagram are affine. 
%We denote $\MO_V(V)$ by $\MO_V$ and let $\varphi^*:\MO_W \to \MO_V$ be the induced 
%ring homomorphism induced by a morphism $\varphi:V \to W$. 
Then it suffices to prove that the following  two induced ring homomorphisms coincide: 
\[
\zeta^*, h_1^* \circ \widetilde g^*:\MO_{X \times_k k^{1/p^{\ell}}}(X \times_k k^{1/p^{\ell}}) \to \MO_X(X), 
\]
where each $(-)^*$ denotes the induced ring homomorphism. 
Since $\MO_X(X)$ is reduced, it is enough to prove that 
$\zeta^*(a^{p^{\ell}}) =  h_1^* \circ \widetilde g^* (a^{p^{\ell}})$ for any 
$a \in \MO_{X \times_k k^{1/p^{\ell}}}(X \times_k k^{1/p^{\ell}})$. 
By $a^{p^{\ell}} \in \pi^*(\MO_X(X)) ={\rm Im}(\pi^*)$, 
it is sufficient to show that $\zeta^* \circ \pi^* = h_1^* \circ \widetilde g^* \circ \pi^*$, 
which follows from a diagram chase. Thus,  the diagram with the triangle is commutative. 
%As $\pi^*$ is injective, it is sufficient to prove 
%$\pi^* \circ \zeta^* \circ \pi^* = \pi^* \circ h_1^* \circ \widetilde g^* \circ \pi^*$, 

%We have shown that the above diagram is commutative. 
Set $\widetilde h:=h_1 \circ \pi$. 
Since the above diagram is commutative, 
we have the required factorisation (\ref{e1-Frob-factor1}). 
Taking the base change $(\ref{e1-Frob-factor1}) \times_{k^{1/p^{\ell}}} k^{1/p^{\infty}}$, 
we obtain 
\[
F^{\ell}_{X \times_k k^{1/p^{\ell}}}\times_{k^{1/p^{\ell}}} k^{1/p^{\infty}}
:Y \to (X \times_k k^{1/p^{\infty}})_{\red}^N \to 
X \times_k k^{1/p^{\infty}},
\]
where $Y$ is obtained as the base change and 
the isomorphism 
\[
(X \times_k k^{1/p^{\ell}})_{\red}^N \times_{k^{1/p^{\ell}}} k^{1/p^{\infty}} 
\simeq (X \times_k k^{1/p^{\infty}})_{\red}^N
\]
follows from Remark \ref{r-F-length}. 
Note that the $\ell$-th iterated absolute Frobenius morphism 
$F^{\ell}_{X \times_k k^{1/p^{\infty}}}$ 
factors through the base change of the $\ell$-th iterated absolute Frobenius
$
F^{\ell}_{X \times_k k^{1/p^{\ell}}}\times_{k^{1/p^{\ell}}} k^{1/p^{\infty}},
$
as desired. 
\end{proof}

\begin{rem}
The author does not know whether $\ell_F(X/k)$ 
coincides with the minimum non-negative integer $\nu$ such that 
the $\nu$-th iterated absolute Frobenius morphism of $X \times_k k^{1/p^{\infty}}$ 
factors through the induced morphism 
$g:(X \times_k k^{1/p^{\infty}})_{\red}^N \to X \times_k k^{1/p^{\infty}}$. 
\end{rem}

\subsection{$\Q$-Gorenstein index}\label{ss5-lF}

In this subsection, 
we give a bound of $\Q$-Gorenstein index of $(X \times_k k^{1/p^{\infty}})_{\red}^N$ 
by using $\ell_F(X/k)$. 
The following proposition is a key result.

\begin{prop}\label{p-Qgor-index} 
Let $k$ be a field of characteristic $p>0$ and 
let $X$ be a normal variety over $k$ such that $k$ is algebraically closed in $K(X)$. 
Set $\ell:=\ell_F(X/k)$ and $X_n:=(X \times_k k^{1/p^n})_{\red}^N$ for any $n \in \Z_{\geq 0}$.
Fix $\kappa \in \{k^{1/p^{\infty}}, \overline k\}$ and 
let $Y:= (X \times_k \kappa)_{\red}^N$.  
Let 
\[
\begin{CD}
X_0 @<\varphi_1 << X_1 @<\varphi_2 << \cdots @<\varphi_{\ell} << X_{\ell} @<\widetilde{\varphi}<< Y\\
@VV\alpha_0 V @VV\alpha_1 V @. @VV\alpha_{\ell} V @VV\beta V\\
\Spec\,k @<\psi_1<< \Spec\,k^{1/p} @<\psi_2<< \cdots @<\psi_{\ell}<< \Spec\,k^{1/p^{\ell}} @<\widetilde{\psi}<< \Spec\,k^{1/p^{\infty}}
\end{CD}
\] 
be the commutative diagram consisting of the induced morphisms. 
Then the following hold. 
\begin{enumerate}
\item 
For any integer $i$ with $1 \leq i \leq \ell$, there exists an effective Weil divisor $C_i$ on $X_i$ such that the linear equivalence 
\[
K_{X_i} + (p-1) C_i \sim \varphi_i^* K_{X_{i-1}}
\]
holds. 
Furthermore, each $C_i$ can be chosen to be nonzero 
if $X$ is geometrically reduced or $X$ is proper over $k$. 
\item 
Let $D_i$ be the pullback of $C_i$ to $Y$. 
Then the linear equivalence 
\[
K_Y + (p-1) \sum_{i=1}^{\ell}D_i \sim f^*K_X
\]
holds, where $f:Y \to X$ denotes the induced morphism. 
\item 
For any integer $i$ with $1 \leq i \leq \ell$, 
the $i$-th iterated absolute Frobenius morphism $F^i_{X_i}$ of $X_i$ 
factors through the induced morphism $X_i \to X_0=X$: 
\[
F^i_{X_i}:X_i \to X \to X_i. 
\]
%In particular, if $X$ is regular, then $p^iC_i$ is Cartier. 
\end{enumerate}
\end{prop}

\begin{proof}
Let us show (1). 
Fix an integer $i$ with $1 \leq i \leq \ell$. 
By Definition \ref{d-F-length}, $X_{i-1}$ is not geometrically normal over $k^{1/p^{i-1}}$. 
Hence, $X_{i-1} \times_{k^{1/p^{i-1}}} k^{1/p^i}$ is not normal. 
Therefore, the assertion (1) holds by Theorem \ref{t-T-PW}. 

Let us show (2). 
It follows from Remark \ref{r-F-length} that the rightmost square 
in the commutative diagram in the statement is cartesian. 
Hence, (1) implies (2). 

Let us show (3). 
By Proposition \ref{p-bc-comp-field}, we have $K(X_i)^{p^i} \subset K(X) \subset K(X_i)$. 
Therefore, the $i$-th iterated absolute Frobenius morphism $F^i_{X_i}$ of $X_i$ 
factors through the induced morphism $X_i \to X_0=X$: 
\[
F^i_{X_i}:X_i \to X \to X_i. 
\]
Thus (3) holds. 
\end{proof}

\begin{thm}\label{t-Qgor-index} 
Let $k$ be a field of characteristic $p>0$ and 
let $X$ be a normal variety over $k$ such that $k$ is algebraically closed in $K(X)$. 
Set $\ell:=\ell_F(X/k)$. 
Fix $\kappa \in \{k^{1/p^{\infty}}, \overline k\}$ and 
let $Y:= (X \times_k \kappa)_{\red}^N$.  
Assume that there exists a positive integer $n$ such that 
if $E$ is a prime divisor $E$ on $X$, then $nE$ is a Cartier divisor 
(e.g. if $X$ is regular, then $n$ can be chosen to be $1$). 
Then $np^{\ell} K_Y$ is a Cartier divisor. 
\end{thm}

\begin{proof}
We use the same notation as in Proposition \ref{p-Qgor-index}. 
By Proposition \ref{p-Qgor-index}(2), we have 
\[
K_Y + (p-1) \sum_{i=1}^{\ell}D_i \sim f^*K_X.
\]
Since the Weil divisor $np^iC_i$ is the pullback of $nC_i$ by $F^i_{X_i}$, 
$np^iC_i$ is a pullback of a Cartier divisor on $X$ (Proposition \ref{p-Qgor-index}(3)). 
Hence, $np^iC_i$ is Cartier. 
Since $nK_X$ and $np^{\ell}C_i$ are Cartier for any $i$, 
it holds that also $np^{\ell}K_Y$ is Cartier. 
\end{proof}

\section{Frobenius length of geometric non-reducedness $m_F(X/k)$}

In this section, 
we first introduce Frobenius length of geometric non-reducedness $m_F(X/k)$ 
in Subsection \ref{ss1-mF} (Definition \ref{d-mF}). 
We then discuss behaviour of $m_F(X/k)$ 
under base changes (Subsection \ref{ss2-mF}) and morphisms (Subsection \ref{ss3-mF}). 
In Subsection \ref{ss4-mF}, 
we establish Frobenius factorisation property.

\subsection{Definition and basic properties}\label{ss1-mF}

In this subsection, 
we introduce Frobenius length of geometric non-reducedness $m_F(X/k)$ 
(Definition \ref{d-mF}). 
By definition, we have $m_F(X/k) \leq \ell_F(X/k)$. 
After removing suitable closed subsets, 
the equality holds (Proposition \ref{p-gene-lm}). 

\begin{dfn}\label{d-mF}
Let $k$ be a field of characteristic $p>0$. 
Let $X$ be a variety over $k$ 
such that $k$ is algebraically closed in $K(X)$. 
Set 
\[
m_F(X/k):=
\]
\[
\min\{m \in \Z_{\geq 0}\,|\, 
(X \times_k k^{1/p^{m}})_{\red} \text{ is geometrically reduced over }k^{1/p^{m}}\},
\]
whose existence is guaranteed by Remark \ref{r-mF2}. 
We call $m_F(X/k)$ the {\em Frobenius length of geometric non-reducedness} of $X/k$. 
\end{dfn}

\begin{rem}\label{r-mF2}
We use notation as in Definition \ref{d-mF}. 
Lemma \ref{l-N-descend} implies that $(X \times_k k^{1/p^m})_{\red}$ 
is geometrically reduced over $k^{1/p^m}$ for any $m \gg 0$. 
Furthermore, we have that 
\[
(X \times_k k^{1/p^{m_1}})_{\red} \times_{k^{1/p^{m_1}}} k^{1/p^{m_2}} 
\simeq (X \times_k k^{1/p^{m_2}})_{\red}
\]
for any $m_F(X/k) \leq m_1 < m_2 \leq \infty$. 
\end{rem}

\begin{rem}\label{r-mF1}
We use notation as in Definition \ref{d-mF}. 
Proposition \ref{p-geom-ht1} implies 
that the following are equivalent. 
\begin{enumerate}
\item $X$ is geometrically reduced over $k$. 
\item $m_F(X/k)=0$. 
\end{enumerate}
\end{rem}

\begin{rem}\label{r-mF3}
We use notation as in Definition \ref{d-mF}. 
If $X'$ is a non-empty open subset of $X$, 
then it holds that $m_F(X/k) = m_F(X'/k)$. 
Indeed, for any non-negative integer $n$, 
Proposition \ref{p-geom-red-birat} implies that the following are equivalent. 
\begin{itemize}
\item 
$(X \times_k k^{1/p^n})_{\red}$ is geometrically reduced over $k^{1/p^n}$. 
\item 
$(X' \times_k k^{1/p^n})_{\red}$ is geometrically reduced over $k^{1/p^n}$. 
\end{itemize}
\end{rem}

\begin{rem}\label{r-m-ell}
Let $k$ be a field of characteristic $p>0$. 
Let $X$ be a normal variety over $k$ 
such that $k$ is algebraically closed in $K(X)$. 
Then it follows from Definition \ref{d-F-length} and Definition \ref{d-mF} 
 that $m_F(X/k) \leq \ell_F(X/k)$. 
\end{rem}

\begin{prop}\label{p-gene-lm}
Let $k$ be a field of characteristic $p>0$. 
Let $X$ be a variety over $k$ 
such that $k$ is algebraically closed in $K(X)$. 
Then there exists a non-empty open subset $Y$ of $X$ 
such that $Y$ is normal and 
the equation 
\[
m_F(X/k) = m_F(Z/k) = \ell_F(Z/k)
\] 
holds 
for any non-empty open subset $Z$ of $Y$.  
\end{prop}

\begin{proof}
For $m:=m_F(X/k)$, 
$(X \times_k k^{1/p^{m}})_{\red}$ is geometrically reduced over $k^{1/p^{m}}$. 
Note that $(X \times_k k^{1/p^{n}})_{\red}$ is an integral scheme 
for any $n \in \Z_{\geq 0}$ (cf. Proposition \ref{p-prelim-girre}). 
Since $(X \times_k k^{1/p^m})_{\red}$ is geometrically reduced over $k^{1/p^m}$, 
there exists a non-empty open subset $Y$ of $X$ such that 
$Y$ is regular and $(Y \times_k k^{1/p^m})_{\red}$ is smooth over $k^{1/p^m}$. 
%and $(Y \times_k k^{1/p^n})_{\red}$ is regular for any $0 \leq n \leq m$. 
Fix a non-empty open subset $Z$ of $Y$. 
Since $(Z \times_k k^{1/p^m})_{\red}$ is smooth over $k^{1/p^m}$ 
we have $\ell_F(Z/k) \leq m$. 
On the other hand, Remark \ref{r-mF3} deduces that $m_F(X/k) = m_F(Z/k)$. 
To summarise, it holds that 
\[
\ell_F(Z/k) \leq m = m_F(X/k) = m_F(Z/k) \leq \ell_F(Z/k), 
\]
where the last inequality follows from Remark \ref{r-m-ell}. 
\end{proof}

\subsection{Base changes}\label{ss2-mF}

\begin{prop}\label{p-length-sep-bc2}
Let $k$ be a field of characteristic $p>0$. 
Let $X$ be a variety over $k$ such that $k$ is algebraically closed in $K(X)$. 
Let $k \subset l$ be a  (not necessarily algebraic) separable field extension. 
Then $m_F(X/k) =m_F(X\times_k l/l)$. 
\end{prop}

\begin{proof}
The assertion follows from 
the same argument as in 
Proposition \ref{p-length-sep-bc}.% and Proposition \ref{p-gene-lm}.
\end{proof}

\begin{prop}\label{p-length-alg-bc2}
Let $k$ be a field of characteristic $p>0$. 
Let $X$ be a variety over $k$ such that $k$ is algebraically closed in $K(X)$. 
Let $k \subset k'$ be a field extension 
such that $k'$ is algebraically closed in $K(X')$, where $X':=(X \times_k k')_{\red}$. 
Then the inequality $m_F(X/k) \geq m_F(X'/k')$ holds. 
\end{prop}

\begin{proof}
By Proposition \ref{p-length-sep-bc2}, we may assume that $k \subset k'$ is a purely inseparable extension. 
Then the assertion follows from Proposition \ref{p-length-alg-bc} and Proposition \ref{p-gene-lm}.
\end{proof}

\begin{prop}\label{p-length-Frob-bc2}
Let $k$ be a field of characteristic $p>0$. 
Let $X$ be a variety over $k$ such that $k$ is algebraically closed in $K(X)$. 
Let $r$ be a non-negative integer such that $m_F(X/k) \geq r$. 
Set $X':=(X \times_k k^{1/p^r})_{\red}$. 
Then it holds that $m_F(X'/k^{1/p^r}) =m_F(X/k)-r$. 
\end{prop}

\begin{proof}
The assertion follows from Proposition \ref{p-length-Frob-bc} and Proposition \ref{p-gene-lm}.
\end{proof}

\subsection{Behaviour under morphisms}\label{ss3-mF}

\begin{prop}
Let $k$ be a field of characteristic $p>0$ 
and let $f:X \to Y$ be a smooth $k$-morphism of $k$-varieties such that 
$k$ is algebraically closed in $K(Y)$. 
Let $k_X$ be the algebraic closure of $k$ in $K(X)$. 
Then it holds that $m_F(X/k_X) = m_F(Y/k)$
\end{prop}

\begin{proof}
Thanks to Remark \ref{r-mF3} and Proposition \ref{p-gene-lm}, 
we may assume that $X$ and $Y$ are normal, 
$m_F(X/k_X)=\ell_F(X/k_X), m_F(Y/k)=\ell_F(Y/k)$, and 
$f$ is surjective. 
Then the assertion follows from Proposition \ref{p-ell-sm-des}(2). 
\end{proof}

\begin{prop}
Let $k$ be a field of characteristic $p>0$. 
Let $X$ be a variety over $k$ 
such that $k$ is algebraically closed in $K(X)$. 
Let $X'$ be a variety over $k$ which is birational to $X$. 
Then it holds that $m_F(X/k)=m_F(X'/k)$. 
\end{prop}

\begin{proof}
The assertion follows from Remark \ref{r-mF3}. 
\end{proof}

\subsection{Frobenius factorisation}\label{ss4-mF}

\begin{thm}\label{t-Frob-factor2}
Let $k$ be a field of characteristic $p>0$. 
Let $X$ be a variety over $k$ such that $k$ is algebraically closed in $K(X)$. 
Set $m:=m_F(X/k)$. 
Then the $m$-th iterated absolute Frobenius morphism of $X \times_k k^{1/p^{\infty}}$ 
factors through the induced morphism 
$g:(X \times_k k^{1/p^{\infty}})_{\red} \to X \times_k k^{1/p^{\infty}}$: 
\[
F^m_{X \times_k k^{1/p^{\infty}}}:X \times_k k^{1/p^{\infty}}\to (X \times_k k^{1/p^{\infty}})_{\red} \xrightarrow{g} X \times_k k^{1/p^{\infty}}.
\]
\end{thm}

\begin{proof}
The same argument as in Theorem \ref{t-Frob-factor1} works. 
\end{proof}

\section{Thickening exponents $\epsilon(X/k)$}

In this section, 
we introduce thickening exponent $\epsilon(X/k)$ in Subsection \ref{ss1-epsilon} 
(Definition \ref{d-epsilon}). 
We then discuss behaviour of $\epsilon(X/k)$ 
under base changes (Subsection \ref{ss2-epsilon}) and morphisms (Subsection \ref{ss3-epsilon}).

\subsection{Definition and basic properties}\label{ss1-epsilon}

In this subsection, we define thickening exponent $\epsilon(X/k)$. 
To this end, we prove Theorem \ref{t-def-epsilon}, 
which asserts that $\length_R R$ is a power of $p$ for the local ring $R$ at the generic point of $X \times_k k^{1/p^{\infty}}$. 
Although this result itself is known (cf. \cite[Ch. I, Proposition 7.1.1]{Kol96}), 
we give an alternative proof, 
which simultaneously enables us to establish a relation to elemental extensions (Proposition \ref{p-epsilon-elem}). 
%The strategy of the proof is induction by using elemental extensions. 
%Hence, we study the behaviour of $\length_R R$ under elemental extensions (Proposition \ref{p-epsilon-elem}). 
We start with the following auxiliary result.

\begin{lem}\label{l-length-dim}
Let $(R, \m)$ be an artinian local ring and 
let $F$ be a coefficient field of $R$, i.e. $F$ is a field such that 
$F$ is a subring of $R$ and the composite ring homomorphism $F \hookrightarrow R \to R/\m$ is an isomorphism. 
Let $M$ be a finitely generated $R$-module. 
Then 
\[
\length_R\,M=\dim_F M.
\]
\end{lem}

\begin{proof}
Since $R$ is an artinian local ring, 
it follows from \cite[Theorem 6.4]{Mat89} that there is a composition sequence of $R$-submodules of $M$ 
\[
M=:M_0 \supset M_1 \supset \cdots \supset M_n=0
\]
such that $M_i/M_{i+1} \simeq R/\m$ for any $i$. 
Then we have that $\dim_F M=n$. 
On the other hand, the Jordan--H\"older theorem deduces that $\length_R\,M=n$, as desired.  
\end{proof}

\begin{prop}\label{p-epsilon-elem}
Let $k$ be a field of characteristic $p>0$ and 
let $X$ be a variety such that $k$ is algebraically closed in $K(X)$. 
Let $k \subset l$ be a purely inseparable field extension such that $X \times_k l$ 
is an integral scheme. 
Let $X'$ be the normalisation of $X \times_k l$ and let $k'$ be the algebraic closure of $l$ in $K(X')$. 
For $\kappa:=k^{1/p^{\infty}}$, we consider $l$ and $k'$ 
as subfields of $\kappa$. 
Set $R$ and $R'$ to be the local rings of $X \times_k \kappa$ and $X' \times_{k'} \kappa$ 
at the generic points, respectively. 
Then it holds that 
\[
\length_R\,R = [k':l] \length_{R'}\,R'.
\]
\end{prop}

\begin{proof}
Set $K:=K(X)$ and $K':=K(X')$. 
Since $X \times_k l$ is an integral scheme and $X'$ is its normalisation, 
we obtain 
\[
K \otimes_k l =K(X) \otimes_k l = K(X \times_k l) = K(X')=K'. 
\]
Let $S$ be the $k'$-algebra 
such that we set $S:=K' \otimes_l k'$ as an $l$-algebra and 
we consider $S$ as a $k'$-algebra via 
the with the second coprojection $k' \to K' \otimes_l k'$. 
Set $S':=K' \otimes_{k'} k'$, which is a $k'$-algbera. 
It holds that  
\begin{itemize}
\item $R'= K' \otimes_{k'} \kappa = S' \otimes_{k'} \kappa$, and 
\item $R= K \otimes_k \kappa= (K \otimes_k l) \otimes_l \kappa=K' \otimes_l \kappa 
=S \otimes_{k'} \kappa$. 
\end{itemize} 
%where we consider $K' \otimes_l k'$ as a $k'$-algebra 
%via the second coprojection $k' \to K' \otimes_l k'$. 
We have the natural surjective $k'$-algebra homomorphism 
\[
\pi:S=K' \otimes_l k' \to K' \otimes_{k'} k'=S', \qquad x \otimes_l y \mapsto  x \otimes_{k'} y. 
\] 
Applying the tensor product $(-) \otimes_{k'} \kappa$, we obtain the surjective ring homomorphism 
\[
\pi \otimes_{k'} \kappa: R=S \otimes_{k'} \kappa \to S' \otimes_{k'} \kappa =R'. 
\]
We consider $S'$ and $R'$ as an $S$-algebra and an $R$-algebra  
via $\pi$ and $\pi \otimes_{k'} \kappa$, respectively. 

Since $S'$ is a field and $S$ is an artinian local ring, 
we obtain the $S$-algebra isomorphism $S' \simeq S_{\red}$, i.e. $\Ker(\pi)$ is equal to the nilradical of $S$. 
We have $l$-algebra homomorphisms: 
\[
{\rm id}:K' \xrightarrow{\varphi} K' \otimes_l k' = S \xrightarrow{\pi} S' = K' \otimes_{k'} k' \xrightarrow{\theta, \simeq} K', 
\]
where $\varphi$ is the first coprojection and 
$\theta$ is defined by $\theta(x \otimes y)=xy$ for $x \in K'$ and $y \in k'$. 
Set $F:=\varphi(K')$. 
Then it follows from Lemma \ref{l-length-dim} that $\length_S\,S=\dim_F S=[k':l]$. 
Since $S$ is an artinian local ring, there exists a composition sequence of $S$-submodules 
\[
S=:M_0 \supset M_1 \supset \cdots \supset M_{[k':l]}
\]
with $S$-module isomorphisms $M_i/M_{i+1} \simeq S_{\red} \simeq S'$ for any $i$. 
For $N_i:=M_i \otimes_{k'} \kappa$, we obtain a sequence of $R$-submodules: 
\[
R=N_0 \supset N_1 \supset \cdots \supset N_{[k':l]}
\]
with $R$-module isomorphisms $N_i/N_{i+1} \simeq S' \otimes_{k'} \kappa \simeq R'$ for any $i$. 
Therefore, it holds that  
\[
\length_R R = [k':l] \length_R R' =[k':l] \length_{R'} R',
\]
where the second equality holds because the $R$-module structure of $R'$ is defined by 
a surjective ring homomorphism $\pi \otimes_{k'} \kappa:R \to R'$ of artinian local rings. 
\end{proof}

\begin{thm}\label{t-def-epsilon}
Let $k$ be a field of characteristic $p>0$ and 
let $X$ be a variety such that $k$ is algebraically closed in $K(X)$. 
Set $R$ to be the local ring of $X \times_k k^{1/p^{\infty}}$ at the generic point. 
Then the following hold. 
\begin{enumerate}
\item 
There exists a non-negative integer $\epsilon$ such that  
\[
\length_R\,R = p^{\epsilon}.
\]
\item 
For the non-negative integer $\epsilon$ as in (1), 
$\epsilon=0$ if and only if $X$ is geometrically reduced over $k$. 
\end{enumerate}
\end{thm}

\begin{proof}
Replacing $X$ by a suitable open subset of $X$, we may assume that $X$ is normal. 
Let us prove (1) by induction on $\length_R\,R$. 
If $\length_R\,R=1$, then there is nothing to show. 
Assume that $\length_R\,R>1$. 
Then $R$ is not a field. 
In particular, $X$ is not geometrically reduced over $k$. 
%After replacing $X$ by a suitable non-empty open subset of $X$, 
It follows from Theorem \ref{t-schroer} that there exists 
a purely inseparable field extension $k \subset l$ such that 
\begin{itemize}
\item $X \times_k l$ is an integral scheme, 
\item $k':=k^{1/p}$ is equal to the algebraic closure of $l$ in $K(X \times_k l)$, 
\item $l \subsetneq k'$, and 
\item $(k \subset l \subset k', \varphi:X' \to X)$ is an elemental extension for 
$X':=(X \times_k l)^N$ and the induced morphism $\varphi$. 
\end{itemize}
By Proposition \ref{p-epsilon-elem}, 
we have 
\[
\length_R\,R = [k':l] \length_{R'}\,R',
\]
where $R'$ denotes the local ring of $X' \times_{k'} k'^{1/p^{\infty}}$ at the generic point. 
By the induction hypothesis, 
we obtain $\length_{R'}\,R'=p^{\epsilon'}$ for some non-negative integer $\epsilon'$. 
Since $l \subset k'$ is a finite purely inseparable extension (Lemma \ref{l-el-ext}(2)), (1) holds. 

Let us show (2). 
If $X$ is geometrically reduced over $k$, then $R$ is a field. 
Hence we have that $\epsilon=0$. 
Conversely, if $\epsilon=0$, then $R$ is a field. 
Then it follows from Lemma \ref{l-geom-red-birat} that $X$ is geometrically reduced over $k$. 
Thus, (2) holds. 
\end{proof}

\begin{dfn}\label{d-epsilon}
Let $k$ be a field of characteristic $p>0$ and 
let $X$ be a variety such that $k$ is algebraically closed in $K(X)$. 
Set $R$ to be the local ring of $X \times_k k^{1/p^{\infty}}$ at the generic point. 
We define $\epsilon(X/k)$ as the non-negative integer such that 
\[
\length_R\,R = p^{\epsilon(X/k)},
\]
whose existence is guaranteed by Theorem \ref{t-def-epsilon}(1). 
We call $\epsilon(X/k)$ the {\em thickening exponent} of $X/k$. 
\end{dfn}

\begin{rem}\label{r-def-epsilon}
We use notation as in Definition \ref{d-epsilon}. 
By Theorem \ref{t-def-epsilon}(2), 
the following are equivalent. 
\begin{enumerate}
\item $X$ is geometrically reduced over $k$. 
\item $\epsilon(X/k)=0$. 
\end{enumerate}
\end{rem}

%\begin{rem}\label{r-def-epsilon2}
%We use notation as in Definition \ref{d-epsilon}. 
%If $X'$ is a non-empty open subset of $X$, 
%then it holds that $\epsilon(X/k) = \epsilon(X'/k)$. 
%\end{rem}

\subsection{Base changes}\label{ss2-epsilon}

\begin{prop}\label{p-epsilon-des-bc}
Let $k$ be a field of characteristic $p>0$ and 
let $X$ be a variety over $k$ 
such that $k$ is algebraically closed in $K(X)$. 
For a field extension $k \subset k'$, 
assume that $X \times_k k'$ is a variety over $k'$ 
and $k'$ is algebraically closed in $K(X \times_k k')$. 
Then it holds that $\epsilon(X/k)=\epsilon(X \times_k k'/k')$. 
\end{prop}

\begin{proof}
Set $R:=K(X) \otimes_k k^{1/p^{\infty}}$ and $R':=K(X \times_k k') \otimes_{k'} k'^{1/p^{\infty}}$.  
Take a composition sequence of $R$-submodules: 
\[
R=:M_0 \supset M_1 \supset \cdots \supset M_{p^{\epsilon(X/k)}}=0. 
\]
It holds that $M_i/M_{i+1} \simeq R_{\red}$ for any $i$.

We now apply the base change $(-) \otimes_{k^{1/p^{\infty}}} k'^{1/p^{\infty}}$ to the above sequence. 
Note that we have 
\[
R \otimes_{k^{1/p^{\infty}}} k'^{1/p^{\infty}} 
=(K(X) \otimes_k k^{1/p^{\infty}}) \otimes_{k^{1/p^{\infty}}} k'^{1/p^{\infty}} 
=K(X) \otimes_k k'^{1/p^{\infty}}
\]
\[
=(K(X) \otimes_k k') \otimes_{k'} k'^{1/p^{\infty}} =K(X \times_k k') \otimes_{k'} k'^{1/p^{\infty}}
=R'. 
\]
Therefore, for $M'_i := M_i \otimes_{k^{1/p^{\infty}}} k'^{1/p^{\infty}}$, 
we obtain a sequence of $R'$-submodules 
\[
R'=:M'_0 \supset M'_1 \supset \cdots \supset M'_{p^{\epsilon(X/k)}}=0 
\]
such that $M'_i/M'_{i+1} \simeq R_{\red} \otimes_{k^{1/p^{\infty}}} k'^{1/p^{\infty}}$. 

It is enough to show that 
$R_{\red} \otimes_{k^{1/p^{\infty}}} k'^{1/p^{\infty}} \simeq R'_{\red}$. 
Applying $(-) \otimes_{k^{1/p^{\infty}}} k'^{1/p^{\infty}}$ 
to the surjective ring homomorphism $\pi:R \to R_{\red}$, we obtain a 
surjective ring homomorphism 
\[
R' \to R_{\red} \otimes_{k^{1/p^{\infty}}} k'^{1/p^{\infty}}. 
\]
Since $k^{1/p^{\infty}}$ is a perfect field, 
$R_{\red} \otimes_{k^{1/p^{\infty}}} k'^{1/p^{\infty}}$ is reduced (cf. Lemma \ref{l-sep-rn-commu}(2)). 
Since $R'$ is a local artinian ring, 
we have that $R_{\red} \otimes_{k^{1/p^{\infty}}} k'^{1/p^{\infty}} \simeq R'_{\red}$. 
\end{proof}

\begin{prop}\label{p-epsilon-sep-bc}
Let $k$ be a field of characteristic $p>0$ and 
let $X$ be a variety over $k$ 
such that $k$ is algebraically closed in $K(X)$. 
Let $k \subset k'$ be a (not necessarily algebraic) separable field extension. 
Set $X':=X \times_k k'$. 
Then it holds that $\epsilon(X/k)=\epsilon(X'/k')$. 
\end{prop}

\begin{proof}
The assertion follows from Proposition \ref{p-epsilon-des-bc}. 
\end{proof}

\begin{prop}\label{p-epsilon-compute}
Let $k$ be a field of characteristic $p>0$ and 
let $X$ be a variety over $k$ 
such that $k$ is algebraically closed in $K(X)$. 
Let $k \subset l$ be a purely inseparable field extension. 
Set $X':=(X \times_k l)^N_{\red}$ and 
let $k'$ be the algebraic closure of $l$ in $K(X')$. 
Then the following hold. 
\begin{enumerate}
\item 
If $X \times_k l$ is reduced, then it holds that 
\[
\epsilon(X/k) = \epsilon(X'/k') + \log_p [k':l].
\]
\item 
$\epsilon(X/k) \geq \epsilon(X'/k')$. 
\item 
The following are equivalent. 
\begin{enumerate}
\item 
$\epsilon(X/k) = \epsilon(X'/k')$
\item 
$X \times_k l$ is reduced and $l=k'$. 
\end{enumerate}
\end{enumerate}
\end{prop}

\begin{proof}
The assertion (1) follows from Proposition \ref{p-epsilon-elem}. 
The assertion (2) holds by (1), Theorem \ref{t-dec-elem}, and Proposition \ref{p-epsilon-des-bc}. 
Let us show (3). 
It follows from (1) that (b) implies (a). 
By (1) and Theorem \ref{t-dec-elem}, (a) implies (b). 
\end{proof}

\begin{prop}\label{p-epsilon-alg-bc}
Let $k$ be a field of characteristic $p>0$ and 
let $X$ be a variety over $k$ 
such that $k$ is algebraically closed in $K(X)$. 
Let $k \subset l$ be a field extension and 
set $X':=(X \times_k l)^N_{\red}$ and 
let $k'$ be the algebraic closure of $l$ in $K(X')$. 
Then $\epsilon(X/k) \geq \epsilon(X'/k')$. 
\end{prop}

\begin{proof}
The assertion follows from 
Proposition \ref{p-epsilon-sep-bc} and Proposition \ref{p-epsilon-compute}. 
\end{proof}

\subsection{Behaviour under morphisms}\label{ss3-epsilon}

\begin{prop}\label{p-epsilon-birat}
Let $k$ be a field of characteristic $p>0$. 
Let $X$ be a variety over $k$ 
such that $k$ is algebraically closed in $K(X)$. 
Let $X'$ be a variety over $k$ which is birational to $X$. 
Then it holds that $\epsilon(X/k)=\epsilon(X'/k)$. 
\end{prop}

\begin{proof}
The assertion follows from Definition \ref{d-epsilon}. 
\end{proof}

\begin{prop}\label{p-epsilon-sm}
Let $k$ be a field of characteristic $p>0$ 
and let $f:X \to Y$ be a smooth $k$-morphism of $k$-varieties such that 
$k$ is algebraically closed in $K(Y)$. 
Let $k_X$ be the algebraic closure of $k$ in $K(X)$. 
Then it holds that $\epsilon(X/k_X)=\epsilon(Y/k)$. 
\end{prop}

\begin{proof}
By Proposition \ref{p-epsilon-birat}, 
we may assume that $f$ is surjective. 

Let us prove the assertion by induction on $\epsilon(Y/k)$. 
If $\epsilon(Y/k)=0$, then $Y$ is geometrically reduced over $k$ (Remark \ref{r-def-epsilon}). 
Since $f:X \to Y$ is smooth, $X$ is geometrically reduced over $k$. 
As $k \subset k_X$ is a finite separable extension (Proposition \ref{p-sm-sep}), 
$X$ is geometrically reduced over $k_X$. 
Hence Remark \ref{r-def-epsilon} implies that $\epsilon(X/k)=0$. 
Therefore, the assertion holds for the case when $\epsilon(Y/k)=0$. 

Assume $\epsilon(Y/k)>0$. 
By Proposition \ref{p-epsilon-birat}, we may assume that $X$ and $Y$ are regular. 
Remark \ref{r-def-epsilon} implies that $Y$ is not geometrically reduced. 
It follows from Theorem \ref{t-schroer}(1) that 
there is an elemental extension $(k \subset l \subset k^{1/p}, \psi:Y' \to Y)$ 
for $Y':=(Y \times_k k^{1/p})_{\red}^N$ and the induced morphism $\psi$. 
Moreover,  Theorem \ref{t-schroer}(2) implies that $[k^{1/p}:l]>1$. 
By Proposition \ref{p-el-ext-sm}, 
if $l_X:=l \otimes_k k_X$ and $X':=Y' \times_Y X$, 
then $(k_X \subset l_X \subset k'_X, \varphi:X' \to X)$ is an elemental extension, 
where $k'_X$ denotes the algebraic closure of $l_X$ in $X'$ and $\varphi$ is the induced morphism. 
We have that 
\[
k_X^{1/p} \supset k'_X \supset k_X \otimes_k k^{1/p} =k_X^{1/p}, 
\]
where the first inequality holds by Lemma \ref{l-el-ext}(1), 
the second one follows from the fact that $k_X \otimes_k k^{1/p}$ is a field, 
and the last equality holds, because $k \subset k_X$ is a finite separable extension (Proposition \ref{p-sm-sep}). 
Therefore, $(k_X \subset l_X \subset k_X^{1/p}, \varphi:X' \to X)$ is an elemental extension 
and we have $[k^{1/p}:l] = [k_X^{1/p}:l_X]$. 
It follows from Proposition \ref{p-epsilon-compute} that 
\begin{itemize}
\item 
$\epsilon(X/k_X) = \epsilon((X \times_{k_X} l_X)^N/k_X^{1/p}) + \log_p [k_X^{1/p}:l_X]$ and 
\item 
$\epsilon(Y/k) = \epsilon((Y \times_k l)^N/k^{1/p}) + \log_p [k^{1/p}:l]$. 
\end{itemize}
The induction hypothesis implies that 
\[
\epsilon((X \times_{k_X} l_X)^N/k_X^{1/p})=
\epsilon((Y \times_k l)^N/k^{1/p}).
\]
By $[k^{1/p}:l] = [k_X^{1/p}:l_X]$, we obtain $\epsilon(X/k_X)=\epsilon(Y/k)$. 
\end{proof}

\section{Relation between invariants}

We have introduced four invariants $\gamma (X/k), \ell_F (X/k), m_F (X/k)$, and $\epsilon (X/k)$. 
In \cite{Sch10}, Schr\"{o}er studied another invariant, 
called geometric generic embedding dimension. 
In this section, we summarise relations between these five invariants. 
We first recall geometric generic embedding dimensions by Schr\"{o}er (Subsection \ref{ss-gged}). 
We then  compare invariants related to geometric non-reducedness (Subsection \ref{ss-g-non-red}). 
In Subsection \ref{ss-g-non-norm}, we compare invariants related to geometric non-normality. 
In Subsection \ref{ss-rel-bet-inv}, we summarise relations we have established.

\subsection{Geometric generic embedding dimension}\label{ss-gged}

In this subsection, 
we recall definition and some properties of geometric generic embedding dimension 
introduced by Schr\"{o}er.

\begin{dfn}\label{d-gedim}
Let $k$ be a field of characteristic $p>0$ and 
let $X$ be a variety over $k$. 
We set $\gedim (X/k)$ to be the embedding dimension of $K(X) \otimes_k k^{1/p^{\infty}}$. 
We call $\gedim (X/k)$ the {\em geometric generic embedding dimension} of $X/k$. 
\end{dfn}

\begin{rem}\label{r-gedim}
The geometric generic embedding dimension is introduced by Schr\"{o}er in 
\cite[immediately before Theorem 2.3]{Sch10}, 
although he does not use the notation $\gedim (X/k)$. 
We here summarise some results established there. 
\begin{enumerate}
\item $\gedim (X/k)=\edim (K(X) \otimes_k k^{1/p^{\infty}})=
\edim (K(X) \otimes_k k^{1/p})$ (\cite[Proposiiton 2.1]{Sch10}), 
where $\edim(R)$ denotes the embedding dimension of a local ring $R$. 
\item If $X$ is not geometrically reduced over $k$, then it holds that 
\[
\gedim (X/k) < \log_p [k:k^p], 
\]
where $\log_p [k:k^p]$ is called the degree of imperfection of $k$ 
(\cite[Theorem 2.3]{Sch10}).
\end{enumerate}
\end{rem}

\subsection{Geometric non-reducedness}\label{ss-g-non-red}

\subsubsection{$\epsilon(X/k)$ vs $\gedim (X/k)$}

\begin{prop}\label{p-edim-vs-epsilon}
Let $k$ be a field of characteristic $p>0$ and 
let $X$ be a variety over $k$ such that $k$ is algebraically closed in $K(X)$. 
Then the following hold. 
\begin{enumerate}
\item $\gedim (X/k)+1 \leq p^{\epsilon(X/k)}$. 
\item The equation $\gedim (X/k)+1 = p^{\epsilon(X/k)}$ holds if and only if $\m^2=0$, 
where $\m$ denotes the maximal ideal of $K(X) \otimes_k k^{1/p^{\infty}}$. 
\end{enumerate}
\end{prop}

\begin{proof}
Set 
$R:=K(X) \otimes_k k^{1/p^{\infty}}$. 
Note that $R$ is an artinian local $k^{1/p^{\infty}}$-algebra. 
Let $F$ be a coefficient field of $R$, i.e. $F$ is a subring of $R$ such that 
the natural composite ring homomorphism $F \to R \to R/\m$ is an isomorphism. 
The existence of such $F$ is guaranteed by \cite[Theorem 6.4]{Mat89}. 
We have that 
\begin{itemize}
\item $\gedim (X/k) = \dim_F (\m/\m^2)$ (Remark \ref{r-gedim}(1)) and  
\item $p^{\epsilon (X/k)}=\length_R R = \dim_F R$ (Lemma \ref{l-length-dim}, Definition \ref{d-epsilon}). 
\end{itemize}
Therefore, it holds that 
\[
\gedim (X/k)+1 = \dim_F (\m/\m^2)+1 \leq \dim_F \m+1= \dim_F R = p^{\epsilon (X/k)}. 
\]
Therefore, (1) and (2) hold.  
\end{proof}

\subsubsection{$\epsilon(X/k)$ vs $m_F(X/k)$}

\begin{prop}\label{p-epsilon-ell}
Let $k$ be a field of characteristic $p>0$ such that $1<[k:k^p]<\infty$. 
Let $X$ be a variety over $k$ 
such that $k$ is algebraically closed in $K(X)$. 
Then it holds that 
\[
\epsilon(X/k) \leq \epsilon((X \times_k k^{1/p})_{\red}/k^{1/p})+\log_p [k:k^p]-1. 
\]
\end{prop}

\begin{proof}
Recall that  $\epsilon(X/k)=0$ if and only 
if $X$ is geometrically reduced over $k$ (Remark \ref{r-def-epsilon}). 
Thus, if $X$ is geometrically reduced over $k$, then 
we obtain $\epsilon(X/k)=0$ and we are done. 

We now handle the case when $X$ is not geometrically reduced over $k$. 
It follows from Theorem \ref{t-schroer}(1) that 
there is an intermediate field $k \subset l \subset k^{1/p}$ such that 
$X \times_k l$ is an integral scheme and 
the algebraic closure of $k$ in $K(X \times_k l)$ coincides with $k^{1/p}$. 
Then it holds that $K(X \times_k l)=K((X \times_k k^{1/p})_{\red})$, 
i.e. $X \times_k l$ is birational to $(X \times_k k^{1/p})_{\red}$ 
(Proposition \ref{p-bc-comp-field}). 
%Moreover, Theorem \ref{t-schroer}(2) implies that $[k^{1/p}:l]$
We have that 
\begin{eqnarray*}
\epsilon(X/k) 
&=& \epsilon((X \times_k l)^N/k^{1/p})+\log_p [k^{1/p}:l]\\
&=& \epsilon((X \times_k k^{1/p})_{\red}/k^{1/p})+\log_p [k^{1/p}:l], 
\end{eqnarray*}
where the first equality holds by Proposition \ref{p-epsilon-compute}(1) 
and the second one follows from Proposition \ref{p-epsilon-birat}. 
The assertion follows from $[k^{1/p}:l] < [k^{1/p}:k]=[k:k^p]$, 
where the inequality holds by $[k:k^p] \neq 1$, 
\end{proof}

\begin{thm}\label{t-epsilon-m}
Let $k$ be a field of characteristic $p>0$ such that $[k:k^p]<\infty$. 
Let $X$ be a variety over $k$ 
such that $k$ is algebraically closed in $K(X)$. 
Then it holds that 
\[
\epsilon(X/k) \leq  m_F(X/k) (\log_p [k:k^p]-1). 
\]
\end{thm}

\begin{proof}
If $X$ is geometrically reduced, 
then it holds that $\epsilon(X/k) =  m_F(X/k)=0$ 
(Remark \ref{r-mF1}, Remark \ref{r-def-epsilon}). 

From now on, we treat the case when $X$ is not geometrically reduced. 
In particular, $[k:k^p]>1$. 
For any $n \in \Z_{\geq 0}$, we set $X_n:=(X \times_k k^{1/p^n})^N_{\red}$. 
For any positive integer $n$, 
it follows from Proposition \ref{p-epsilon-ell} that 
\[
\epsilon(X_{n-1}/k^{1/p^{n-1}}) \leq \epsilon(X_n/k^{1/p^n})+(\log_p [k:k^p]-1). 
\]
Set $m:=m_F(X/k)$. 
Since $X_m$ is geometrically reduced over $k^{1/p^m}$ (Definition \ref{d-mF}), 
we have $\epsilon(X_m/k^{1/p^m})=0$ (Remark \ref{r-def-epsilon}). 
By induction, we have that 
\[
\epsilon(X/k)=\epsilon(X_0/k) \leq 
\epsilon(X_m/k^{1/p^m})+m (\log_p [k:k^p]-1)
\]
\[
=m_F(X/k) (\log_p [k:k^p]-1),
\]
as desired. 
\end{proof}

\subsubsection{Characterisation of geometric reducedness}

\begin{prop}\label{p-char-gred}
Let $k$ be a field of characteristic $p>0$. 
Let $X$ be a variety over $k$ 
such that $k$ is algebraically closed in $K(X)$. 
Then the following are equivalent. 
\begin{enumerate}
\item $X$ is geometrically reduced over $k$. 
\item $\gedim (X/k)=0$. 
\item $\epsilon (X/k)=0$. 
\item $m_F(X/k)=0$. 
\end{enumerate}
\end{prop}

\begin{proof}
By Definition \ref{d-gedim}, (1) and (2) are equivalent. 
%Since $\gedim (X/k)+1 \leq p^{\epsilon(X/k)}$ (Proposition \ref{p-edim-vs-epsion}), (2) implies (3). 
By Remark \ref{r-mF1} and Remark \ref{r-def-epsilon}, 
(1), (3), and (4) are equivalent.
\end{proof}

\subsection{Geometric non-normality}\label{ss-g-non-norm}

\subsubsection{$\gamma(X/k)$ vs $\ell_F(X/k)$}

\begin{prop}\label{p-ell-gamma}
Let $k$ be a field of characteristic $p>0$. 
Let $X$ be a normal variety over $k$ such that $k$ is algebraically closed in $K(X)$. 
Then it holds that $\gamma(X/k) \geq \ell_F(X/k)$. 
\end{prop}

\begin{proof}
The assertion follows from definitions of $\gamma(X/k)$ (Definition \ref{d-gamma}) 
and $\ell_F(X/k)$ (Definition \ref{d-F-length}). 
\end{proof}

\begin{prop}\label{p-ell-gamma2}
Let $k$ be a field of characteristic $p>0$ such that $[k:k^p]=p$. 
Let $X$ be a normal variety over $k$ such that $k$ is algebraically closed in $K(X)$. 
Then it holds that $\gamma(X/k) = \ell_F(X/k)$. 
\end{prop}

\begin{proof}
The assertion follows from the fact that $k^{1/p}$ is 
the unique purely inseparable extension of $k$ of degree $p$. 
\end{proof}

\subsubsection{Characterisation of geometric normality}

\begin{prop}\label{p-char-gnor}
Let $k$ be a field of characteristic $p>0$. 
Let $X$ be a normal variety over $k$ 
such that $k$ is algebraically closed in $K(X)$. 
Then the following are equivalent. 
\begin{enumerate}
\item $X$ is geometrically normal over $k$. 
\item $\gamma (X/k)=0$. 
\item $\ell_F(X/k)=0$. 
\end{enumerate}
\end{prop}

\begin{proof}
The assertion follows from Proposition \ref{p-def-gamma} and Remark \ref{r-F-length2}. 
\end{proof}

\subsection{Summary of relations}\label{ss-rel-bet-inv}

In this subsection, we summarise the relations between the invariants.  
Let $k$ be a field of characteristic $p>0$. 
Let $X$ be a variety over $k$ 
such that $k$ is algebraically closed in $K(X)$. 
Then the following hold. 
\begin{enumerate}
\item $\gedim (X/k)+1 \leq p^{\epsilon(X/k)}$ (Proposition \ref{p-edim-vs-epsilon}(1)). 
\item If $[k:k^p] <\infty$, then the following holds (Theorem \ref{t-epsilon-m}): 
\[
\epsilon(X/k) \leq  m_F(X/k) (\log_p [k:k^p]-1).
\]
\item 
If $X$ is normal, then 
$m_F(X/k) \leq \ell_F(X/k) \leq \gamma(X/k)$ 
(Remark \ref{r-m-ell}, Proposition \ref{p-ell-gamma}). 
\end{enumerate}

\section{Examples}

In this section, 
we compute our invariants for some explicit examples. 
In Subsection \ref{ss2-ex}, we introduce complete intersections 
obtained by Fermat-like hypersurfaces. 
We then consider curves of genus zero (Subsection \ref{ss3-ex}) 
and genus one (Subsection \ref{ss4-ex}).

\subsection{$q$-Fermat complete intersections}\label{ss2-ex}

On study of fibrations in positive characteristic, 
a typical example is as follows:
\[
V:=\{s_0x_0^2+s_1x_1^2+s_2x_2^2=0\} \subset \mathbb P^2_x \times \mathbb A^3_s,
\]
where $\mathbb P^2_x:=\mathbb P^2_{\F}$ and $\mathbb A^3_s:=\mathbb A^3_{\F}$ 
denote the projective plane and the three-dimensional affine space over an algebraically closed field $\F$ of characterstic two 
whose coordinates are $x:=(x_0, x_1, x_2)$ and $s:=(s_0, s_1, s_2)$, respectively. 
This is a variant of wild conic bundles, 
i.e. if the affine coordinate $s$ is replaced by a homogeneous coordinate, 
then it becomes a wild conic bundle (cf. \cite{MS03}).  
Since $V$ is smooth over $\F$, 
its generic fibre $X:=\Proj\,k[x_0, x_1, x_2]/(s_0x_0^2+s_1x_1^2+s_2x_2^2)$ is regular, 
where $k:=\F(s_0, s_1, s_2)$. 
On the other hand, $X$ is not geometrically reduced over $k$. 
More generally, for an algebraically closed field $\F$ of characteristic $p>0$ and 
the rational function field $k:=\F(s_0, ..., s_N)$ over $\F$ of degree $N+1$, 
\[
\Proj\,k[x_0, ..., x_N]/(s_0x_0^p+\cdots+s_Nx_N^p)
\]
is a regular projective variety over $k$ which is not geometrically reduced. 
In \cite[Section 3]{Sch10}, such varieties are called $p$-Fermat hypersurfaces. 
In this subsection, we study slightly generalised examples: 
$q$-Fermat complete intersections. 
One of the main results is Theorem \ref{t-pFCI-main}, 
which assures that if $r$ is a positive integer, then 
there exists an example $(k, X)$ 
with $\ell_F(X/k)=m_F(X/k)=1$ and $\gamma(X/k)=\epsilon(X/k)=r$. 

\begin{dfn}\label{d-pFCI}
Let $k$ be a field of characteristic $p>0$. 
Fix a non-negative integer $e$ and set $q:=p^e$. 
Let $r$ and $N$ be integers such that $0<  r < N$. 
Assume that we have a subset $S:=\{ s_{ij}\,|\, 1 \leq i \leq r, 0 \leq j \leq N\} \subset k$ 
such that 
\[
[k^p(S):k^p] = p^{|S|} = p^{r(N+1)}, 
\]
where $k^p(S)$ denotes the minimum subfield of $k$ 
that contains $k^p \cup S$. 
We set 
\begin{eqnarray*}
f_1&:=& s_{10}x_0^q+\cdots+s_{1N}x_N^q\\
f_2&:=& s_{20}x_0^q+\cdots+s_{2N}x_N^q\\
&\cdots&\\
f_r&:=& s_{r0}x_0^q+\cdots+s_{rN}x_N^q\\
\end{eqnarray*}
and  
\[
X:=\Proj\,k[x_0, ..., x_N]/(f_1, ..., f_r).  
\]
Then $X$ is called a $q$-{\em Fermat complete intersection} 
({\em in} $\mathbb P^N_k$ {\em of codimension} $r$). 
\end{dfn}

\begin{rem}\label{r-pFCI}
For convenience, we call $\mathbb P^N_k$ 
a $q$-{\em Fermat complete intersection} ({\em in} $\mathbb P^N_k$) 
of codimension zero, which is corresponding to the case when $r=0$. 
This notion is useful when we apply induction on $r$. 
\end{rem}

\begin{rem}\label{r-pFCI-Fp-model}
We use the same notation as in Definition \ref{d-pFCI}. 
Set $k_0:=\F_p(S)$ and 
\[
X_0:=\Proj\,k_0[x_0, ..., x_N]/(f_1, ..., f_r).  
\]
The condition $[k^p(S):k^p] = p^{|S|}$ implies that $[k^p_0(S):k^p_0]=p^{|S|}$. 
This implies that ${\rm tr.deg}_{\F_p}\,k_0=|S|$. 
Since $X_0 \times_{k_0} k =X$, some problems are reduced to the case when $k=k_0$. 
\end{rem}

\begin{lem}\label{l-pFCI1}
We use the same notation as in Definition \ref{d-pFCI}. 
Then the following hold. 
\begin{enumerate}
\item $\dim X=N-r$. 
\item 
$X$ is a regular integral scheme. 
\item 
For any integer $d$ with $0 \leq d \leq e$, 
$(X \times_k k^{1/p^d})_{\red}$ is $k^{1/p^d}$-isomorphic to a 
$p^{e-d}$-Fermat complete intersection in $\mathbb P^N_{k^{1/p^d}}$ 
of codimension $r$. 
\item 
$(X \times_k k^{1/q})_{\red} \simeq \mathbb P^{N-r}_{k^{1/q}}$. 
\item 
$H^0(X, \MO_X)=k$. 
\item 
%For any $i$, 
$D_+(x_0) \cap \cdots \cap D_+(x_N) \neq \emptyset$, where each $D_+(x_i)$ denotes the open subset of $X$ 
that is equal to the non-vanishing locus of $x_i$ (cf. \cite[Ch. II, Proposition 2.5]{Har77}). 
\end{enumerate} 
\end{lem}

\begin{proof}
Let us show (1). Consider the coefficient matrix: 
\[
\left(
\begin{array}{c}
f_1\\
\vdots\\
f_r
\end{array}
\right) = A
\left(
\begin{array}{c}
x_0^q\\
\vdots\\
x_r^q
\end{array}
\right), \qquad A:=\left(
\begin{array}{ccc}
s_{10} & \cdots & s_{1N}\\
\vdots & \ddots & \vdots\\
s_{r0} & \cdots & s_{rN}\\
\end{array}
\right).
\]
It is enough to prove that $\rank\,A = r$ 
(indeed, applying elementary row operations, $A$ becomes a matrix in reduced row echelon form). 
It suffices to prove this equation as a matrix over $k_0:=\F_p(S)=\F_p(\{s_{ij}\}_{i,j})$. 
This follows from the fact that ${\rm tr.deg}_{\F_p} k_0=|S|$ (Remark \ref{r-pFCI-Fp-model}). 
Thus, (1) holds.

Let us show (2). 
By (1), $X$ is a complete intersection in $\mathbb P^N_k$. 
In particular, $X$ is connected. 
Hence, it suffices to show that $X$ is regular. 
By symmetry, it is enough to prove that the affine open subset 
\[
D_+(x_0) \simeq \Spec\,k[x_1, ..., x_n]/(f'_1, ..., f'_r)
\]
of $X$ is regular, where 
\[
f'_i:=s_{i0}+s_{i1}x_1^p+\cdots+s_{iN}x_N^p.
\]
By \cite[Theorem 26.5]{Mat89} and $[k^p(S):k^p]=p^{|S|}$, 
there is an $\F_p$-derivation  
\[
\partial_{s_{ij}}:k[x_1, ..., x_N] \to k[x_1, .., x_N]
\]
such that $\partial_{s_{ij}}(s_{ij})=1, \partial_{s_{ij}}(x_{\ell})=0, \partial_{s_{ij}}(s_{i'j'})=0$ 
for any $\ell$ and $(i', j') \neq (i, j)$. 
We now apply Jacobian criterion \cite[Proposition 22.6.7(iii)]{EGAIV1} 
to $D_+(x_0)$. 
Consider the following minor of the Jacobian matrix:  
\[
J:=
\left(
\begin{array}{cccc}
\partial_{s_{10}} f_1 & \partial_{s_{20}}f_1 & \cdots & \partial_{s_{r0}} f_1\\
\partial_{s_{10}} f_2 & \partial_{s_{20}}f_2 & \cdots & \partial_{s_{r0}} f_2\\
\vdots & \vdots & \ddots & \vdots\\
\partial_{s_{10}} f_r & \partial_{s_{20}}f_r & \cdots & \partial_{s_{r0}} f_r\\
\end{array}
\right)
=
\left(
\begin{array}{cccc}
1 & 0 & \cdots & 0\\
0 & 1 & \cdots & 0\\
\vdots & \vdots & \ddots & \vdots\\
0 & 0 & \cdots & 1\\
\end{array}
\right)
\]
Since $\rank J=r$, the Jacobian criterion implies that $D_+(x_0)$ is regular. 
Thus, (2) holds.

Let us show (3). 
For 
\begin{eqnarray*}
g_i&:=& s^{1/p^d}_{i0}x_0^{p^{e-d}}+\cdots +s^{1/p^d}_{iN}x_N^{p^{e-d}},
\end{eqnarray*}
we set 
\[
Y:=\Proj\,k^{1/p^d}[x_0, ..., x_N]/(g_1, ..., g_r).  
\]
For $S^{1/p^d}:=\{ s^{1/p^d}\,|\, s\in S\}$, we obtain  
\[
[(k^{1/p^d})^p(S^{1/p^d}):(k^{1/p^d})^p]=[k^p(S):k^p]=p^{|S|}. 
\]
Therefore, $Y$ is a $p^{e-d}$-Fermat complete intersection in $\mathbb P^N_{k^{1/p^d}}$ of codimension $r$. 
By (2), $Y$ is reduced, hence it holds that $(X \times_k k^{1/p^d})_{\red} \simeq Y$. 
Thus (3) holds. 
The assertion (4) follows from (1) and (3).

Let us show (5). 
Let $Z$ be a complete intersection in $\mathbb P^N_k$ of codimension $c \leq N-1$. 
By induction on codimension $c$, 
we can check that $H^0(Z, \MO_Z)=k$ and $H^i(Z, \MO_Z(m))=0$ 
for any $m \in \Z$ and $0<i < N-c$. 
Since $\dim X=N-r \geq 1$, we are done. 

Let us show (6). 
Since $X \neq \emptyset$, 
there exists $j$ such that $X \cap D_+(x_j) \neq \emptyset$. 
By symmetry, we have $X \cap D_+(x_i) \neq \emptyset$ for any $i$. 
Since $X$ is irreducible by (2), the assertion (6) holds. 
\end{proof}

\begin{prop}\label{p-pFCI2}
We use the same notation as in Definition \ref{d-pFCI}. 
Assume $e>0$. 
Set 
\begin{itemize}
\item $l:=k(s_{11}^{1/q}, ..., s_{1N}^{1/q})$, and 
\item $k':=k(s_{10}^{1/q}, s_{11}^{1/q}, ..., s_{1N}^{1/q})$. 
\end{itemize}
Then the following hold. 
\begin{enumerate}
\item 
There is an $l$-linear ring isomorphism
\[
\sigma:l[x_0, ..., x_N] \to l[y_0, ..., y_N]
\] 
preserving degrees, $t_{10} \in l,$ and 
$T:=\{t_{ij}\,|\, 2 \leq i \leq r, 1 \leq j \leq N\} \subset l$ 
such that if 
\[
\widetilde{\sigma}:\Proj\,l[y_0, ..., y_N] \to
\Proj\,l[x_0, ..., x_N]
\]
is the $l$-isomorphism induced by $\sigma$ 
and we define $g_1, ..., g_r$ by the following equations 
\begin{equation}\label{e1-pFCI2}
\begin{split}
g_1 &:= t_{10}y_0^q+y_1^q\\
g_2 &:= t_{21}y_1^q+\cdots +t_{2N}y_N^q\\
&\cdots \\
g_r &:= t_{r1}y_1^q+\cdots +t_{rN}y_N^q,\\
\end{split}
\end{equation}
then 
\begin{enumerate}
\item 
$\sigma(x_0)=y_0$, $\sigma(x_N)=y_N$, 
\item 
$\sigma((f_1, ..., f_r))=(g_1, ..., g_r)$, and 
\item 
$[l^p(\{t_{10}\} \cup T):l^p]=p^{|T|+1}$. 
\end{enumerate}
\item For $R:=l[y_0, ..., y_N]/(g_1, ..., g_r)$, 
$R$ is an integral domain and $X \times_k l$ is $l$-isomorphic to $\Proj\,R$. 
In particular, $X \times_k l$ is an integral scheme. 
\item 
$X \times_k l$ is not normal. 
%\item 
%For the normalisation of $X' \to X \times_k l$, $H^0(X', \MO_{X'})=k'$.  
\item 
For the normalisation $X' \to X \times_k l$ of $X \times_k l$,
$X'$ is $k'$-isomorphic to a $q$-Fermat complete intersection 
in $\mathbb P^{N-1}_{k'}$ of codimension $r-1$. 
\end{enumerate}
\end{prop}

\begin{proof}
Let us show (1). 
We define an $l$-algebra homomorphism $\sigma$ by 
\begin{equation}\label{e2-pFCI2}
\begin{split}
\tau:l[y_0, ..., y_N] &\to l[x_0, ..., x_N]\\
y_i&\mapsto  x_i \quad \text{if}\quad 0 \leq i \leq N\,\,{\rm and}\,\,i \neq 1\\
y_1 &\mapsto  x_1+s_{11}^{-1/q}s_{12}^{1/q}x_2\cdots +s_{11}^{-1/q}s_{1N}^{1/q}x_N.\\
\end{split}
\end{equation}
Then $\tau$ is an $l$-algebra isomorphism preserving degrees, 
hence $\sigma:=\tau^{-1}$ is an $l$-algebra isomorphism preserving degrees. 
It holds that (\ref{e2-pFCI2}) implies (a). 
We define $t_{ij}$ as follows. 
\begin{enumerate}
\renewcommand{\labelenumi}{(\roman{enumi})}
\item $t_{10}:=s_{11}^{-1} s_{10}$. 
\item $t_{i1}:=s_{i1}-s_{i0}t_{10}^{-1}$ for any integer $i$ such that $2 \leq i \leq r$. 
\item $t_{ij}:=s_{ij}-s_{i1}s_{11}^{-1}s_{1j}$ 
for any integers $i, j$ such that $2 \leq i \leq r, 2 \leq j \leq N$. 
\end{enumerate}
For $g_1, ..., g_r$ defined as in (\ref{e1-pFCI2}), 
we now show $\sigma((f_1, ..., f_r))=(g_1, ..., g_r)$. We have 
\begin{eqnarray*}
\sigma\left(s_{11}^{-1}f_1\right)
&=&\sigma\left(s_{11}^{-1}( s_{10}x_0^q+s_{11}x_1^q+s_{12}x_2^q+\cdots +s_{1N}x_N^q)\right)\\
&=&\sigma\left(s_{11}^{-1} s_{10}x_0^q+x_1^q+s_{11}^{-1}s_{12}x_2^q\cdots +s_{11}^{-1}s_{1N}x_N^q\right)\\
&=&\sigma\left(s_{11}^{-1} s_{10}x_0^q+(x_1+s_{11}^{-1/q}s_{12}^{1/q}x_2\cdots +s_{11}^{-1/q}s_{1N}^{1/q}x_N)^q\right)\\
&=&t_{10}y_0^q+y_1^q\\
&=& g_1. 
\end{eqnarray*}
For any integer $i$ such that $2 \leq i \leq r$, it holds that 
\begin{eqnarray*}
&&\sigma\left( f_i\right)\\
&=&\sigma\left( s_{i0}x_0^q+s_{i1}x_1^q+s_{i2}x_2^q+\cdots+s_{iN}x_N^q\right)\\
&=&\sigma\left( s_{i0}x_0^q+s_{i1}(\tau(y_1)-s_{11}^{-1/q}s_{12}^{1/q}x_2-\cdots -s_{11}^{-1/q}s_{1N}^{1/q}x_N)^q+ s_{i2}x_2^q+\cdots +s_{iN}x_N^q\right)\\
&=&\sigma\left( s_{i0}x_0^q+s_{i1}(\tau(y_1)^q-s_{11}^{-1}s_{12}x^q_2-\cdots -s_{11}^{-1}s_{1N}x_N^q)+ s_{i2}x_2^q+\cdots +s_{iN}x_N^p\right)\\
&=&\sigma\left( s_{i0}x_0^q+s_{i1}\tau(y_1)^q+(s_{i2}-s_{i1}s_{11}^{-1}s_{12})x^q_2+ \cdots 
+(s_{iN}-s_{i1}s_{11}^{-1}s_{1N})x^q_N\right)\\
&=&s_{i0}y_0^q+(t_{i1}+s_{i0}t_{10}^{-1})y_1^q+t_{i2}y^q_2+ \cdots 
+t_{iN}y^q_N\\
&=&s_{i0}t_{10}^{-1}g_1 +g_i. 
\end{eqnarray*}
Therefore, it holds that $\sigma((f_1, ..., f_r))=(g_1, ..., g_r)$. 
Thus (b) holds. 

In order to complete the proof of (1), 
it is enough to prove (c), i.e. $[l^p(\{t_{10}\} \cup T):l^p]=p^{|T|+1}$. 
It is clear that $[l^p(\{t_{10}\} \cup T):l^p] \leq p^{|T|+1}$. 
Thus it suffices to show $[l^p(\{t_{10}\} \cup T):l^p] \geq p^{|T|+1}$. 
Set $S':=\{s_{11}, ..., s_{1N}\}$ and $S'':=\{s_{20}, ..., s_{r0}\}$.  
%we set $|S_1|=N+r-1$. 
We have that $l^p=k^p(S')$. 

We now show that 
\begin{equation}\label{e3-pFCI2}
\kappa_1:=k^p(S' \cup S'' \cup \{t_{10}\} \cup T) = k^p(S)=:\kappa_2. 
\end{equation}
The inclusion $\kappa_1  \subset \kappa_2$ is obvious. 
Let us show $\kappa_1  \supset \kappa_2$. 
By $t_{10}=s_{11}^{-1} s_{10}$, $t_{10} \in \kappa_1$, and 
$s_11 \in S' \subset \kappa_1$, we obtain $s_{10} \in \kappa_1$. 
For any integers $i,j$ such that $2 \leq i \leq r$ and $2 \leq j \leq N$, 
it follows from (ii), (iii), and $S' \cup S''$ that $s_{i1}, s_{ij} \in \kappa_1$. 
Therefore, we obtain $\kappa \subset \kappa_2$. 
This completes the proof of (\ref{e3-pFCI2}). 

It holds that 
\[
p^{|S|}=[k^p(S):k^p]=[k^p(S):k^p(S' \cup S'')][k^p(S' \cup S''):k^p]
\]
\[
\leq [l^p(\{t_{10}\} \cup T):l^p][k^p(S' \cup S''):k^p] \leq [l^p(\{t_{10}\} \cup T):l^p]\cdot p^{|S' \cup S''|},
\]
where the first inequality follows from 
$k^p(S' \cup S'') \supset k^p(S')= l^p$ and 
$k^p(S)=(k^p(S' \cup S''))(\{t_{10}\} \cup T)$. 
Therefore, in order to prove $[l^p(\{t_{10}\} \cup T):l^p] \geq p^{|T|+1}$, 
it is sufficient to show that $|S|-|S' \cup S''| = |T|+1$. 
This equation holds by the following computation: 
\[
r(N+1)=|S|=|S' \cup S'' \cup \{t_{10}\} \cup T|
\leq |S' \cup S''|+| \{t_{10}\}|+|T| 
\]
\[
\leq (N+(r-1))+1+N(r-1)=r(N+1). 
\]
Thus, (1) holds.

Let us show (2). 
By (1),  $X \times_k l$ is $l$-isomorphic to $\Proj\,R$. 
Then it follows from Lemma \ref{l-pFCI1}(2) that $\Proj\,R$ is irreducible. 
By \cite[Proposition 3.1.14(ii)]{Gro61}, 
it is enough to show that $R$ and $\Proj\,R$ are reduced. 
Both the proofs are the same, hence we only show that $R$ is reduced. 
Since $\Spec\,R$ is a complete intersection in $\mathbb A^{N+1}_k$, $R$ is Cohen--Macaulay. 
Therefore, it is enough to prove that $R$ is regular in codimension zero, 
i.e. to find a regular point of $\Spec\,R$. 
Consider the following minor of a Jacobian matrix: 
\[
\left(
\begin{array}{cccc}
\partial_{t_{10}} g_1 & \partial_{t_{2N}}g_1 & \cdots & \partial_{t_{rN}} g_1\\
\partial_{t_{10}} g_2 & \partial_{t_{2N}}g_2 & \cdots & \partial_{t_{rN}} g_2\\
\vdots & \vdots & \ddots & \vdots\\
\partial_{t_{10}} g_r & \partial_{t_{2N}}g_r & \cdots & \partial_{t_{rN}} g_r\\
\end{array}
\right)
=
\left(
\begin{array}{cccc}
y_0^q & 0 & \cdots & 0\\
0 & y_N^q & \cdots & 0\\
\vdots & \vdots & \ddots & \vdots\\
0 & 0 & \cdots & y_N^q\\
\end{array}
\right).
\]
By (a) of (1) and Lemma \ref{l-pFCI1}(6), 
there exists a regular point of $\Spec\,R$. 
Thus, (2) holds.

Let us show (3). 
Assuming that $X \times_k l$ is normal, 
let us derive a contradiction. 
If $l_0:=\F_p(\{t_{10}\} \cup T)$ and $Y_0:=l_0[y_0, ..., y_N]/(g_1, ..., g_r)$, 
then we have $X \times_k l \simeq Y_0 \times_{l_0} l$. 
The usual Jacobian criterion shows that 
the hyperplane section $H:=\{y_0=0\}$ on $Y_0$ 
is contained in the non-regular locus of $Y_0$. 
Since $\dim H \geq \dim Y_0 -1$, $Y_0$ is not normal. 
Then the faithfully flat descent for normality 
implies that $X \times_k l$ is not normal. 
Thus, (3) holds.

Let us show (4). 
We have 
\[
R=l[y_0, ..., y_N]/(t_{10}y_0^q+y_1^q, g_2, ..., g_r). 
\]
We set 
\[
R':=k'[y_1, ..., y_N]/(g_2, ..., g_r) = l(t_{10}^{1/q})[y_1, ..., y_N]/(g_2, ..., g_r). 
\]
We consider an $l$-algebra homomorphism of graded $l$-algebras preserving degrees: 
\begin{eqnarray*}
\varphi:R &\to& R'\\
y_0 &\mapsto& -t_{10}^{-1/q}y_1\\
y_i &\mapsto& y_i\qquad \text{if}\quad 1 \leq i \leq r. 
\end{eqnarray*}
Note that $X':=\Proj\,R'$ is a $q$-Fermat complete intersection 
in $\mathbb P^{N-1}_{k'}$ of codimension $r-1$. 
By Lemma \ref{l-pFCI1}(2), $\Proj\,R'$ is normal. 
It holds that $R'$ is a finitely generated $R$-module. 
Since $R'$ is an integral domain and $\dim R=\dim R'$, 
$\varphi$ is injective. 
Then, by the construction, 
$\varphi$ induces a finite birational morphism $\Proj\,R' \to \Proj\,R$. 
Therefore, $\varphi$ is the normalisation of $X \times_k l$. 
Thus, (4) holds. 
\end{proof}

\begin{cor}\label{c-pFCI-elem}
Let $k$ be a field of characteristic $p>0$. 
Fix a positive integer $e$ and set $q:=p^e$. 
Let $r$ and $N$ be integers such that $0 <  r < N$. 
Let $X$ be a $q$-Fermat complete intersection in $\mathbb P^N_k$ 
of codimension $r$. 
Then there exists an elemental extension $(k \subset l \subset k', \varphi:X' \to X)$ 
of $(k, X)$ 
such that 
\begin{enumerate}
\item $[k':l]=q$, $[l:k]=q^N$, 
\item $X \times_k l$ is not normal, and 
\item $X'$ is $k'$-isomorphic to 
a $q$-Fermat complete intersection in $\mathbb P^{N-1}_{k'}$ of codimension $r-1$ 
(note that if $r=1$, this assertion means that $X'$ is $k'$-isomorphic to $\mathbb P^{N-1}_{k'}$). 
\end{enumerate}
\end{cor}

\begin{proof}
The assertion follows from Proposition \ref{p-pFCI2}. 
\end{proof}

\begin{thm}\label{t-pFCI-main}
Let $k$ be a field of characteristic $p>0$. 
Fix a positive integer $e$ and set $q:=p^e$. 
Let $r$ and $N$ be integers such that $0 \leq  r < N$. 
Let $X$ be a $q$-Fermat complete intersection in $\mathbb P^N_k$ 
of codimension $r$. 
Then the following hold. 
\begin{enumerate}
\item 
If $r>0$, then $\ell_F(X/k)=m_F(X/k)=e$. 
\item 
$\epsilon(X/k)=er$ and $\gamma(X/k) \geq r$. 
\item 
If $e=1$, then $\gamma(X/k)=r$. 
\end{enumerate}
\end{thm}

\begin{proof}
The assertion (1) follows from Lemma \ref{l-pFCI1}(2)(3). 
We prove (2) by induction on $r$. 
If $r=0$, then we have $X \simeq \mathbb P^{N-r}_k$ (Remark \ref{r-pFCI}), 
hence the assertion follows from Proposition \ref{p-char-gred} and Proposition \ref{p-char-gnor}. 
Thus, we may assume that $r>0$. 
Applying Corollary \ref{c-pFCI-elem}, 
there exists an elemental extension $(k \subset l \subset k', \varphi:X' \to X)$ 
of $(k, X)$ 
such that $[k':l]=q$ and 
$X'$ is $k'$-isomorphic to 
a $q$-Fermat complete intersection in $\mathbb P^{N-1}_{k'}$ of codimension $r-1$. 
By the induction hypothesis, 
it holds that  $\epsilon(X'/k')=e(r-1)$ and $\gamma(X'/k') \geq r-1$. 
Then Proposition \ref{p-epsilon-compute} implies 
$\epsilon(X/k)=\epsilon(X'/k')+\log_p[k':l]=er$. 
By Definition \ref{d-gamma}, we have $\gamma(X/k) \geq \gamma(X'/k')+1=r$. 
Thus, (2) holds.

Let us show (3). 
By (2), it is enough to show that $\gamma(X/k) \leq r$. 
Set $Y:=(X \times_k k^{1/p^{\infty}})_{\red}$. It holds that 
\[
Y = \Proj\,k^{1/p^{\infty}}[x_0, ..., x_N]
/(f'_1, ..., f'_r) \simeq \mathbb P^{N-r}_{k^{1/p^{\infty}}}
\]
with $f_i=f'^p_i$. 
For the induced morphism $f:Y \to X$, we have that 
\[
K_Y+(p-1) \sum_{i=1}^{\gamma(X/k)} D_i \sim f^*K_X
\]
for some nonzero effective divisors $D_1, ..., D_{\gamma(X/k)}$ on $Y$ (Theorem \ref{t-cano-gamma}). 

Fix 
$H \in |\MO_{\mathbb P_k^N}(1)|$, 
$H_X \in |\MO_{\mathbb P_k^N}(1)|_X|$, and 
$H_Y \in |\MO_{\mathbb P_{k^{1/p}}^N}(1)|_Y|$. 
We have $H_Y \sim f^*H_X$. 
It holds that  
\begin{eqnarray*}
\left(K_Y+(p-1) \sum_{i=1}^{\gamma(X/k)} D_i\right) \cdot H_Y^{N-r-1} 
&\geq&K_Y\cdot H_Y^{N-r-1} + (p-1)\gamma(X/k) \\
&=& -(N-r+1)+(p-1) \gamma(X/k).  
\end{eqnarray*}
Since $\epsilon(X/k)=r$, we have that  
\begin{eqnarray*}
f^*K_X \cdot H_Y^{N-r-1} 
&=& \frac{1}{p^{\epsilon(X/k)}} K_X \cdot H_X^{N-r-1}\\
&=& \frac{1}{p^r}  (K_{\mathbb P_k^N}+rpH) \cdot (pH)^r \cdot H^{N-r-1}\\
&=& -(N+1)+rp, 
\end{eqnarray*}
where the first equality holds by \cite[Example 1 in page 299]{Kle66}. 
To summarise, it holds that 
\begin{eqnarray*}
-(N+1)+rp &=& f^*K_X \cdot H_Y^{N-r-1}  \\
&=& \left(K_Y+(p-1) \sum_{i=1}^{\gamma(X/k)} D_i\right) \cdot H_Y^{N-r-1} \\
&\geq&  -(N-r+1)+(p-1) \gamma(X/k). 
\end{eqnarray*}
Therefore, we get $r \geq \gamma(X/k)$, as desired.
\end{proof}

We close this subsection with the following result. 
This example shows that the conductor-like effective divisor 
$C$ that satisfies $K_Y+C \sim f^*K_X$ 
is not canonically determined. 

\begin{prop}\label{p-Fermat}
Let $k$ be a field of characteristic $p>0$. 
Let $s, t \in k$ be elements 
such that $[k^p(s, t):k^p]=p^2$. Set 
\[
X:=\Proj\,k[x, y, z]/(sx^p+ty^p+z^p). 
\]
Then the following hold. 
\begin{enumerate}
\item $X$ is a regular projective curve over $k$ such that $H^0(X, \MO_X)=k$. 
\item $\ell_F(X/k)=m_F(X/k)=\gamma(X/k)=\epsilon(X/k)=1.$
\end{enumerate}
Furthermore, given finitely many closed points $P_1, ..., P_n$ of $X$, 
there exist purely inseparable extensions $k \subset l \subset k'$ 
which satisfy the following properties. 
\begin{enumerate}
\setcounter{enumi}{2}
\item $[k':l]=[l:k]=p$. 
\item $X \times_k l$ is an integral scheme which is not normal. 
\item The normalisation $Y$ of $X \times_k l$ is isomorphic to $\mathbb P^1_{k'}$. 
\item The effective divisor $D$ on $Y$ defined by the conductor of 
the normalisation $\nu:Y \to X \times_k l$ can be written by 
$D=(p-1)Q$ for some $k'$-rational point $Q$ of $\mathbb P^1_{k'}$ and it holds that 
\[
K_Y+(p-1)Q \sim \varphi^*K_X 
\]
for the induced morphism $\varphi:Y \to X$. 
\item $Q$ is not contained in $\varphi^{-1}(\{P_1, ..., P_n\})$. 
\end{enumerate}
In particular, $(k \subset l \subset k', \varphi:Y \to X)$ is an elemental extension of $(k, X)$. 
\end{prop}

\begin{proof}
Let us show (1) and (2). 
Fix a rational function field $k':=k(u')$ over $k$ of degree one. 
Set $X':=X \times_k k'=\Proj\,k'[x, y, z]/(sx^p+ty^p+z^p)$. 
For $s':=u's$ and $t':=u't$, it holds that 
\[
X'=\Proj\,k'[x, y, z]/(sx^p+ty^p+z^p) \simeq \Proj\,k'[x, y, z]/(s'x^p+t'y^p+u'z^p). 
\]
Since $[k'^p(s', t', u'):k'^p]=p^3$, 
$X'$ is isomorphic to a $p$-Fermat complete intersection in $\mathbb P^2_{k'}$ of codimension one 
(cf. Definition \ref{d-pFCI}). 
Therefore, (1) holds by Lemma \ref{l-pFCI1}(2)(5). 
It follows from Theorem \ref{t-pFCI-main} that 
\[
\ell_F(X'/k')=m_F(X'/k')=\gamma(X'/k')=\epsilon(X'/k')=1.
\]
Since $\gamma(X/k) \geq 1$ (Proposition \ref{p-def-gamma}), 
we get 
$\ell_F(X/k)=m_F(X/k)=\gamma(X/k)=\epsilon(X/k)=1$ 
(Proposition \ref{p-gamma-sep-bc}, Proposition \ref{p-length-sep-bc}, Proposition \ref{p-length-sep-bc2}, Proposition \ref{p-epsilon-sep-bc}). 
Hence (2) holds.

Fix an arbitrary element $a \in k$ and set $l:=k((s-a^pt)^{1/p})$
(we shall prove that (3)--(6) hold for any $a \in k$ and 
that if $a$ is sufficiently general, then also (7) holds). 
We have 
\[
sx^p+ty^p+z^p = ((s-a^pt)^{1/p}x)^p+a^ptx^p+ty^p+z^p
\]
\[
=((s-a^pt)^{1/p}x+z)^p+t(ax+y)^p.
\]
For $x':=x$, $y':=ax+y$, and $z':=(s-a^pt)^{1/p}x+z$, we have that 
\[
X \times_k l = \Proj\,l[x, y, z]/(sx^p+ty^p+z^p) 
\simeq \Proj\,l[x', y', z']/(ty'^p+z'^p).
\]
Then (3)--(5) hold (cf. the proof of Proposition \ref{p-pFCI2}).  
We have 
\[
K_Y + (p-1)C \sim \varphi^*K_X
\]
for some nonzero effective divisor $C$ on $Y$. 
By computing the degrees of the both hand sides (cf. the proof of Theorem \ref{t-pFCI-main}), 
we obtain $\deg_{k^{1/p}} C=1$, i.e. $C=:Q$ is a single $k^{1/p}$-rational point. 
Thus, (6) holds.

Hence it is enough to prove (7). 
By Jacobian criterion for regularity, 
$[x':y':z']=[1:0:0]$ is the unique non-regular closed point of $X \times_k l$. 
Note that the image $[x:y:z]$ of this point to $X$ satisfies $x \neq 0$ and 
\[
[x:y:z]=[1:a:(s-a^pt)^{1/p}].
\]
Since $k$ is an imperfect field, $k$ is an infinite field. 
Hence, the set 
\[
\{[1:a:(s-a^pt)^{1/p}] \in X \times_k k^{1/p}\,|\, a\in k\}
\] 
is an infinite set. 
Therefore, there exists $a \in k$ such that (7) holds.  
\end{proof}

\subsection{Curves of genus zero}\label{ss3-ex}

The purpose of this subsection is 
to determine non-smooth regular curves of genus zero (Theorem \ref{t-conic}). 
%Although this topic might be well known, 
Although some of the arguments are extracted from the proof of \cite[Lemma 2.12]{Tan1}, 
we include a self-contained proof for the reader's convenience.

\begin{lem}\label{l-conic}
Let $k$ be a field of characteristic $p>0$. 
Take 
\[
f=ax^2+by^2+cz^2+\alpha yz+\beta zx+ \gamma xy \in k[x, y, z]. 
\]
Assume that $X:=\Proj\,k[x, y, z]/(f)$ is regular and of dimension one.  
Then one of the following holds. 
\begin{enumerate}
\item $X$ is smooth over $k$. 
\item $p=2$, $\alpha=\beta=\gamma=0$, $a \neq 0, b \neq 0, c \neq 0$, and $[k((a/c)^{1/2}, (b/c)^{1/2}):k]=4$. 
\end{enumerate}
\end{lem}

\begin{proof}
The scheme $X$ is covered by the following three affine open subsets: 
\begin{eqnarray*}
X_1 &:=& \Spec\,k[y, z]/(a+by^2+cz^2+\alpha yz+\beta z+ \gamma y)\\
X_2 &:=& \Spec\,k[x, z]/(ax^2+b+cz^2+\alpha z+\beta xz+ \gamma x)\\
X_3 &:=& \Spec\,k[x, y]/(ax^2+by^2+c+\alpha y+\beta x+ \gamma xy).\\
\end{eqnarray*}
The proof consists of four steps. 

\setcounter{step}{0}

\begin{step}\label{s1-conic}
If $p \neq 2$, then $X$ is smooth over $k$. 
\end{step}

\begin{proof}(of Step \ref{s1-conic}) 
Since $p \neq 2$, 
there is an orthogonal basis for any symmetric bilinear form. 
Therefore, we obtain 
\[
X \simeq \Proj\,k[x, y, z]/(a'x^2+b'y^2+c'z^2)
\] 
for some $a', b', c' \in k$. 
If $a' =0$, then it is easy to see that $\dim X \neq 1$ or $X$ is not regular. 
Hence, we have $a' \neq 0$. 
By symmetry, we obtain $b' \neq 0$ and $c' \neq 0$. 
It follows from the Jacobian criterion for smoothness 
that $X$ is smooth over $k$. 
This completes the proof of Step \ref{s1-conic}. 
\end{proof}

\begin{step}\label{s2-conic}
If $p=2$ and $\alpha \neq 0$, then $X_1$ is smooth over $k$. 
\end{step}

\begin{proof}(of Step \ref{s2-conic})
Passing to the base change to the separable closure, 
we may assume that $k$ is separably closed. 
By $p=2$ and $\alpha \neq 0$,  we have $by^2+\alpha yz + c z^2 = l_1 l_2$ 
for some $l_1:=\lambda_1 y+\mu_1 z, l_2:=\lambda_2 y+\mu_2 z \in k[y, z]$ 
such that $(l_1, l_2)=(y, z)$. 
The problem is reduced to the case when $b=c=0$. 
Then the defining equation of $X_1$ satisfies 
\[
a+by^2+cz^2+\alpha yz+\beta z+ \gamma y
=a+\alpha yz+\beta z+ \gamma y =\alpha (y+\zeta )(z+\xi ) +\eta
\]
for some $\zeta, \xi, \eta \in k$. 
Since $X$ is an integral scheme, so is $X_1$, hence $\eta \neq 0$. 
Then we have $X_1 \simeq \Spec\,k[y', z']/(y'z'+\eta')$ for some $\eta' \in k \setminus \{0\}$. 
Therefore, $X_1$ is smooth over $k$. 
This completes the proof of Step \ref{s2-conic}. 
\end{proof}

\begin{step}\label{s3-conic}
If $p=2$ and $(\alpha, \beta, \gamma) \neq (0, 0, 0)$, then $X$ is smooth over $k$. 
\end{step}

\begin{proof}(of Step \ref{s3-conic}) 
By symmetry, we may assume that $\alpha \neq 0$. 
It follows from Step \ref{s2-conic} that $X_1$ is smooth. 
Again by symmetry, it is enough to show that $X_2$ is smooth. 
Applying Step \ref{s2-conic} after replacing $(X_1, \alpha)$ by $(X_2, \beta)$, 
we may assume that $\beta =0$. 
By $p=2$, 
it follows from the Jacobian criterion for smoothness 
that $X_2$ is smooth. 
This completes the proof of Step \ref{s3-conic}. 
\end{proof}

\begin{step}\label{s4-conic}
If $p=2$ and $\alpha=\beta=\gamma=0$, then $a \neq 0, b \neq 0, c \neq 0$, and $[k((a/c)^{1/2}, (b/c)^{1/2}):k]=4$. 
\end{step}

\begin{proof}(of Step \ref{s4-conic}) 
Assume $p=2$ and $\alpha=\beta=\gamma=0$. 
It is easy to see that $a \neq 0, b \neq 0, c \neq 0$. 
It is enough to show that $[k(a'^{1/2}, b'^{1/2}):k]=4$ for $a':=a/c$ and $b':=b/c$. 
We have $[k(a'^{1/2}, b'^{1/2}):k] \in \{1, 2, 4\}$. 
Note that 
\[
X \simeq \Proj\,k[x, y, z]/(a'x^2+b'y^2+z^2). 
\]
If $[k(a'^{1/2}, b'^{1/2}):k]=1$, i.e. $a'^{1/2} \in k$ and $b'^{1/2} \in k$, 
then we have that 
\[
a'x^2+b'y^2+z^2 = (a'^{1/2}x+b'^{1/2} y+z)^2. 
\]
In this case, $X$ is not reduced, which contradicts the assumption that $X$ is regular. 
Hence we obtain $[k(a'^{1/2}, b'^{1/2}):k]\neq 1$. 

Suppose that $[k(a'^{1/2}, b'^{1/2}):k] = 2$. 
Let us derive a contradiction. 
By symmetry, 
we may assume that $a'^{1/2} \not\in k$. This implies that $k(a'^{1/2})=k(a'^{1/2}, b'^{1/2})$. 
Then we have $b'^{1/2} \in k(a'^{1/2}) = k + k a'^{1/2}$, 
hence we can write $b'^{1/2}=r +sa'^{1/2}$ for some $r, s \in k$. 
In particular, $b'=r^2+s^2a'$. 
Then we obtain
\[
a'x^2+b'y^2+z^2=a'x^2+(r^2+s^2a')y^2+z^2=a'(x+sy)^2+(z+ry)^2 
=a'x'^2+z'^2
\]
for $x':=x+sy$ and $z':=z+ry$. 
Since $a'^{1/2} = z'/x' \in K(X)$ and $X$ is normal, we have that $a'^{1/2} \in H^0(X, \MO_X)$. 
This contradicts $a'^{1/2} \not\in k$ and $H^0(X, \MO_X)=k$. 
This completes the proof of Step \ref{s4-conic}. 
\end{proof}
Step \ref{s1-conic}, Step \ref{s3-conic}, and Step \ref{s4-conic} complete the proof of 
Lemma \ref{l-conic}. 
\end{proof}

\begin{thm}\label{t-conic}
Let $k$ be a field of characteristic $p>0$. 
Let $X$ be a regular proper curve over $k$ 
such that $H^0(X, \MO_X)=k$ and $H^1(X, \MO_X)=0$. 
Then the following hold. 
\begin{enumerate}
\item $X$ is a conic curve in $\mathbb P^2_k$. 
\item $\deg_k K_X=-2$. 
\item $X$ has a $k$-rational point if and only if $X \simeq \mathbb P^1_k$. 
\item If $X$ is geometrically reduced over $k$, then $X$ is smooth over $k$. 
\item If $p \geq 3$, then $X$ is smooth over $k$. 
\item Assume that $p=2$ and $X$ is not smooth over $k$. 
Then a $k$-isomorphism 
\[
X \simeq \Proj\,k[x, y, z]/(sx^2+ty^2+z^2), 
\]
holds for some $s, t \in k$ such that $[k^2(s, t):k^2]=4$. 
In particular, the properties in Proposition \ref{p-Fermat} hold. 
\end{enumerate}
\end{thm}

\begin{proof}
By \cite[Lemma 10.6(3)]{Kol13}, (1) holds. 
Then (1) implies (2). 
Let us show (3). 
Assume that $X$ has a $k$-rational point. 
Then it follows from Corollary \ref{c-rat-red} that 
$X$ is geometrically reduced over $k$. 
Then Lemma \ref{l-conic} implies that $X$ is smooth over $k$. 
Hence, $X$ is a Severi--Brauer curve, 
i.e. $X \times_k \overline k \simeq \mathbb P^1_{\overline k}$. 
As $X$ has a $k$-rational point, 
Ch\^{a}telet's theorem (cf. \cite[Theorem 5.1.3]{GS06}) 
implies that $X \simeq \mathbb P^1_k$. 
Therefore, (3) holds. 

The assertion (4) follows from (1) and the fact that geometrically integral conic is smooth. 
The remaining assertions (5) and (6) follow from (1) and Lemma \ref{l-conic}. 
\end{proof}

\subsection{Curves of genus one}\label{ss4-ex}

In this subsection, we study non-smooth regular curves $X$ of genus one. 
We first restrict possibilities for invariants 
$\gamma(X/k), \ell_F(X/k)$, and $m(X/k)$ (Proposition \ref{p-g-1-inv}). 
We then exhibit several examples 
(Example \ref{e-quasi-ell}, Example \ref{e-Fermat-cubic}, Example \ref{e-Fermat-22}).

\begin{prop}\label{p-g-1-inv}
Let $k$ be a field of characteristic $p>0$. 
Let $X$ be a regular projective curve over $k$ with $H^0(X, \MO_X)=k$ and 
$\dim_k H^1(X, \MO_X)=1$. 
Assume that $X$ is not smooth over $k$. 
Take $\kappa \in \{k^{1/p^{\infty}}, \overline k\}$ and 
set $Y:=(X \times_k \kappa)_{\red}^N$. 
Then the following hold. 
\begin{enumerate}
\item There exist nonzero effective Cartier divisors $D_1, ..., D_{\gamma(X/k)}$ on $Y$ such that 
\[
K_Y+(p-1)\sum_{i=1}^{\gamma(X/k)} D_i \sim 0. 
\]
\item 
It holds that $p \leq 3$. 
\item $-K_Y$ is ample. 
\item If $p=3$, then $\gamma(X/k)=\ell_F(X/k)=1$, and 
$m_F(X/k) \in \{0, 1\}$.% and $\epsilon(X/k) \leq m_F(X/k)(\log_3[k:k^3]-1)$. 
\item 
Assume $p=2$. 
Then the triple $(\gamma(X/k), \ell_F(X/k), m_F(X/k))$ 
satisfies one of the possibilities in the following table. 
\begin{table}[H]
\caption{$p=2$ case}
\begin{tabular}{|c|c|c|} \hline
$\gamma(X/k)$ & $\ell_F(X/k)$ & $m_F(X/k)$ \\ \hline \hline
$1$ & $1$ & $0$  \\ \cline{2-3}
    & $1$ & $1$  \\ 
\hline
$2$ & $1$ & $1$  \\ \cline{2-3}
    & $2$ & $1$  \\ \cline{3-3}
    &     & $2$  \\ 
\hline
\end{tabular}
\end{table}
%\noindent
%Furthermore, the inequality 
%\[
%\epsilon(X/k) \leq m_F(X/k)(\log_2[k:k^2]-1)
%\] 
%holds. 
\end{enumerate}
\end{prop}

\begin{proof}
By Subsection \ref{ss-rel-bet-inv}(3), we have 
\begin{equation}\label{e1-g-1-inv}
m_F(X/k) \leq \ell_F(X/k) \leq \gamma(X/k). 
\end{equation}
%and 
%\begin{equation}\label{e2-g-1-inv}
%\epsilon(X/k) \leq m_F(X/k)(\log_p[k:k^p]-1). 
%\end{equation}
Since $X$ is not geometrically normal, 
it follows from Proposition \ref{p-char-gnor} that 
\begin{equation}\label{e2-g-1-inv}
1 \leq \ell_F(X/k) \leq \gamma(X/k). 
\end{equation}

The assertion (1) follows from $K_X \sim 0$ and Theorem \ref{t-cano-gamma}. 
By $\deg \left(\sum_{i=1}^{\gamma(X/k)} D_i\right) \geq 1$, 
(1) implies that $-K_Y$ is ample and $\deg_{\kappa} K_Y=-2$. 
Thus (2) and (3) hold (note that (2) can be proven by \cite[Theorem 1.1]{PW}, 
although they do not write it explicitly). 
By (1), we have 
\begin{equation}\label{e3-g-1-inv}
(p-1)\sum_{i=1}^{\gamma(X/k)} \deg_{\kappa}D_i =2. 
\end{equation}
Thus, if $p=3$, then $\gamma(X/k)=1$. 
Therefore, (\ref{e1-g-1-inv}) and (\ref{e2-g-1-inv}) imply (4). 

Let us show (5). 
Assume $p=2$. 
Then (\ref{e3-g-1-inv}) implies that $\gamma(X/k) \leq 2$. 
Moreover, if $\gamma(X/k)=1$, then we obtain $\ell_F(X/k)=1$ 
and $m_F(X/k) \in \{0, 1\}$ by (\ref{e1-g-1-inv}) and (\ref{e2-g-1-inv}). 
Thus, the assertion (5) holds for the case when $\gamma(X/k)=1$. 
Therefore, we may assume that $\gamma(X/k)=2$. 
Again by (\ref{e1-g-1-inv}) and (\ref{e2-g-1-inv}), 
it is enough to show that $m_F(X/k) \neq 0$. 
Suppose that $m_F(X/k)=0$, i.e. $X$ is geometrically reduced over $k$ 
(Proposition \ref{p-char-gred}). 
Let us derive a contradiction. 
Take an arbitrary purely inseparable field extension $k \subset k'$ 
such that $X \times_k k'$ is not normal. 
Then $X \times_k k'$ is an integral scheme. 
For the normalisation $X' \to X \times_k k'$ of $X \times_k k'$ 
and the algebraic closure $k''$ of $k'$ in $K(X')$, 
it holds that $k'=k''$ (Proposition \ref{p-epsilon-compute}). 
In particular, $X'$ is geometrically integral over $k'$ (Lemma \ref{l-geom-red-birat}). 
Since $-K_{X'}$ is ample (Theorem \ref{t-T-PW}), 
$X'$ is smooth over $k'$ (Theorem \ref{t-conic}(4)). 
By Definition \ref{d-gamma}, it holds that $\gamma(X/k) \leq 1$, which is absurd. 
Thus, (5) holds. 
\end{proof}

\begin{rem}
We use the same notation as in Proposition \ref{p-g-1-inv}. 
If $p=3$ or $(p, \gamma(X/k))=(2, 2)$, 
then the above argument shows that $Y$ has a $\kappa$-rational point. 
In particular, we have $Y \simeq \mathbb P^1_{\kappa}$ for these cases (Theorem \ref{t-conic}(3)).
\end{rem}

\begin{cor}\label{c-g-1-inv}
Let $k$ be a field of characteristic $p>0$. 
Let $X$ be a projective regular curve over $k$ with $H^0(X, \MO_X)=k$ and 
$\dim_k H^1(X, \MO_X)=1$. 
If $X$ is geometrically reduced over $k$ and not smooth over $k$, 
then $\epsilon(X/k)=m_F(X/k)=0$ and $\gamma(X/k)=\ell_F(X/k)=1$. 
\end{cor}

\begin{proof}
Since $X$ is geometrically reduced over $k$, 
we obtain $\epsilon(X/k)=m_F(X/k)=0$ by Proposition \ref{p-char-gred}. 
Then Proposition \ref{p-g-1-inv} implies $\gamma(X/k)=\ell_F(X/k)=1$. 
\end{proof}

\begin{ex}\label{e-quasi-ell}
Let $\mathbb F$ be an algebraically closed field of characteristic $p \in \{2, 3\}$. 
Let $\pi:S \to B$ be a quasi-elliptic fibration 
(cf. \cite{BM76}, \cite[Section 7]{Bad01}), 
i.e. $\pi$ is a projective $\mathbb F$-morphism 
from a smooth $\mathbb F$-surface $S$ to a smooth $\mathbb F$-curve 
such that $\pi_*\MO_S=\MO_B$ and 
the generic fibre $X$ of $\pi$ is a non-smooth regular projective curve over $k:=K(B)$ 
of genus one. 
In this case, $X$ is geometrically integral over $k$ (cf. \cite[Corollary 7.3]{Bad01}). 
Hence, it follows from Corollary \ref{c-g-1-inv} 
that $\epsilon(X/k)=m_F(X/k)=0$ and $\gamma(X/k)=\ell_F(X/k)=1$. 
\end{ex}

\begin{ex}\label{e-Fermat-cubic}
Let $k$ be a field of characteristic three. 
Let $s, t \in k$ be elements such that $[k^3(s, t):k^3]=9$. 
Then Proposition \ref{p-Fermat} implies that 
\[
X:=\Proj\,k[x, y, z]/(sx^3+ty^3+z^3)
\]
is a regular projective curve $X$ in $\mathbb P^2_k$ of genus one  such that 
$H^0(X, \MO_X)=k$ and  
\[
\ell_F(X/k)=m_F(X/k)=\epsilon(X/k)=\gamma(X/k)=1. 
\]
\end{ex}

\begin{ex}\label{e-Fermat-22}
Let $k$ be a field of characteristic two. 
Let 
\[
S:=\{s_0, s_1, s_2, s_3, t_0, t_1, t_2, t_3\} \subset k
\]
be a subset such that $[k^2(S):k^2]=2^{|S|}=256$. 
Then Lemma \ref{l-pFCI1} and Theorem \ref{t-pFCI-main} imply that 
\[
X:=\Proj\,\frac{k[x_0, x_1, x_2, x_3]}{(s_0x_0^2+s_1x_1^2+s_2x_2^2+s_3x_3^2, t_0x_0^2+t_1x_1^2+t_2x_2^2+t_3x_3^2)}
\]
is a regular projective curve $X$ in $\mathbb P^3_k$ such that 
$H^0(X, \MO_X)=k$,
\[
\ell_F(X/k)=m_F(X/k)=1,\quad\text{and}\quad \quad \epsilon(X/k)=\gamma(X/k)=2. 
\]
\end{ex}

\section{Genus change formula}

The purpose of this section is to establish a genus change formula for regular curves (Theorem \ref{t-genus-change}). 
We then deduce some criteria for geometric reducedness of curves (Corollary \ref{c-genus-change}). Recall that 
if $X$ is a regular projective curve over a field $k$, then 
the genus $g_X$ of $X$ is defined by 
\[
g_X:= \frac{\dim_k H^1(X, \MO_X)}{\dim_k H^0(X, \MO_X)}. 
\]
In particular, if $H^0(X, \MO_X)=k$, then $g_X=\dim_k H^1(X, \MO_X)$. 
%We start with definition of genus. 

\begin{thm}\label{t-genus-change}
Let $k$ be a field of characteristic $p>0$ 
and let $X$ be a regular projective curve over $k$ with $H^0(X, \MO_X)=k$. 
Let $(\kappa, Y)$ be one of 
$(k^{1/p^{\infty}}, (X \times_k k^{1/p^{\infty}})_{\red}^N)$ and 
$(\overline k, (X \times_k \overline k)_{\red}^N)$. 
Assume that $X$ is not smooth over $k$. 
Then the following hold. 
\begin{enumerate}
\item 
There exist nonzero effective Cartier divisors $D_1, ..., D_{\gamma(X/k)}$ on $Y$ 
such that 
the linear equivalence 
\begin{equation}\label{e1-genus-change}
K_Y+ (p-1)\sum_{i=1}^{\gamma(X/k)} D_i \sim f^*K_X
\end{equation}
holds. 
\item 
For nonzero effective Cartier divisors $D_1, ..., D_{\gamma(X/k)}$ 
satisfying (\ref{e1-genus-change}), 
the equation 
\[
2g_Y-2+(p-1)\sum_{i=1}^{\gamma(X/k)}  \deg_{\kappa} D_i = \frac{2g_X-2}{p^{\epsilon(X/k)}}
\]
holds. 
\end{enumerate}
\end{thm}

\begin{proof}
For the separable closure $k'$ of $k$ and the base change $X':=X \times_k k'$, 
it holds that $\gamma(X'/k') \geq \gamma(X/k)$ (Proposition \ref{p-gamma-sep-bc}) 
and $\epsilon(X'/k')=\epsilon(X/k)$ (Proposition \ref{p-epsilon-sep-bc}). 
Therefore, it is enough to treat the case when 
$(\kappa, Y) = (k^{1/p^{\infty}}, (X \times_k k^{1/p^{\infty}})_{\red}^N)$. 
Then the assertion (1) follows from Theorem \ref{t-cano-gamma}.

Let us show (2). 
Taking the degree $\deg_{k^{1/p^{\infty}}}(-)$ of the both hand sides of (\ref{e1-genus-change}), 
we obtain 
\[
2g_Y-2 + (p-1)\sum_{i=1}^{\gamma(X/k)} \deg_{k^{1/p^{\infty}}} D_i = \deg_{k^{1/p^{\infty}}} (f^*K_X). 
\]
We have the following factorisation of $f$ consisting of the induced morphisms: 
\[
f:Y = (X \times_k k^{1/p^{\infty}})_{\red}^N \xrightarrow{f_1} 
(X \times_k k^{1/p^{\infty}})_{\red} 
\xrightarrow{f_2} X \times_k k^{1/p^{\infty}} \xrightarrow{f_3} X. 
\]
Using this factorisation, 
the degree $\deg_{k^{1/p^{\infty}}} (f^*K_X) = \deg_{k^{1/p^{\infty}}} (f^*\omega_X)$ can be computed as follows: 
\begin{eqnarray*}
\deg_{k^{1/p^{\infty}}} (f^*\omega_X) 
&=& \deg_{k^{1/p^{\infty}}} (f_1^*f_2^*f_3^*\omega_X)\\
&=& \deg_{k^{1/p^{\infty}}} (f_2^*f_3^*\omega_X)\\
&=& \frac{1}{p^{\epsilon(X/k)}}\deg_{k^{1/p^{\infty}}} (f_3^*\omega_X)\\
&=& \frac{1}{p^{\epsilon(X/k)}}\deg_{k} (\omega_X),\\
\end{eqnarray*}
where the third equality follows from 
Definition \ref{d-epsilon} and \cite[page 299]{Kle66}, 
and the last equality holds by the flat base change theorem. 
To summarise, (2) holds. 
\end{proof}

\begin{cor}\label{c-genus-change}
Let $k$ be a field of characteristic $p>0$ 
and let $X$ be a regular projective curve over $k$ with $H^0(X, \MO_X)=k$. 
Then the following hold. 
\begin{enumerate}
\item 
If $2g_X-2$ is not divisible by $p$, then $X$ is geometrically integral over $k$. 
\item 
If $X$ is not geometrically integral over $k$, then the inequality
\[
g_X  \geq \frac{(p-1)(p-2)}{2}.
\]
holds. 
\end{enumerate}
\end{cor}

\begin{proof}
Let us show (1). 
Assume that $2g_X-2$ is not divisible by $p$. 
Theorem \ref{t-genus-change}(2) implies that $\epsilon(X/k)=0$. 
It follows from Remark \ref{r-def-epsilon} that $X$ is geometrically reduced over $k$. 
Then (1) holds, since $X$ is geometrically irreducible (Proposition \ref{p-prelim-girre}). 

Let us show (2). Assume that $X$ is not geometrically integral over $k$. 
Then $X$ is not geometrically reduced over $k$ (Proposition \ref{p-prelim-girre}). 
If $g_X=0$, then it follows from Theorem \ref{t-conic} that $p=2$. 
Then the assertion is clear. 
We may assume that $g_X \geq 1$. 
Again by Remark \ref{r-def-epsilon}, we have $\epsilon(X/k) \geq 1$. 
It follows from Proposition \ref{p-def-gamma} that $\gamma(X/k) \geq 1$. 
Therefore, Theorem \ref{t-genus-change}(2) implies that 
\[
-2 + (p-1) \leq 2g_Y-2+(p-1)\sum_{i=1}^{\gamma(X/k)}  \deg_{\kappa} D_i = \frac{2g_X-2}{p^{\epsilon(X/k)}} 
\leq \frac{2g_X-2}{p}.
\]
Hence, we have $p(p-3) \leq 2g_X-2$, which implies $g_X  \geq \frac{(p-1)(p-2)}{2}$, as desired. 
\end{proof}

\begin{rem}
The equality in Corollary \ref{c-genus-change}(2) is attained 
by the examples in Proposition \ref{p-Fermat}. 
Indeed, the inequalities in the proof of Corollary \ref{c-genus-change}(2) 
become equalities when $g_Y=0$ and $\gamma(X/k)=\epsilon(X/k)=1$. 
\end{rem}

\section{Boundedness of regular curves}\label{s-bdd}

In this section, 
we discuss boundedness of regular curves over imperfect fields. 
It is not difficult to show the boundedness of regular curves of genus $g \neq 1$ 
(Proposition \ref{p-g0-curve-bdd}, Theorem \ref{t-g2-curve-bdd}). 
For the case of genus one, 
non-smooth regular curves appear only in characteristic two and three (Proposition \ref{p-g-1-inv}(2)). 
We prove that such curves are bounded when we fix a base field $k$ 
such that $[k:k^p]<\infty$ (Theorem \ref{t-g1-bdd-F-fin}). 
These results will be established in Subsection \ref{ss2-bdd}. 
To this end, we first recall a criterion for very ampleness in Subsection \ref{ss1-bdd}. 
In Subsection \ref{ss3-bdd}, 
we also prove that if $X$ is a regular curve of genus one 
which is geometrically integral and not smooth, 
then $X$ is a cubic curve when $p=3$ and 
$X$ is a complete intersection of two quadric surfaces when $p=2$ (Theorem \ref{t-gen-g-1}).

\subsection{Very ampleness}\label{ss1-bdd}

The purpose of this subsection is to establish 
a criterion for very ampleness (Theorem \ref{t-va-crit}). 
The strategy is to use Mumford's regularity, 
which is a standard argument at least for varieties over algebraically closed fields. 
We only show how to fill gaps between references and what we need. 
For details, we refer to \cite[Section 5.2]{FGAex} or \cite[Section 1.8.B]{Laz04}.

\begin{dfn}\label{d-m-reg}
Let $k$ be a field and let $X$ be a projective scheme over $k$. 
Let $H$ be an ample globally generated invertible sheaf on $X$. 
Fix $m \in \Z$. 
A coherent sheaf $F$ on $X$ is $m$-{\em regular} with respect to $A$ if 
it holds that 
\[
H^i(X, F \otimes_{\MO_X} A^{m-i})=0
\]
for any $i>0$. 
If $D$ is a Cartier divisor on $X$, 
then $D$ is $m$-{\em regular} with respect to $A$ if so is the associated invertible sheaf $\MO_X(D)$. 
\end{dfn}

\begin{lem}\label{l-Mumford-gen}
Let $k$ be a field and let $X$ be a projective scheme over $k$. 
Let $H$ be an ample globally generated invertible sheaf on $X$. 
For $m \in \Z$, 
let $F$ be a coherent sheaf $F$ which is $m$-regular 
with respect to $A$. 
Then the induced map 
\[
H^0(X, A) \otimes_k H^0(X, F \otimes A^{m+r}) \to H^0(F \otimes A^{m+r+1})
\]
is surjective for any $r \in \Z_{\geq 0}$. 
\end{lem}

\begin{proof}
Let $f:X \to \mathbb P^n_k$ be the projective morphism induced by $|H|$. 
Note that $f$ is a finite morphism, 
hence we may assume that $X=\mathbb P^n_k$ and $H=\MO_{\mathbb P^n_k}(1)$. 
In this case, the assertion follows from \cite[Lemma 5.1(a)]{FGAex}. 
\end{proof}

\begin{prop}\label{p-m-reg}
Let $k$ be a field and let $X$ be a projective scheme over $k$. 
Let $H$ be an an ample globally generated invertible sheaf on $X$. 
Let $L$ be an invertible sheaf on $X$. 
If $L$ is $0$-regular with respect to $H$, then $L \otimes H$ is very ample over $k$. 
\end{prop}

\begin{proof}
Taking the base change to the algebraic closure of $k$, 
we may assume that $k$ is algebraically closed. 
Fix a closed point $x \in X$. 
It is enough to show that $L \otimes H \otimes \m_x$ is globally generated, 
where $\m_X$ denotes the coherent ideal sheaf corresponding to $x$. 
This follows from the same argument as in \cite[Lemma 3.7]{Wit17}. 
\end{proof}

\begin{thm}\label{t-va-crit}
Let $X$ be a regular projective curve of genus $g \geq 1$ such that $H^0(X, \MO_X)=k$. 
Let $D$ be a Cartier divisor on $X$. 
Then the following holds. 
\begin{enumerate}
\item 
If $\deg_k D \geq g$, then $H^0(X, D) \neq 0$. 
\item 
If $\deg_k D \geq 2g-1$, then $\MO_X(2D)$ is globally generated. 
\item 
If $\deg_k D \geq 2g-1$, then $|5D|$ is very ample over $k$. 
\end{enumerate}
\end{thm}

\begin{proof}
The assertion (1) follows from 
\[
h^0(X, D) \geq \chi(X, D)=\deg_k D + 1 - g \geq 1, 
\]
where the equality follows from the Riemann--Roch theorem. 

Let us show (2).  
Since $2g -1 \geq g$, (1) implies $H^0(X, D) \neq 0$. 
Hence, we may assume that $D$ is a nonzero effective Cartier divisor. 
Consider the induced exact sequence:  
\[
0 \to \MO_X(-D) \to \MO_X \to \MO_D \to 0. 
\]
We have 
\begin{equation}\label{e1-va-crit}
H^1(X, \MO_X(D))=0
\end{equation}
by Serre duality and $\deg_k D \geq 2g-1$. 
Thus, the induced map 
\[
H^0(X, \MO_X(2D)) \to H^0(D, \MO_X(2D)|_D)
\]
is surjective. 
Therefore, $\MO_X(2D)$ is globally generated. 
Thus, (2) holds. 

Let us show (3). 
By (2), $\MO_X(2D)$ is globally generated. 
By Proposition \ref{p-m-reg}, 
it is enough to show that 
$3D$ is $0$-regular with respect to $\MO_X(2D)$, which follows from 
\[
H^1(X, \MO_X(3D) \otimes_{\MO_X} \MO_X(-2D)) \simeq 
H^1(X, \MO_X(D)) =0, 
\]
where the equality holds by (\ref{e1-va-crit}). 
Thus (3) holds. 
\end{proof}

\subsection{Boundedness}\label{ss2-bdd}

In this subsection, 
we establish several results on boundedness of regular curves. 
Although the case of genus $g \neq 1$ might be known for experts, 
we include the proof for the sake of completeness. 
An essentially new result is Theorem \ref{t-g1-bdd-F-fin}, 
which establishes the boundedness of 
non-smooth regular curves of genus one over a fixed field $k$ with $[k:k^p]<\infty$.

\subsubsection{Genus zero} 

It is well known that regular curves of genus zero is a conic curve. 

\begin{prop}\label{p-g0-curve-bdd}
Let $k$ be a field. 
Let $X$ be a regular projective curve over $k$ such that $H^0(X, \MO_X)=k$ and $-K_X$ is ample. 
Then $|-K_X|$ is very ample and it induces a closed immersion $\Phi_{|-K_X|}:X \hookrightarrow \mathbb P^2_k$ whose image is a conic curve.  
\end{prop}

\begin{proof}
See \cite[Lemma 10.6]{Kol13}. 
\end{proof}

\subsubsection{Higher genus} 

%For the case when $g \geq 2$, it holds that $|10K_X|$ is very ample. 

\begin{thm}\label{t-g2-curve-bdd}
Let $X$ be a regular projective curve of genus $g_X \geq 2$ 
such that $H^0(X, \MO_X)=k$. 
Then the following hold. 
\begin{enumerate}
\item $\MO_X(4K_X)$ is globally generated. 
\item $|10K_X|$ is very ample over $k$ and it induces a closed immersion  to $\mathbb P_k^{19g-20}$: 
\[
\Phi_{|10K_X|}:X \hookrightarrow \mathbb P_k^{19g-20}.
\] 
\end{enumerate}
\end{thm}

\begin{proof}
We can apply Theorem \ref{t-va-crit} for $D:=2K_X$, since 
$\deg\,D = \deg\,(2K_X) = 4g-4\geq 2g-1$. 
Hence, $\MO_X(4K_X)$ is globally generated and 
$|10K_X|$ is very ample over $k$. 
Then the remaining assertion follows from 
\[
h^0(X, 10K_X) = \chi(X, 10K_X) = \deg_k (10K_X)+1-g=19(g-1), 
\]
where the first and second equalities hold by Serre duality and the Riemann--Roch formula, respectively. 
\end{proof}

\subsubsection{Genus one}

We now focus on the case of genus one. 
The main result is Theorem \ref{t-g1-bdd-F-fin}, 
i.e. we want to show that 
if the base field $k$ satisfies $[k:k^p] < \infty$, 
then the boundedness holds for 
non-smooth regular $k$-curves of genus one. 
To this end, we give a similar boundedness result in terms of the invariant $\epsilon(X/k)$ 
(Theorem \ref{t-g1-bdd-div}). 
We start with the following result.

\begin{lem}\label{l-Frob-symm}
Let $\zeta:Y \to Z$ be a finite bijective morphism of $\F_p$-schemes, 
where $Y$ is reduced. 
Let $\ell$ be a positive integer and assume that 
the $\ell$-th iterated absolute Frobenius morphism $F^{\ell}_Z$ factors through $\zeta$: 
\[
F^{\ell}_Z:Z \xrightarrow{\xi} Y \xrightarrow{\zeta} Z.
\]
Then we obtain another factorisation $F^{\ell}_Y:Y \xrightarrow{\zeta} Z \xrightarrow{\xi} Y$, 
i.e. $\xi \circ \zeta=F^{\ell}_Y$.  
\end{lem}

\begin{proof}
The problem is local, hence we may assume that $Z=\Spec\,A$ and $Y=\Spec\,B$. 
Our assumption induces the following factorisation of ring homomorphisms: 
\[
F_A^{\ell}:A \xrightarrow{\varphi} B \xrightarrow{\psi} A,
\]
i.e. the equation $\psi \circ \varphi(a)=a^{p^{\ell}}$ holds for any $a \in A$. 
Fix $b \in B$. 
It is enough to prove that $\varphi \circ \psi(b) = b^{p^{\ell}}$. 
We have 
\[
\psi(\varphi \circ \psi(b))= 
\psi \circ \varphi \circ \psi (b) =  \psi \circ \varphi(\psi (b))=\psi(b)^{p^{\ell}} = \psi(b^{p^{\ell}}).
\]
Hence, we obtain $\varphi \circ \psi(b) - b^{p^{\ell}} \in \Ker (\psi)$. 
Therefore, it is sufficient to show that $\psi:B \to A$ is injective. 
This follows from the fact that the corresponding morphism $\xi:Z=\Spec\,A \to \Spec\,B=Y$ is 
an affine surjective morphism to a reduced scheme $Y$. 
\end{proof}

\begin{thm}\label{t-g1-bdd-div}
Let $k$ be a field of characteristic $p>0$. 
Let $X$ be a regular projective curve over $k$ with $H^0(X, \MO_X)=k$ and 
$g_X=1$. 
Assume that $X$ is not smooth over $k$. 
Then the following hold. 
\begin{enumerate}
\item If $p=2$, then there exists an ample Cartier divisor $D$ on $X$ 
such that $\deg_k D \leq 2^{\epsilon(X/k)+2}$. 
\item If $p=2$ and $X$ is geometrically reduced over $k$, 
then there exists a Cartier divisor $D$ on $X$ such that $\deg_k D =4$. 
\item If $p=3$, then there exists an ample Cartier divisor $D$ on $X$ 
such that $\deg_k D \leq 3^{\epsilon(X/k)+1}$. 
\item If $p=3$ and $X$ is geometrically reduced over $k$, 
then there exists a Cartier divisor $D$ on $X$ such that $\deg_k D =3$. 
\end{enumerate}
\end{thm}

\begin{proof}
We first reduce the problem to the case when $k$ is a finitely generated field over $\F_p$. 
There exists a subfield $k_0 \subset k$ and a regular projective curve $X_0$ over $k_0$ 
such that $k_0$ is finitely generated over $\F_p$ and $X_0 \times_{k_0} k \simeq X$. 
Note that all the assumptions on $(k, X)$ hold also for $(k_0, X_0)$. 
Furthermore, we have that $\epsilon(X_0/k_0)=\epsilon(X/k)$ 
(Proposition \ref{p-epsilon-des-bc}), 
hence the problem is reduced to the case when $k$ is finitely generated over $\F_p$.

Theorem \ref{t-schroer} enables us to find an intermediate field $k \subset l \subset k^{1/p}$ 
such that $X \times_k l$ is an integral scheme and the algebraic closure of $l$ in $K(X \times_k l)$ is equal to $k^{1/p}$. 
For the normalisation $X'$ of $X \times_k l$, 
we have the following commutative diagram: 
\[
\begin{CD}
X @<\varphi <<  X' @<\psi << X\\
@VVV @VVV @VVV\\
\Spec\,k @<<< \Spec\,k^{1/p} @<\simeq <<\Spec\,k,
\end{CD}
\]
where the horizontal composite arrows are the absolute Frobenius morphisms. 
It follows from Lemma \ref{l-Frob-symm} that also $\psi \circ \varphi:X' \to X'$ is the absolute Frobenius morphism. 
By Theorem \ref{t-T-PW}, 
there exists a nonzero effective Cartier divisor $C$ on $X'$ such that 
\[
K_{X'}+(p-1) C \sim \varphi^*K_X \sim 0.  
\]
If $p=2$ (resp. $p=3$), then $A:=-K_{X'}$ (resp. $A:=C$) 
is a Cartier divisor on $X'$ such that  $\deg_{k^{1/p}} A=2$ (resp. $\deg_{k^{1/p}} A=1$). 
For $D:=\psi^*A$, let us prove that (1)--(4) hold.

We have a commutative diagram consisting of the induced morphisms 
\[
\begin{CD}
X' @<\psi << X @<<< X \times_k l @<<< X'\\
@VVV @VVV @VVV @VVV\\
\Spec\,k^{1/p} @<\simeq <<\Spec\,k @<<< \Spec\,l @<<< \Spec\,k^{1/p}
\end{CD}
\]
where the composite horizontal morphisms are the absolute Frobenius morphisms. 
Since 
\[
\deg_k D=(\deg \psi) \cdot (\deg_{k^{1/p}} A), 
\]
let us compute $\deg \psi$. 
It follows that 
\begin{itemize}
\item $[k^{1/p}:k]=[k^{1/p}:l][l:k]$ and
\item $[K(X')^{1/p}:K(X')]=\deg \psi \times [l:k]$, 
\end{itemize}
where the latter equation holds, 
because the isomorphism $K(X) \otimes_k l \simeq K(X \times_k l)$ implies 
$[K(X \times_k l):K(X)]=[l:k]$. 
Since $[k^{1/p}:k]=p^{{\rm tr.deg}_{\F_p} k}$ and 
$[K(X')^{1/p}:K(X')] =p^{{\rm tr.deg}_{\F_p} K(X')}$, 
it holds that $[K(X')^{1/p}:K(X')]=p \cdot [k^{1/p}:k]$. 
We then obtain 
\[
\deg \psi = \frac{[K(X')^{1/p}:K(X')]}{[l:k]}=\frac{p \cdot [k^{1/p}:k]}{[l:k]} 
=p \cdot [k^{1/p}:l].
\]
Therefore, it holds that 
\[
\deg_k D = (\deg \psi) \cdot (\deg_{k^{1/p}} A) = p \cdot [k^{1/p}:l] \cdot (\deg_{k^{1/p}} A). 
\]
It follows from Proposition \ref{p-epsilon-compute}(1)  
that $[k^{1/p}:l] \leq p^{\epsilon(X/k)}$. 
Hence, we obtain 
\[
\deg_k D = p \cdot [k^{1/p}:l] \cdot (\deg_{k^{1/p}} A) \leq p^{\epsilon(X/k)+1} \cdot (\deg_{k^{1/p}} A).
\]
Thus, (1) and (3) hold. 

To show (2) and (4), we assume that $X$ is geometrically reduced. 
Then Remark \ref{r-def-epsilon} implies that $\epsilon(X/k)=0$. 
By $[k^{1/p}:l] \leq p^{\epsilon(X/k)}$, we obtain $[k^{1/p}:l]=1$. 
Therefore, it holds that  
\[
\deg_k D = p \cdot [k^{1/p}:l] \cdot (\deg_{k^{1/p}} A) =p \cdot (\deg_{k^{1/p}} A).
\]
Thus, (2) and (4) hold. 
\end{proof}

\begin{thm}\label{t-g1-bdd-F-fin}
Let $k$ be a field of characteristic $p>0$ such that $[k:k^p]<\infty$. 
Let $X$ be a regular projective curve over $k$ with $H^0(X, \MO_X)=k$ and 
$g_X=1$. 
Assume that $X$ is not smooth over $k$. 
Then the following hold. 
\begin{enumerate}
\item 
If $p=2$, then there exists a very ample Cartier divisor $H$ on $X$ such that 
$\deg_k H \leq 5([k:k^2])^2$. 
%Its complete linear system induces a closed immersion to $\mathbb P_k^N$ for $N:=3[k:k^2]$ 
\item 
If $p=3$, then 
there exists a very ample Cartier divisor $H$ on $X$ such that 
$\deg_k H \leq \frac{5}{3}[k:k^3]$. 
\end{enumerate}
\end{thm}

\begin{proof}
Let $D$ be an ample Cartier divisor as in (1) and (3) of Theorem \ref{t-g1-bdd-div}. 
Set $H:=5D$. 
Then it follows from Theorem \ref{t-va-crit} that $H$ is very ample. 
It is enough to prove the inequalities as in the statement. 

Since the proofs are very similar, 
let us treat only the case when $p=2$. We have that 
\[
\epsilon(X/k) \leq m_F(X/k)(\log_2 [k:k^2]-1) \leq 2\log_2 [k:k^2]-2,
\]
where the first inequality follows from Theorem \ref{t-epsilon-m} and 
the second one holds by $m_F(X/k) \leq 2$ (Proposition \ref{p-g-1-inv}(5)). 
Therefore, it holds that 
\[
\deg_k D \leq 2^{\epsilon(X/k)+2} \leq 2^{2\log_2[k:k^2]} =([k:k^2])^2.
\]
Hence, we obtain $\deg_k H = 5 \deg_k D \leq 5([k:k^2])^2$, as desired. 
%We now treat the case when $p=3$. 
%Since $m_F(X/k) \leq 1$, we have that 
%\[
%\epsilon(X/k) \leq m_F(X/k)(\log_3 [k:k^3]-1) \leq \log_3 [k:k^3]-1.
%\]
%Therefore, 
%\[
%\deg_k D \leq 3^{\epsilon(X/k)+1} \leq 3^{\log_3 [k:k^3]-1} =\frac{1}{3}[k:k^3].
%\]
%Hence, $3D$ is very ample and $\deg_k H \leq \frac{5}{3}[k:k^3]$. 
\end{proof}

\subsection{Non-smooth geometrically reduced curves of genus one}\label{ss3-bdd}

The purpose of this subsection is to determine 
non-smooth geometrically integral curves of genus one 
(Theorem \ref{t-gen-g-1}). 
To this end, we shall prove 
Lemma \ref{l-bs-pt-g1}, Proposition \ref{p-gen-g-1}, and Corollary \ref{c-gen-g-1}. 
The proofs of Lemma \ref{l-bs-pt-g1} and Proposition \ref{p-gen-g-1} 
are based on discussion of the author with Fabio Bernasconi.

\begin{lem} \label{l-bs-pt-g1}
Let $k$ be a field. 
Let $X$ be a geometrically integral Gorenstein projective curve over $k$ 
such that $\dim_k H^1(X, \MO_X)=1$. 
Let $D$ be a Cartier divisor on $X$.
\begin{enumerate}
\item If $\deg_k (D) \geq 1$, then $H^1(X, \MO_X(D))=0$ and $\dim_k H^0(X, \MO_X(D)) 
= \deg_k D$. 
\item If $\deg_k (D) =1$, then $\Bs(L)=\{P\}$ for some $k$-rational point $P$, 
where $\Bs(L):= \bigcap_{s \in H^0(X, \MO_X(D))} \{Q \in X\,|\,s|_Q=0\}$ 
is the base locus of $|D|$. 
Furthermore, $X$ is smooth around $P$. 
\item If $\deg_k(D) \geq 2$, then $\MO_X(D)$ is globally generated.
\end{enumerate}
\end{lem}

\begin{proof}
Let us show (1). 
It follows from Serre duality that $h^1(X, \MO_X(D))=h^0(X, \MO_X(-D))=0$.
By the Riemann--Roch theorem, we obtain 
\begin{equation}\label{e1-bs-pt-g1}
\dim_k H^0(X, \MO_X(D)) =  \deg_k D. 
\end{equation}
Thus, (1) holds. 

Let us show (2). 
The equation (\ref{e1-bs-pt-g1}) implies that $\Bs(L)=\{P\}$ for some $k$-rational point $P$. 
Let us prove that $X$ is smooth around $P$. 
We may assume that $k$ is algebraically closed. 
Since a closed point $P$ is an effective Cartier divisor, 
the maximal ideal $\m_P$ of the stalk $\MO_{X, P}$ is generated by a single element. 
Hence, $X$ is regular around $P$, as desired. 
Thus, (2) holds. 

Let us show (3). 
We may assume that $k$ is algebraically closed. 
Set $d:=\deg_k D$. 
For a smooth closed point $Q$ of $C$, 
it follows from (2) that $D-(d-1)Q \sim P$ for some smooth closed point $P$ on $X$ 
(note that $P$ is possibly equal to $Q$). 
Therefore, the base locus $\Bs (D)$ of $D$ is contained in $\{P, Q\}$. 
For any smooth closed point $R$ of $X$, it holds that 
\[
H^0(X, \MO_X(D)) \to H^0(R, \MO_X(D)|_R) \to H^1(X, \MO_X(D-R))=0,
\]
where the last equality follows from Serre duality and $\deg_k D \geq 2$. 
Hence, $\Bs (D)$ does not contain $R$.  
Applying this for $R \in \{P, Q\}$, we have that $\MO_X(D)$ is globally generated. 
Hence (3) holds. 
\end{proof}

\begin{prop} \label{p-gen-g-1}
Let $k$ be a field. 
Let $X$ be a geometrically integral Gorenstein projective curve over $k$ 
such that $\dim_k H^1(X, \MO_X)=1$. 
Let $D$ be a Cartier divisor on $X$ and 
set $R(X, D):=\bigoplus_{m \geq 0} H^0(X, \MO_X(mD))$ 
to be the graded $k$-algebra. 
Then the following hold. 
\begin{enumerate}
\item If $\deg_k D =1$, then $R(X, D)$ is generated 
by $\bigoplus_{1 \leq j \leq 3} H^0(X, \MO_X(jD) )$ as a $k$-algebra. 
\item If $\deg_k D =2$, then $R(X, D)$ is generated 
by $H^0( X, \MO_X(D) ) \oplus H^0(X, \MO_X(2D))$ as a $k$-algebra. 
\item If $\deg_k D \geq 3$, then $R(X, D)$ is generated by $H^0( X, \MO_X(D) )$ as a $k$-algebra. 
\item If $\deg_k D \geq 3$, then $D$ is very ample over $k$. 
\end{enumerate}
\end{prop}

\begin{proof}
We may assume that $k$ is algebraically closed. 
We first show that 
\begin{enumerate}
\item[(2)'] 
if $\deg_k D \geq 2$, then $R(X, \MO_X(D))$ is generated 
by $H^0( X, \MO_X(D) ) \oplus H^0(X, \MO_X(2D))$ as a $k$-algebra. 
\end{enumerate}
Assume $\deg_k D \geq 2$. 
By Lemma \ref{l-bs-pt-g1}(3), $\MO_X(D)$ is an ample globally generated invertible sheaf. 
We have 
\[
H^1(X, \MO_X(D) \otimes \MO_X(D)^{m-1})=0 
\]
for $m=1$. 
Then $\MO_X(D)$ is $1$-regular with respect to $\MO_X(D)$ 
(Definition \ref{d-m-reg}). 
It follows from Lemma \ref{l-Mumford-gen} that the induced map 
\[
H^0(X, \MO_X(D)) \otimes_k H^0(X, \MO_X((r+2)D)) \to H^0(X, \MO_X((r+3)D))
\]
is surjective for any $r \in \Z_{\geq 0}$. 
Therefore, (2)' holds. 
In particular, also (2) holds. 
We omit the proof of (1), as a similar argument can be applied.

Let us show (3). 
Assume $d:=\deg_k D \geq 3$. 
By (2)', it is enough to show that the induced map 
\[ 
\beta:H^0(X, \MO_X(D)) \otimes H^0(X, \MO_X(D)) \rightarrow H^0(X, \MO_X(2D)) 
\]
is surjective. 
Since $\MO_X(D)$ is globally generated, 
we obtain $D \sim \sum_{i=1}^d P_i$ 
for some closed points $P_1, ..., P_d$ around which $X$ is smooth. 
Consider the Cartier divisor $D':=P_1+P_2$ and the closed subscheme $T$ 
defined by the nonzero effective Cartier divisor $\sum_{i=3}^d P_i$.
We have the following exact sequence
\[ 
0 \to   \MO_X(D') \to \MO_X(D) \to \MO_T \to 0. 
\]
Note that $H^1(X, \MO_X(D'))=0$ by Serre duality. 
Then we obtain a commutative diagram
\begin{equation*}
\small{
\xymatrix{
0 \ar[r] & H^0(D') \otimes H^0(D) \ar[r] \ar[d]^{\alpha} & H^0(D) \otimes H^0(D) \ar[d]^{\beta} \ar[r] & H^0(\mathcal{O}_T) \otimes H^0(D) \ar[r] \ar[d]^{\gamma} & 0 \\
0 \ar[r] & H^0(D'+D) \ar[r] & H^0(2D) \ar[r] & H^0(\MO_X(D)|_T) \ar[r] & 0,
}
}
\end{equation*}
where the horizontal sequences are exact and 
the vertical arrows are multiplicative maps. 
The map $\alpha$ is surjective because $D$ is $0$-regular with respect to $\MO_X(D')$ (Definition \ref{d-m-reg}). 
The map $\gamma$ is surjective because $\mathcal{O}_T$ is $0$-regular with respect to $\MO_X(D)$. 
Therefore, the snake lemma implies that $\beta$ is surjective. 
Hence, (3) holds. 
The assertion (4) follows from (3). 
\end{proof}

\begin{cor} \label{c-gen-g-1}
Let $k$ be a field. 
Let $X$ be a geometrically integral Gorenstein projective curve over $k$ 
such that $\dim_k H^1(X, \MO_X)=1$. 
\begin{enumerate}
\item 
If $D$ is a Cartier divisor on $X$ such that $\deg_k D=3$, 
then 
$|D|$ is very ample and it induces a closed immersion $\Phi_{|D|}:X \hookrightarrow \mathbb P^2_k$ whose image is a cubic curve.  
\item 
If $D$ is a Cartier divisor on $X$ such that $\deg_k D=4$, 
then 
$|D|$ is very ample and it induces a closed immersion $\Phi_{|D|}:X \hookrightarrow \mathbb P^3_k$ whose image is a complete intersection $S_1 \cap S_2$ of two quadric hypersurfaces $S_1$ and $S_2$.  
\end{enumerate}
\end{cor}

\begin{proof}
Let $D$ be a Cartier divisor on $X$ such that $d:=\deg_k D \geq 3$. 
By Lemma \ref{l-bs-pt-g1}(1), we have 
\[
\dim_k H^0(X, \MO_X(D))= \deg_k D. 
\]
Furthermore, it follows from Proposition \ref{p-gen-g-1}(4) that 
$|D|$ is very ample. 
Hence, it induces a closed immersion $\Phi_{|D|}:X \hookrightarrow \mathbb P^{d-1}_k$. 

If $d=\deg_k D=3$, then $\Phi_{|D|}(X)$ is a cubic curve because $\omega_X \simeq \MO_X$. 
Thus (1) holds. 
Let us show (2). 
Assume $d=\deg_k D=4$. 
We obtain 
\[
\dim_k H^0(\mathbb P^3_k, \MO_{\mathbb P^3_k}(2))- \dim_k H^0(X, \MO_X(2D))=10-8=2. 
\]
Thus, (2) holds. 
\end{proof}

\begin{thm}\label{t-gen-g-1}
Let $k$ be a field of characteristic $p>0$. 
Let $X$ be a geometrically integral regular projective curve over $k$ of genus one. 
Assume that $X$ is not smooth over $k$. 
Then the following hold. 
\begin{enumerate}
\item 
If $p=3$, then there exists a closed immersion $j:X \hookrightarrow \mathbb P^2_k$ whose image is a cubic curve. 
\item 
If $p=2$, then there exists a closed immersion $j:X \hookrightarrow \mathbb P^3_k$ whose image is a complete intersection $S_1 \cap S_2$ of two quadric surfaces $S_1$ and $S_2$.  
\item 
Assume $p=2$. 
Then the following are equivalent. 
\begin{enumerate}
\item There exists a closed immersion $j:X \hookrightarrow \mathbb P^2_k$ whose image is a cubic curve. 
\item $X$ has a $k$-rational point. 
\end{enumerate}
\end{enumerate}
\end{thm}

\begin{proof}
If $p=3$ (resp. $p=2$), 
then Theorem \ref{t-g1-bdd-div} enables us to find 
a Cartier divisor $D$ on $X$ such that $\deg_k D=3$ (resp. $\deg_k D=4$). 
Thus, the assertions (1) and (2) follow from Corollary \ref{c-gen-g-1}. 

Let us prove (3). 
Assume $p=2$. 
Consider the following conditions. 
\begin{enumerate}
\item[(a)'] There is a Cartier divisor $A$ on $X$ such that $\deg_k A=3$. 
\item[(b)'] There is a Cartier divisor $B$ on $X$ such that $\deg_k B=1$. 
\end{enumerate}
Since we have a Cartier divisor $D$ on $X$ such that  $\deg_k D=4$, 
the conditions (a)' and (b)' are equivalent. 
By Corollary \ref{c-gen-g-1}(1), (a) is equivalent to (a)'. 
It follows from Lemma \ref{l-bs-pt-g1}(2) that (b) is equivalent to (b)'. 
To summarise, (a) is equivalent to (b). 
\end{proof}

As the following example shows, 
if $p=2$, then there exists a geometrically integral non-smooth regular projective curve 
of genus one which is not $k$-isomorphic to a cubic curve in $\mathbb P^2$.

\begin{ex}
Let $\F$ be an algebraically closed field of characteristic two. 
There exists an Enriques surface $S$ over $\F$ 
with a quasi-elliptic fibration  $\pi:S \to B:=\mathbb P^1$ (cf. \cite[Section 11]{KKM}). 
Then $\pi$ has multiple fibres (cf. \cite[Remark 8.3(c)]{Bad01}). 
Set $k:=K(B)$ and $X:=S \times_B \Spec\,K(B)$. 
Then $X$ is a geometrically integral non-smooth regular projective curve of genus one. 
Since $\pi$ has multiple fibres, $X$ has no $k$-rational points. 
By Theorem \ref{t-gen-g-1}(3), $X$ 
is not $k$-isomorphic to a cubic curve in $\mathbb P^2$. 
%On the other hand, Theorem \ref{t-g1-bdd-div} implies that 
%there exists a Cartier divisor $D$ on $X$ such that $\deg_k D=4$. 
%We now show that $X$ is not $k$-isomorphic to a cubic curve.  
%Suppose that $X$ is $k$-isomorphic to a cubic curve. 
%Let us derive a contradiction. 
%There exists a Cartier divisor $E$ on $X$ such that $\deg_k E=3$. 
%Then $A:=D-E$ is a Cartier divisor with $\deg_k A=1$. 
%Since $H^0(X, \MO_X(A)) \neq 0$, there exists a $k$-rational point on $X$. 
%This is a contradiction. 
\end{ex}

\begin{rem}
In Subsection \ref{ss4-ex}, 
we exhibit several examples of non-smooth regular curves of genus one 
(Example \ref{e-quasi-ell}, Example \ref{e-Fermat-cubic}, Example \ref{e-Fermat-22}). 
However, by Theorem \ref{t-gen-g-1}, it turns out that all the examples satisfy either 
\begin{enumerate}
\item $p=3$ and $X$ is a cubic curve in $\mathbb P^2$, or 
\item $p=2$ and $X$ is a complete intersection of two quadric surfaces in $\mathbb P^3$. 
\end{enumerate}
The author does not know whether 
there exists a non-smooth regular curve of genus one 
which satisfies neither (1) nor (2). 
\end{rem}

\end{document}